\documentclass[a4paper,11pt,BCOR1cm,parskip]{scrartcl}

\usepackage{german}
\usepackage[iso]{umlaute}
\usepackage{theorem}
\usepackage{amsmath}
\usepackage{amssymb}
\usepackage{makeidx}
\usepackage[matrix,arrow]{xy}
\usepackage{epsfig}
\usepackage{boxedminipage}
\usepackage{epic}
\usepackage{longtable}
\usepackage{ifthen}
\usepackage{tabularx}
\usepackage{curves}
\usepackage{parskip}

\newcommand{\R}{\mathrm{I\!R}}
\newcommand{\N}{\mathrm{I\!N}}
\newcommand{\HH}{\mathrm{I\!H}}

\newcommand{\K}{\mathrm{I\!K}}
\newcommand{\PP}{\mathrm{I\!P}}

\newcommand{\Z}{\mathchoice {\hbox{$\sf\textstyle Z\kern-0.4em
Z$}}{\hbox{$\sf\textstyle Z\kern-0.4em Z$}}{\hbox{$\sf\scriptstyle
Z\kern-0.3em Z$}}{\hbox{$\sf\scriptscriptstyle Z\kern-0.2em Z$}}}

\newcommand{\Q}{\mathchoice {\setbox0=\hbox{$\displaystyle\rm
Q$}\hbox{\raise0.15\ht0\hbox to0pt{\kern0.4\wd0\vrule
height0.8\ht0\hss}\box0}}{\setbox0=\hbox{$\textstyle\rm
Q$}\hbox{\raise0.15\ht0\hbox to0pt{\kern0.4\wd0\vrule
height0.8\ht0\hss}\box0}}{\setbox0=\hbox{$\scriptstyle\rm
Q$}\hbox{\raise0.15\ht0\hbox to0pt{\kern0.4\wd0\vrule
height0.7\ht0\hss}\box0}}{\setbox0=\hbox{$\scriptscriptstyle\rm
Q$}\hbox{\raise0.15\ht0\hbox to0pt{\kern0.4\wd0\vrule
height0.7\ht0\hss}\box0}}}

\newcommand{\id}{\mathrm{id}}

\newcommand{\SO}{\mathrm{SO}}
\newcommand{\End}{\mathrm{End}}

\newcommand{\rk}{\mathrm{rk}}

\newcommand{\eps}{\varepsilon}
\newcommand{\vi}{\varphi}
\newcommand{\vkap}{\varkappa}

\newcommand{\qmq}[1]{\quad\mbox{#1}\quad}
\newcommand{\Menge}[2]{\{\,#1\,|\,#2\,\}}

\renewcommand{\bigoplus}{\mathop{\bigcirc
  \raisebox{-0.22em}{\hskip-0.53em\hbox{\vrule height2.08ex width0.04em}
  \raisebox{ 0.48em}{\hskip-0.75em\hbox{\vrule height0.04em width 0.8em}}
  \hskip- 0.2em}}}

\newcommand {\g}[2]{\langle #1,#2\rangle}

\newcommand {\operp}{\mathbin{\mbox{$\ominus\raisebox{2.9pt}
 {\hskip-0.42em\hbox{\vrule height0.7ex width0.02em}\hskip0.42em }$}}}

\newcommand{\Og}{\mathrm{O}}
\newcommand{\Ug}{\mathrm{U}}
\newcommand{\SU}{\mathrm{SU}}
\newcommand{\Sp}{\mathrm{Sp}}

\newcommand{\Spin}{\mathrm{Spin}}

\newcommand{\RE}{\mathop{\mathrm{Re}}\nolimits}
\newcommand{\IM}{\mathop{\mathrm{Im}}\nolimits}

\newcommand{\Eig}{\mathop{\mathrm{Eig}}\nolimits}

\newcommand{\ad}{\mathop{\mathrm{ad}}\nolimits}

\newcommand{\spn}{\mathop{\mathrm{span}}\nolimits}
\newcommand{\tr}{\mathop{\mathrm{tr}}\nolimits}

\newcommand{\A}{\mathfrak{A}}
\newcommand{\liea}{\mathfrak{a}}

\newcommand{\lieg}{\mathfrak{g}}

\newcommand{\liek}{\mathfrak{k}}

\newcommand{\liem}{\mathfrak{m}}

\newcommand{\liep}{\mathfrak{p}}

\newcommand{\liesp}{\mathfrak{sp}}

\newcommand{\RP}{\ensuremath{\R\mathrm{P}}}
\newcommand{\CP}{\ensuremath{\C\mathrm{P}}}
\newcommand{\HP}{\ensuremath{\HH\mathrm{P}}}
\newcommand{\KP}{\ensuremath{\K\mathrm{P}}}

\newcommand{\frakJ}{\mathfrak{J}}

\newcommand{\bbS}{\mathbb{S}}

\newcommand{\Sph}{\bbS}
\newcommand{\gR}[2]{\langle #1,#2\rangle_{\R}}

\newcommand{\beweis}{\begingroup\footnotesize \emph{Proof. }}
\newcommand{\beweisende}{\strut\hfill $\Box$\par\medskip\endgroup}
\newcommand{\Mengegr}[2]{\{\,#1\,{\bigr |}\,#2\,\}}

\newcommand{\wt}{\widetilde}
\newcommand{\wh}{\widehat}

\newcommand{\C}{\mathchoice {\setbox0=\hbox{$\displaystyle\rm
C$}\hbox{\hbox to0pt{\kern0.4\wd0\vrule
height0.95\ht0\hss}\box0}}{\setbox0=\hbox{$\textstyle\rm C$}\hbox{\hbox
to0pt{\kern0.4\wd0\vrule
height0.95\ht0\hss}\box0}}{\setbox0=\hbox{$\scriptstyle\rm C$}\hbox{\hbox
to0pt{\kern0.4\wd0\vrule
height0.95\ht0\hss}\box0}}{\setbox0=\hbox{$\scriptscriptstyle\rm
C$}\hbox{\hbox to0pt{\kern0.4\wd0\vrule height0.95\ht0\hss}\box0}}}

\newcommand{\hfilll}{\hskip 0pt plus 1filll}

\theoremstyle{plain} 

\theorembodyfont{\rmfamily\mdseries\itshape}
\newtheorem{Def}{Definition}[section]
\newtheorem{Prop}[Def]{Proposition}
\newtheorem{Theorem}[Def]{Theorem}
\newtheorem{Lemma}[Def]{Lemma}

\theorembodyfont{\rmfamily\mdseries\upshape}

\newtheorem{Remark}[Def]{Remark}
\newtheorem{Remarks}[Def]{Remarks}

\begin{document}
\selectlanguage{english}

\title{Totally geodesic submanifolds of the complex \\ and the quaternionic 2-Grassmannians}
\author{Sebastian Klein${}^1$}
\date{September 17, 2007}
\maketitle
\footnotetext[1]{This work was supported by a fellowship within the Postdoc-Programme of the German Academic Exchange Service (DAAD).}

\abstract{\textbf{Abstract.} In this article, I classify the totally geodesic submanifolds in the complex \,$2$-Grassmannian \,$G_2(\C^{n+2})$\,
and in the quaternionic \,$2$-Grassmannian \,$G_2(\HH^{n+2})$\,. It turns out that for both of these spaces, 
the earlier classification of maximal totally geodesic submanifolds in Riemannian symmetric spaces of rank \,$2$\, 
published by \textsc{Chen} and \textsc{Nagano} in \cite{Chen/Nagano:totges2-1978} is incomplete.
For example, \,$G_2(\HH^{n+2})$\, with \,$n \geq 5$\, contains totally geodesic submanifolds isometric to a \,$\HP^2$\,, its metric scaled such that
the \emph{minimal} sectional curvature is \,$\tfrac15$\,; they are maximal in \,$G_2(\HH^7)$\,. 
\,$G_2(\C^{n+2})$\, with \,$n \geq 4$\, contains totally geodesic submanifolds which are isometric to a \,$\CP^2$\, contained in the \,$\HP^2$\, mentioned
above; they are maximal in \,$G_2(\C^6)$\,. Neither submanifolds are mentioned in \cite{Chen/Nagano:totges2-1978}.}
\bigskip

\textbf{Author's address.} \\
Sebastian Klein \\ Department of Mathematics, Aras na Laoi \\ University College Cork \\ Cork City \\ Ireland \\
\texttt{s.klein@ucc.ie}

\bigskip


\bigskip

\textbf{MS classification numbers:} 53C35 (Primary); 53C17

\section{Introduction}
\label{Se:intro}

The classification of the totally geodesic submanifolds in the most important Riemannian symmetric
spaces of rank \,$2$\,, namely the 2-Grassmannians, is an interesting and significant problem of Riemannian geometry.
The totally geodesic submanifolds in the oriented real \,$2$-Grassmannians \,$G_2^+(\R^{n+2})$\,
(equivalent to the complex quadrics \,$Q^n \subset \CP^{n+1}$\,) have already been classified in my paper \cite{Klein:2007-claQ};
in the present paper I will solve the classification of totally geodesic submanifolds in the  complex and the quaternionic 2-Grassmannians. 

\bigskip

It should be mentioned that 
already \textsc{Chen} and \textsc{Nagano} gave what they claimed to be a complete classification of the isometry types of maximal totally geodesic submanifolds in all
Riemannian symmetric spaces of rank \,$2$\, in \S 9 of their paper \cite{Chen/Nagano:totges2-1978} based on their \,$(M_+,M_-)$-method.
However, as it will turn out in the present paper,
their classification is faulty: 
several types of totally geodesic submanifolds are missing from their list, both for the space \,$G_2(\C^{n+2})$\, and for the space \,$G_2(\HH^{n+2})$\,;
see Remarks~\ref{R:claH:missing} and \ref{Rs:claC:claC}(a) for a more detailed description. Even apart from these omissions, Chen's and Nagano's investigation
is not satisfactory, as they name only the isometry type of the totally geodesic submanifolds, without giving any description of their position
in the ambient space. (Such a description might, for example, be constituted by giving explicit totally geodesic, isometric embeddings for the various congruence classes
of totally geodesic submanifolds, or at least by describing the tangent spaces of the totally geodesic submanifolds (i.e.~the Lie triple systems)
as subspaces of the tangent space of the ambient symmetric space in an explicit way.)

Besides the results of Chen and Nagano,
various other \emph{partial} results concerning totally geodesic submanifolds in Grassmann manifolds have been obtained: For example, \textsc{Wolf} has
obtained a classification of the totally geodesic submanifolds of the \,$1$-Grassmannians (i.e.~the projective spaces) in 
\cite{Wolf:1963-elliptic}. In \cite{Wolf:1963-spheres} and \cite{Wolf:1963-elliptic} he also classified those totally geodesic submanifolds
of a Grassmannian manifold \,$G_r(\K^n)$\, in which any two distinct elements have zero intersection as subspaces of \,$\K^n$\,; it turns out
that such totally geodesic submanifolds are necessarily of rank \,$1$\,. We will use the latter classification result by Wolf here to handle one specific case
of the general classification of totally geodesic submanifolds in \,$G_2(\HH^n)$\, (namely the case where the submanifold has constant ``characteristic angle''
\,$\tfrac\pi4$\,; it will turn out that it then satisfies the hypothesis of the classification by Wolf).

Some specific types of totally geodesic submanifolds have been classified in all Riemannian symmetric spaces \,$M$\,. The two most important results of this
kind are the classification of reflective submanifolds (i.e.~those submanifolds of \,$M$\, which are connected components of the fixed point set 
of an involutive isometry on \,$M$\,) due to \textsc{Leung} (see \cite{Leung:reflective-1979}) and the classification of totally geodesic submanifolds
of maximal rank due to \textsc{Zhu} and \textsc{Liang} (see \cite{Zhu/Liang:maxrank-2004}). Moreover, \textsc{Nagano} and \textsc{Sumi} gave a classification
of totally geodesic spheres in Riemannian symmetric spaces in \cite{Nagano/Sumi:1991-spheres}.

These partial results notwithstanding, the problem of classifying all totally geodesic submanifolds of arbitrary rank
in Riemannian symmetric spaces still remains open, even for the symmetric spaces of rank \,$2$\,. 

\bigskip

The usual strategy for the classification of totally geodesic submanifolds in a Riemannian symmetric space \,$M=G/K$\,,
which is used also here, is as follows.
Let \,$\lieg = \liek \oplus \liem$\, be the decomposition of the Lie algebra of \,$G$\, induced by the
symmetric structure of \,$M$\,. As it is well-known, the Lie triple systems \,$\liem'$\, in \,$\liem$\, (i.e.~the linear
subspaces \,$\liem' \subset \liem$\, which satisfy \,$[[\liem',\liem'],\liem'] \subset \liem'$\,) are in one-to-one
correspondence with the (connected, complete) totally geodesic submanifolds \,$M_{\liem'}$\, of \,$M$\,
running through the ``origin point'' \,$p_0 = eK \in M$\,,
the correspondence being that \,$M_{\liem'}$\, is characterized by \,$p_0 \in M_{\liem'}$\, and
\,$T_{p_0}M_{\liem'} = \tau(\liem')$\,, where \,$\tau: \liem \to T_{p_0}M$\, is the canonical isomorphism.

Thus the task of classifying the totally geodesic submanifolds of \,$M$\, splits into two steps: (1) To classify the Lie triple systems
in \,$\liem$\,, and (2) for each of the Lie triple systems \,$\liem'$\, found in the first step,
to construct a (connected, complete) totally geodesic
submanifold \,$M_{\liem'}$\, of \,$M$\, so that \,$p_0 \in M_{\liem'}$\, and \,$\tau^{-1}(T_{p_0}M_{\liem'}) = \liem'$\, holds.

Herein, step (1) is the one which generally poses the more significant difficulties. As an approach to accomplishing this step, we describe
in Section~\ref{Se:generallts} for an arbitrary Riemannian symmetric space \,$M$\, of compact type
relations between the roots and root spaces of \,$M$\, and the
roots resp.~root spaces of its totally geodesic submanifolds (regarded as symmetric subspaces).
These relations provide conditions which are necessary for a linear subspace \,$\liem'$\, of \,$\liem$\, to be a Lie triple system. However,
these conditions are not generally sufficient, and therefore a specific investigation needs to be made to see which of the linear subspaces of \,$\liem$\,
satisfying the conditions are in fact Lie triple systems. This investigation is carried out for \,$G_2(\HH^{n+2})$\, in Section~\ref{Se:claH};
it is the laborious part of the proof of the classification theorem (Theorem~\ref{T:claH:claH}). 

It should be emphasized that to carry out this investigation for a given Riemannian symmetric space, it does not suffice to know the 
(restricted) root system (with multiplicities) of that space, or equivalently, the action of the Jacobi operators \,$R(\,\cdot\,,v)v$\,
on the various root spaces. Rather, full control of the curvature tensor is needed. For this reason, we give a description of the curvature tensor,
and associated objects like the Cartan subalgebras, the roots and the root spaces, of \,$G_2(\HH^{n+2})$\, in Section~\ref{Se:curv}. 
We further prepare the classification by giving an explicit description of the orbits of the isotropy action of \,$G_2(\HH^{n+2})$\, in Section~\ref{Se:orbits}.

The central part of the present paper is the classification of the Lie triple systems in \,$G_2(\HH^{n+2})$\,, which is carried out in Section~\ref{Se:claH};
the result is found in Theorem~\ref{T:claH:claH}. By inspection of the root systems (with multiplicities) of the various Lie triple systems found,
we can already tell the local isometry type of the corresponding totally geodesic submanifolds of \,$G_2(\HH^{n+2})$\,. But to determine the global
isometry type, some more considerations are needed. Also it is desirable to describe the actual totally geodesic submanifolds corresponding to the
various types of Lie triple systems as explicitly as possible. This is done in Section~\ref{Se:Hembed}. 

Finally, we classify the Lie triple systems and the totally geodesic submanifolds of \,$G_2(\C^{n+2})$\, in Section~\ref{Se:claC}. 
Because \,$G_2(\C^{n+2})$\, is a totally geodesic submanifold of \,$G_2(\HH^{n+2})$\,, we can obtain this classification simply by checking
which of the totally geodesic submanifolds of \,$G_2(\HH^{n+2})$\, are contained in a given totally geodesic \,$G_2(\C^{n+2}) \subset G_2(\HH^{n+2})$\,. 

\bigskip

The results of the present paper were obtained by me while working at the University College Cork under the advisorship of Professor J.~Berndt. 
I would like to thank him for his dedicated support and guidance, as well as his generous hospitality.

\section{General facts on Lie triple systems}
\label{Se:generallts}

In this section we suppose that \,$M=G/K$\, is any Riemannian symmetric space of compact type. 
We consider the decomposition \,$\lieg = \liek \oplus \liem$\, of the Lie algebra \,$\lieg$\, of \,$G$\, induced by the symmetric
structure of \,$M$\,. Because \,$M$\, is of compact type, the Killing form \,$\vkap: \lieg \times \lieg \to \R,\;(X,Y) \mapsto \tr(\ad(X) \circ \ad(Y))$\,
is negative definite, and therefore \,$\g{\,\cdot\,}{\,\cdot\,} := -c \cdot \vkap$\, gives rise to a Riemannian metric on \,$M$\, for arbitrary \,$c \in \R_+$\,.
In the sequel we suppose that \,$M$\, is equipped with such a Riemannian metric.

Let us fix notations concerning flat subspaces, roots and root spaces of \,$M$\, (for the corresponding theory, see for example \cite{Loos:1969-2}, Section~V.2):
A linear subspace \,$\liea \subset \liem$\, is called \emph{flat} if \,$[\liea,\liea] = \{0\}$\, holds. The maximal flat subspaces of \,$\liem$\,
are all of the same dimension, called the \emph{rank} of \,$M$\, (or \,$\liem$) and denoted by \,$\rk(M)$\, or \,$\rk(\liem)$\,; they are called the
\emph{Cartan subalgebras} of \,$\liem$\,. If a Cartan subalgebra 
\,$\liea \subset \liem$\, is fixed, we put for any linear form \,$\lambda \in \liea^*$\,
$$ \liem_\lambda := \Menge{\,X \in \liem\,}{\,\forall Z \in \liea: \ad(Z)^2 X = -\lambda(Z)^2X\,} $$
and consider the \emph{root system}
$$ \Delta(\liem,\liea) := \Menge{\,\lambda \in \liea^*\setminus \{0\}\,}{\,\liem_\lambda \neq \{0\}\,} $$
of \,$\liem$\, with respect to \,$\liea$\,. The elements of \,$\Delta(\liem,\liea)$\, are called \emph{roots} of \,$\liem$\, with respect to \,$\liea$\,,
for \,$\lambda \in \Delta(\liem,\liea)$\, \,$\liem_\lambda$\, is called the \emph{root space} corresponding to \,$\lambda$\,, and \,$n_{\lambda} := \dim(\liem_{\lambda})$\,
is called the \emph{multiplicity} of the root \,$\lambda$\,. If we fix a system of positive roots \,$\Delta_+ \subset \Delta(\liem,\liea)$\,
(i.e.~we have \,$\Delta_+ \dot{\cup} (-\Delta_+) = \Delta(\liem,\liea)$\,), 
we obtain the \emph{root space decomposition} of \,$\liem$\,:
\begin{equation}
\label{eq:roots:mdecomp}
\liem = \liea \;\oplus\; \bigoplus_{\lambda \in \Delta_+} \liem_\lambda \;.
\end{equation}
The \emph{Weyl group} \,$W(\liem,\liea)$\, is the transformation group on \,$\liea$\, generated by the reflections in the hyperplanes \,$\Menge{v \in \liea}{\lambda(v)=0}$\,
(where \,$\lambda$\, runs through \,$\Delta(\liem,\liea)$\,);
it can be shown that the root system \,$\Delta(\liem,\liea)$\, is invariant under the action of \,$W(\liem,\liea)$\,. 

\bigskip

Let us now consider a Lie triple system \,$\liem' \subset \liem$\,, i.e.~\,$\liem'$\, is a linear subspace of \,$\liem$\, so that
\,$[\,[\liem',\liem'] \,,\, \liem'\,] \subset \liem'$\, holds. In spite of the fact that the symmetric space corresponding to \,$\liem'$\,
does not need to be of compact type (it can contain Euclidean factors), it can be shown easily that the usual statements of the root space theory
for symmetric spaces of compact type carry over to \,$\liem'$\,, see \cite{Klein:2007-claQ}. 

More specifically, the maximal flat subspaces of \,$\liem'$\, are all of the same dimension
(again called the \emph{rank} of \,$\liem'$\,), and they are again called the \emph{Cartan subalgebras} of \,$\liem'$\,. For any Cartan subalgebra \,$\liea'$\,
of \,$\liem'$\,, there exists a Cartan subalgebra \,$\liea$\, of \,$\liem$\, so that \,$\liea' = \liea \cap \liem'$\, holds. With respect to any
Cartan subalgebra \,$\liea'$\, of \,$\liem'$\, we have a root system \,$\Delta(\liem',\liea')$\, (defined analogously as for \,$\liem$\,) 
and the corresponding root space decomposition
\begin{equation}
\label{eq:roots:m'decomp}
\liem' = \liea' \;\oplus\; \bigoplus_{\alpha \in \Delta_+(\liem',\liea')} \liem'_{\alpha}
\end{equation}
(with a system of positive roots \,$\Delta_+(\liem',\liea') \subset \Delta(\liem',\liea')$\,); we also again call \,$n_{\alpha}' := \dim(\liem_{\alpha}')$\,
the multiplicity of \,$\alpha \in \Delta(\liem',\liea')$\,. 
\,$\Delta(\liem',\liea')$\, is again invariant under the action of the corresponding Weyl group \,$W(\liem',\liea')$\,. 
It should be noted, however,
that in the case where a Euclidean factor is present in \,$\liem'$\,, \,$\Delta(\liem',\liea')$\, does not span \,$(\liea')^*$\,.

The following proposition describes the relation between the root space decompositions~\eqref{eq:roots:m'decomp} of \,$\liem'$\, and 
\eqref{eq:roots:mdecomp} of \,$\liem$\,. In particular, it shows the 
extent to which the the position of the individual root spaces \,$\liem'_\alpha$\, of \,$\liem'$\, is adapted to the root space
decomposition~\eqref{eq:roots:mdecomp} of the ambient space \,$\liem$\,. We will base our classification
of the Lie triple systems in \,$G_2(\HH^n)$\, on these relations.

\begin{Prop}
\label{P:cla:subroots:subroots-neu}
Let \,$\liea'$\, be a Cartan subalgebra of \,$\liem'$\,, and let \,$\liea$\, be a Cartan subalgebra of \,$\liem$\, so that \,$\liea' = \liea \cap \liem'$\, holds.
\begin{enumerate}
\item
The roots resp.~root spaces of \,$\liem'$\, and of \,$\liem$\, are related in the following way:
\begin{gather}
\label{eq:cla:subroots:subroots-neu:toshow-Delta}
\Delta(\liem',\liea') \;\subset\; \;\bigr\{\;\lambda|\liea'\; \bigr| \;\lambda \in \Delta(\liem,\liea), \lambda|\liea' \neq 0\; \bigr\} \;.\\
\label{eq:cla:subroots:subroots-neu:toshow-liemalpha}
\textstyle
\forall\alpha \in \Delta(\liem',\liea')\;:\; \liem_\alpha' = \left( \bigoplus_{\substack{\lambda \in \Delta(\liem,\liea) \\ \lambda|\liea' = \alpha}} \liem_\lambda \right) \;\cap\; \liem' \; .
\end{gather}
\item
We have \,$\rk(\liem') = \rk(\liem)$\, if and only if \,$\liea' = \liea$\, holds. If this is the case, then we have
\begin{equation}
\label{eq:cla:subroots:subroots-neu:c}
\Delta(\liem',\liea') \subset \Delta(\liem,\liea) \;,\quad \forall \alpha \in \Delta(\liem',\liea') : \liem'_\alpha = \liem_\alpha \cap \liem' \; . 
\end{equation}
\end{enumerate}
\end{Prop}

\beweis
See \cite{Klein:2007-claQ}, the proof of Proposition~2.1.
\beweisende

For the remainder of the section, we fix a Cartan subalgebra \,$\liea'$\, of \,$\liem'$\,, and let \,$\liea$\, be 
any Cartan subalgebra of \,$\liem$\, so that \,$\liea' = \liea \cap \liem'$\, holds.

\begin{Def}
\label{D:cla:subroots:Elemcomp}
Let \,$\alpha \in \Delta(\liem',\liea')$\, be given.
Recall that by Proposition~\ref{P:cla:subroots:subroots-neu}(a) there exists at least one root \,$\lambda \in \Delta(\liem,\liea)$\, with
\,$\lambda|\liea' = \alpha$\,. We call \,$\alpha$\, 
\begin{enumerate}
\item \emph{elementary}, if there exists only one root \,$\lambda \in \Delta(\liem,\liea)$\, with \,$\lambda|\liea' = \alpha$\,;
\item \emph{composite}, if there exist at least two different roots \,$\lambda, \mu \in \Delta(\liem,\liea)$\, with \,$\lambda|\liea' = \alpha = \mu|\liea'$\,.
\end{enumerate}
\end{Def}

Elementary roots play a special role:
If \,$\alpha \in \Delta(\liem',\liea')$\, is elementary, then the root space \,$\liem_\alpha'$\, 
is contained in the root space \,$\liem_\lambda$\,, where \,$\lambda \in \Delta(\liem,\liea)$\, is the unique root with \,$\lambda|\liea' = \alpha$\,.
As we will see in Proposition~\ref{P:cla:subroots:Comp} below, this property
causes restrictions for the possible positions (in relation to \,$\liea'$\,) of \,$\lambda$\,. The exploitation of these restrictions will play
an important role in the classification of the rank \,$1$\, Lie triple systems of \,$G_2(\HH^n)$\, in Section~\ref{SSe:claH:rk1}. 

It should also be mentioned that in the case \,$\rk(\liem') = \rk(\liem)$\, we have \,$\liea' = \liea$\,, and therefore in that case
every \,$\alpha \in \Delta(\liem',\liea')$\, is elementary (compare Proposition~\ref{P:cla:subroots:subroots-neu}(b)).

For any linear form \,$\lambda \in \liea^*$\, we now denote by \,$\lambda^\sharp$\, the Riesz vector
corresponding to \,$\lambda$\,, i.e.~the vector \,$\lambda^\sharp \in \liea$\, characterized by \,$\g{\,\cdot\,}{\lambda^\sharp} = \lambda$\,. Here \,$\g{\,\cdot\,}{\,\cdot\,} = -c\cdot \vkap$\, is again the inner product obtained from the Killing form \,$\vkap$\,
of \,$\lieg$\,. 

\begin{Prop}
\label{P:cla:subroots:Comp}
Let \,$\alpha \in \Delta(\liem',\liea')$\, be given. 
\begin{enumerate}
\item
If \,$\alpha$\, is elementary 
and \,$\lambda \in \Delta(\liem,\liea)$\, is the unique root with \,$\lambda|\liea' = \alpha$\,, then we have \,$ \lambda^\sharp \in \liea'$\,. 
\item
If \,$\alpha$\, is composite and \,$\lambda, \mu \in \Delta(\liem,\liea)$\, are two different roots with
\,$\lambda|\liea' = \alpha = \mu|\liea'$\,, then \,$\lambda^\sharp - \mu^\sharp$\, is orthogonal to \,$\liea'$\,. 
\end{enumerate}
\end{Prop}

\beweis
For (a) see \cite{Klein:2007-claQ}, the proof of Proposition~2.3(a). (b) is obvious.
\beweisende

\begin{Prop}
\label{P:cla:skew}
Suppose that \,$\alpha \in \Delta(\liem',\liea')$\, is a composite root such that there exist precisely two roots \,$\lambda,\mu \in \Delta(\liem,\liea)$\, 
with \,$\lambda|\liea' = \alpha = \mu|\liea'$\,. Further suppose that \,$\alpha^\sharp$\, can be written as a linear combination \,$\alpha^\sharp = a\,\lambda^\sharp
+ b\,\mu^\sharp$\, with non-zero \,$a,b \in \R$\,. 

Then we have \,$a,b > 0$\,, and
there exists a linear subspace \,$\liem_{\lambda}' \subset \liem_\lambda$\, and an isometric linear map \,$\Phi: \liem_{\lambda}' \to \liem_\mu$\,
so that
\begin{equation}
\label{eq:cla:skew:skew}
\liem_{\alpha}' = \Menge{x+ \sqrt{\tfrac{b}{a}}\, \Phi(x)}{x \in \liem_{\lambda}'}
\end{equation}
holds. In particular we have \,$n_{\alpha}' \leq \min\{n_\lambda,n_\mu\}$\,. 
\end{Prop}

\beweis
First we note that the hypotheses imply that \,$\lambda$\, and \,$\mu$\, are linearly independent: Assume to the contrary that \,$\mu = c\,\lambda$\, holds
with some \,$c \in \R$\,; we have \,$c \not\in \{0,1\}$\,.
We would then have \,$\alpha^\sharp = a\,\lambda^\sharp + b\,\mu^\sharp = (a+bc)\lambda^\sharp = (\tfrac{1}{c}a+b)\mu^\sharp$\,, which because of 
\,$0 \neq \alpha = \lambda|\liea' = \mu|\liea'$\, implies \,$1 = a+bc = \tfrac{1}{c}a+b$\,. The equality \,$a+bc = \tfrac{1}{c}a+b$\, implies \,$a=-cb$\, because
of \,$c\neq 1$\,. Thus we have \,$a+bc=0$\,, in contradiction to \,$1 = a+bc$\,. 

In the sequel we make use of the fact that for every \,$\lambda \in \Delta(\liem,\liea)$\, and every \,$v \in \liem_\lambda$\,, there exists one and only one
vector \,$\wh{v} \in \liek$\, which is ``related'' to \,$v$\, in the sense that 
\begin{equation}
\label{eq:cla:skew:rel1}
\forall H \in \liea \; : \; \bigr( \; \ad(H)v = \lambda(H)\cdot \wh{v} \qmq{and} \ad(H)\wh{v} = -\lambda(H)\cdot v \; \bigr)
\end{equation}
holds; we then also have
\begin{equation}
\label{eq:cla:skew:rel2}
[v,\wh{v}] = \|v\|^2 \cdot \lambda^\sharp \; .
\end{equation}
(For example, see \cite{Loos:1969-2}, Lemma~VI.1.5(a), p.~62.) The analogous statement holds in the Lie triple system \,$\liem'$\,. 

Let us fix \,$H \in \liea'$\, with \,$\alpha(H) \neq 0$\,, and let \,$v \in \liem_\alpha'$\, be given. By Proposition~\ref{P:cla:subroots:subroots-neu}(a)
there exist unique \,$v_\lambda \in \liem_\lambda$\, and \,$v_\mu \in \liem_\mu$\, so that \,$v = v_\lambda + v_\mu$\, holds. We now calculate \,$R(H,v)v$\, in two 
different ways: On the one hand, we have 
\begin{equation}
\label{eq:cla:skew:R1}
R(H,v)v = -[[H,v],v] \overset{\eqref{eq:cla:skew:rel1}}{=} -\alpha(H)\cdot [\wh{v},v] \overset{\eqref{eq:cla:skew:rel2}}{=} \alpha(H) \cdot \|v\|^2 \cdot \alpha^\sharp
= \alpha(H) \cdot \|v\|^2 \cdot (a\,\lambda^\sharp + b\,\mu^\sharp) \; ; 
\end{equation}
in particular we see \,$R(H,v)v \in \liea' \subset \liea$\,. On the other hand we have
$$ R(H,v_\lambda)v_\lambda = \underbrace{\lambda(H)}_{= \alpha(H)} \cdot \|v_\lambda\|^2 \cdot \lambda^\sharp \qmq{and}
R(H,v_\mu)v_\mu = \underbrace{\mu(H)}_{= \alpha(H)} \cdot \|v_\mu\|^2 \cdot \mu^\sharp $$
by the analogous calculation as in \eqref{eq:cla:skew:R1}, and therefore
\begin{align*}
R(H,v)v & = R(H,v_\lambda)v_\lambda + R(H,v_\mu)v_\mu + R(H,v_\lambda)v_\mu + R(H,v_\mu)v_\lambda \\
& = \underbrace{\alpha(H) \cdot (\|v_\lambda\|^2\,\lambda^\sharp + \|v_\mu\|^2\,\mu^\sharp)}_{(*)} + \underbrace{R(H,v_\lambda)v_\mu + R(H,v_\mu)v_\lambda}_{(\dagger)} \; . 
\end{align*}
Both \,$R(H,v)v$\, and the term marked \,$(*)$\, are members of \,$\liea$\,, whereas the term marked \,$(\dagger)$\, is a member of \,$\liem_{\lambda+\mu}
\oplus \liem_{\lambda-\mu}$\,, and is therefore orthogonal to \,$\liea$\, (the linear independence of \,$\lambda$\, and \,$\mu$\, implies \,$\lambda\pm\mu \neq 0$\,).
It follows that \,$(\dagger)$\, vanishes, and thus we have
\begin{equation}
\label{eq:cla:skew:R2}
R(H,v)v = \alpha(H) \cdot (\|v_\lambda\|^2\,\lambda^\sharp + \|v_\mu\|^2\,\mu^\sharp) \; . 
\end{equation}
By comparing Equations~\eqref{eq:cla:skew:R1} and \eqref{eq:cla:skew:R2} we now obtain
$$ \|v\|^2 \cdot (a\,\lambda^\sharp + b\,\mu^\sharp) = \|v_\lambda\|^2 \,\lambda^\sharp + \|v_\mu\|^2 \,\mu^\sharp $$
and therefore because of the linear independence of \,$\lambda$\, and \,$\mu$\,
$$ \|v_\lambda\|^2 = a\cdot \|v\|^2 \qmq{and} \|v_\mu\|^2 = b\cdot \|v\|^2 \; . $$
It follows that we have \,$a,b > 0$\,, and that the linear maps 
$$ \Phi_\lambda : \liem_\alpha' \to \liem_\lambda, \; v = v_\lambda+v_\mu \mapsto \tfrac{1}{\sqrt{a}}\,v_\lambda \qmq{and}
\Phi_\mu : \liem_\alpha' \to \liem_\mu, \; v = v_\lambda+v_\mu \mapsto \tfrac{1}{\sqrt{b}}\,v_\mu $$
are isometric, in particular they are injective. Now consider the linear subspace \,$\liem_\lambda' := \Phi_\lambda(\liem_\alpha')$\, of \,$\liem_\lambda$\,
and the linear isometry \,$\Phi := \Phi_\mu \circ (\Phi_\lambda)^{-1}: \liem_\lambda' \to \liem_\mu$\,. For any \,$v = v_\lambda + v_\mu \in \liem_\alpha'$\, as before,
we have \,$\Phi(v_\lambda) = \Phi_\mu(\Phi_\lambda^{-1}(v_\lambda)) = \Phi_\mu(\sqrt{a}\,v) = \sqrt{\tfrac{a}{b}}\,v_\mu$\, and therefore
\,$v = v_\lambda + \sqrt{\tfrac{b}{a}}\,\Phi(v_\lambda)$\,. Hence we have shown~\eqref{eq:cla:skew:skew}. It follows that
\,$n_\alpha' = \dim(\liem_\alpha') = \dim(\liem_\lambda') \leq \dim(\liem_\lambda) = n_\lambda$\, holds; by exchanging the roles of \,$\lambda$\, and \,$\mu$\,
we also get \,$n_\alpha' \leq n_\mu$\,. 
\beweisende

\section{The curvature tensor of quaternionic 2-Grassmannians}
\label{Se:curv}

\paragraph{Generalities on quaternionic linear spaces.}
We denote by \,$\HH$\, the skew-field of \emph{quaternions}, by \,$\overline{c}$\, the \emph{conjugate} of a quaternion \,$c \in \HH$\, and by
\,$\IM(\HH) := \Menge{c \in \HH}{\overline{c}=-c}$\, the real-3-dimensional space of \emph{purely imaginary} quaternions. A \emph{canonical basis} of \,$\IM(\HH)$\,
is an orthonormal basis \,$(i,j,k)$\, of \,$\IM(\HH)$\, so that \,$k=ij$\, holds. For any \,$c \in \HH$\,, we denote by \,$\RE(c) := \tfrac12(c+\overline{c}) \in \R$\,
its \emph{real part}, and by \,$\IM(c) := \tfrac12(c-\overline{c}) \in \IM(\HH)$\, its \emph{imaginary part}.

Let \,$V$\, be a \emph{symplectic space}, i.~e.~a right-linear space over \,$\HH$\, equipped with a quaternionic inner product 
\,$\g{\,\cdot\,}{\,\cdot\,}: V \times V \to \HH$\,, its homogeneity rule is
\begin{equation}
\label{eq:curv:homo}
\forall v,w \in V, \; c,c' \in \HH \; : \; \g{vc}{wc'} = \overline{c} \cdot \g{v}{w} \cdot c' 
\end{equation}
in accordance with the usual conventions. \,$v \perp w$\, stands for \,$\g{v}{w}=0$\,. 
We denote by \,$\Sph(V) := \Menge{v \in V}{\g{v}{v}=1}$\, the unit sphere in \,$V$\,. 
A \emph{symplectic basis} of \,$V$\, is a basis \,$(a_1,\dotsc,a_n)$\, of \,$V$\, so that \,$\g{a_\mu}{a_\nu} = \delta_{\mu\nu}$\,
holds.
The Lie group \,$\Sp(V) := \Menge{B \in \End(V)}{\forall v,w \in V : \g{Bv}{Bw} = \g{v}{w}}$\, is the \emph{symplectic group} of \,$(V,\g{\,\cdot\,}{\,\cdot\,})$\,; 
its Lie algebra
\,$\liesp(V)$\, is isomorphic to the Lie algebra \,$\Menge{X \in \End(V)}{X^* = -X}$\, of skew-adjoint \,$\HH$-linear endomorphisms on \,$V$\,,
equipped with the commutator \,$[X,Y] := X\circ Y - Y \circ X$\, as Lie bracket.
Finally, we note that the quaternionic
inner product also gives rise to a real inner product \,$\gR{\,\cdot\,}{\,\cdot\,} := \RE(\g{\,\cdot\,}{\,\cdot\,})$\, on \,$V$\, seen as a real linear space;
\,$v \perp_{\R} w$\, stands for \,$\gR{v}{w}=0$\,. 

To introduce on \,$V$\, besides the given right-multiplication also a left-multiplication,
we need to single out a real form \,$V_{\R}$\, of \,$V$\,
(i.~e.~\,$V_{\R}$\, is an \,$\R$-linear subspace of \,$V$\, with \,$\dim_{\R}(V_{\R}) = \dim_{\HH}(V)$\,, such that \,$V_{\R} \cdot i$\, is \,$\R$-orthogonal to \,$V_{\R}$\,
for every \,$i \in \IM(\HH)$\,). Then we define the left-multiplication with some given \,$c \in \HH$\, as the right-$\HH$-linear extension of the \,$\R$-linear map
\,$V_{\R} \to (V_{\R}\cdot c),\; v \mapsto v\,c$\,, note that \,$c \, x = x \, c$\, holds for every \,$x \in V_{\R}$\,. The left-multiplication is described
explicitly by
$$ \forall c \in \HH, \; v \in V \; : \; c\cdot v = A(\,A(v)\cdot \overline{c}\,) \; , $$
where \,$A: V \to V$\, is the \,$\R$-linear involution characterized by \,$A|V_{\R} = \id_{V_{\R}}$\, and \,$A|V_{\R}^\perp = -\id_{V_{\R}^\perp}$\,.

\bigskip

Let \,$V$\, and \,$V'$\, be symplectic spaces. We denote the space of \,$\HH$-right-linear maps \,$V' \to V$\, by \,$L(V',V)$\,, and put \,$\End(V) := L(V,V)$\,. 
For every \,$f \in L(V',V)$\, there is a unique \emph{adjoint map} \,$f^* \in L(V,V')$\, characterized by
$$ \forall v,w \in V' \; : \; \g{f(v)}{w} = \g{v}{f^*(w)} \; ; $$
if \,$(a_1,\dotsc,a_n)$\, is a symplectic basis of \,$V'$\,, then \,$f^*$\, is explicitly given by
\begin{equation}
\label{eq:curv:adjointexplicit}
\forall \, v \in V \; : \; f^*(v) = \sum_\nu a_\nu \cdot \g{f(a_\nu)}{v} \; .
\end{equation}

Let us now suppose that real forms \,$V'_{\R}$\, of \,$V'$\, and \,$V_{\R}$\, of \,$V$\, have been singled out, and let us denote the left-multiplications
defined thereby by \,$L_c': V' \to V', \; v \mapsto c\, v$\, and  \,$L_c: V \to V, \; v \mapsto c\, v$\, (for \,$c \in \HH$\,). Then \,$L(V',V)$\,
becomes a \,$\HH$-right- and \,$\HH$-left-linear space by the definitions (for \,$c \in \HH$\, and \,$f \in L(V',V)$\,): 
$$ f \cdot c := f \circ L_c' \qmq{and} c\cdot f := L_c \circ f  \; . $$
Note that if \,$f$\, maps \,$V_{\R}'$\, into \,$V_{\R}$\,, then we have \,$c\cdot f = f\cdot c$\,. 

\paragraph{The quaternionic 2-Grassmannian and its tangent space.}
Let \,$V$\, be a symplectic space of dimension \,$n \geq 2$\, and \,$V'$\, be another symplectic space of dimension \,$2$\,. In the sequel we will study
the quaternionic 2-Grassmannian \,$G_2(V' \oplus V) \cong G_2(\HH^{n+2})$\,, i.e.~the manifold of \,$2$-dimensional, quaternionic subspaces of \,$V' \oplus V$\,. 
It is well-known that this Grassmannian is an irreducible Riemannian symmetric space of compact type and
rank \,$2$\, with respect to the natural action of \,$\Sp(V' \oplus V)$\, on it.
The isotropy group of this action at the point \,$V' \in G_2(V' \oplus V)$\, is \,$\Menge{B \in \Sp(V' \oplus V)}{B(V')=V'}
\cong \Sp(V') \times \Sp(V)$\,, hence \,$G_2(V' \oplus V)$\, is isomorphic 
to the quotient manifold \,$\Sp(V' \oplus V) / \Sp(V') \times \Sp(V)$\,.

The symmetric structure on \,$G_2(V' \oplus V)$\, is induced by the involutive Lie group automorphism
$$ \sigma: \Sp(V'\oplus V) \to \Sp(V' \oplus V), \; B \mapsto SBS^{-1} \;, $$
where \,$S \in \Sp(V' \oplus V)$\, is the symplectic involution characterized by \,$S|V' = \id_{V'}$\, and \,$S|V = -\id_V$\,. The linearization of \,$\sigma$\,
is a Lie algebra involution on the Lie algebra \,$\liesp(V' \oplus V)$\, of skew-adjoint endomorphisms on \,$V' \oplus V$\,.
It induces the Cartan decomposition \,$\liesp(V' \oplus V) = \liek \oplus \liem$\,
corresponding to the symmetric structure; we have 
\begin{align*}
\liek & = \Menge{X \in \liesp(V' \oplus V)}{X(V') \subset V'} \;\cong\; \liesp(V') \oplus \liesp(V) \;, \\
\liem & = \Menge{X \in \liesp(V' \oplus V)}{X(V') \subset V,\;X(V) \subset V'} \;\cong\; L(V',V) \;, 
\end{align*}
where the isomorphisms are given by \,$\liek \to \liesp(V') \oplus \liesp(V),\;X \mapsto (X|V',\,X|V)$\, and \,$\liem \to L(V',V),\; X \mapsto X|V'$\,. 
In the sequel, we will identify the isotropy algebra \,$\liek$\, with \,$\liesp(V') \oplus \liesp(V)$\,, and 
the tangent space \,$\liem$\, of \,$G_2(V' \oplus V)$\, with \,$L(V',V)$\, in this way.

We equip \,$\liem \cong L(V',V)$\, with the usual quaternionic inner product for spaces of linear maps: Let \,$(e_1,e_2)$\, be any symplectic basis of \,$V'$\,,
then we put for \,$v_1,v_2 \in \liem$\,
$$ \g{v_1}{v_2} := \g{v_1(e_1)}{v_2(e_1)} + \g{v_1(e_2)}{v_2(e_2)} \;; $$
this definition does not depend on the choice of the basis \,$(e_1,e_2)$\,. Moreover, this inner product is invariant under the 
action of the isotropy group \,$\Sp(V') \times \Sp(V)$\, of \,$G_2(V' \oplus V)$\, on \,$\liem$\,, which is given by
\begin{equation}
\label{eq:curv:iso}
\forall\, B=(B_1,B_2) \in \Sp(V') \times \Sp(V), \; v \in \liem \; : \; Bv = B_2 \circ v \circ B_1^* \; .
\end{equation}
Therefore the corresponding real inner product \,$\gR{\,\cdot\,}{\,\cdot\,} := \RE(\g{\,\cdot\,}{\,\cdot\,})$\, gives rise to an invariant
Riemannian metric on \,$G_2(V' \oplus V)$\,; we will view \,$G_2(V' \oplus V)$\, with this metric from now on.

\paragraph{The Lie bracket and the curvature tensor.}
It is now easy to get the following formulas for the Lie bracket, which are valid for all \,$X,Y \in \liek$\,
(say \,$X = (X_1,X_2)$\, and \,$Y = (Y_1,Y_2)$\, with \,$X_1,Y_1 \in \liesp(V')$\, and \,$X_2,Y_2 \in \liesp(V)$\,)
and \,$u,v \in \liem \cong L(V',V)$\,:
\begin{align*}
[X,Y] & = ([X_1,Y_1],\;[X_2,Y_2]) \;\in\; \liek \;, \\
[X,v] & = X_2 \circ v + v \circ X_1^* \;\in\; \liem \;, \\
[u,v] & = (v^* \circ u - u^*\circ v,\; v\circ u^* - u\circ v^*) \;\in\; \liek \; . 
\end{align*}
Therefrom we obtain the following formula for the curvature tensor \,$R$\, of \,$G_2(V' \oplus V)$\, via the well-known relationship \,$R(u,v)w = -[[u,v],w]$\,:
\begin{equation}
\label{eq:curv:R}
R(u,v)w = (uv^* - vu^*)w + w(v^*u - u^*v) \; . 
\end{equation}

\paragraph{Conjugations on \,$\boldsymbol{\liem}$\,.}
Our next aim is to describe all Cartan subalgebras of \,$\liem$\,. 
In order to be able to do so in an efficient way, we introduce the concept of a conjugation on \,$V'$\,: Let
$$ \A := \Menge{A \in \Sp(V')}{A^2 = \id_{V'},\; A \neq \pm \id_{V'}} \; ; $$
we call the \,$A \in \A$\, \emph{conjugations} on \,$V'$\,. 
Any \,$A \in \A$\, is symplectically diagonalizable, and its eigenvalues are \,$1$\, and \,$-1$\,, each with (quaternionic) multiplicity \,$1$\,.
It follows that \,$A$\, induces the decomposition \,$V' = V'_+(A) \operp V'_-(A)$\, into the quaternionic-$1$-dimensional eigenspaces
\,$V'_\pm(A) := \Eig(A,\pm 1)$\,. This decomposition of \,$V'$\, also gives rise to a decomposition
of \,$\liem = L(V',V)$\,, namely, we have \,$\liem = L_+(A) \oplus L_-(A)$\, with the quaternionic-$n$-dimensional spaces
\,$L_\pm(A) := \Menge{v \in \liem}{v \circ A = \pm v} = \Menge{v \in \liem}{v|V_\mp'(A)=0}$\,. 
--- It is also a consequence of the consideration of the eigenvalues of \,$A \in \A$\, that
\,$\A$\, is a conjugacy orbit in \,$\Sp(V')$\,.

It should be noted that in this setting, there is no \emph{canonical} isomorphism between \,$V'_+(A)$\, and \,$V'_-(A)$\,, nor between \,$L_+(A)$\, and \,$L_-(A)$\,. 
To describe such isomorphisms, we define
$$ \frakJ_A := \Menge{J \in \Sp(V')}{J^2 = -\id_{V'},\;J\circ A = -A \circ J} \; . $$
Any \,$J \in \frakJ_A$\, interchanges \,$V_+'(A)$\, and \,$V_-'(A)$\, by an isomorphism which respects the inner product on \,$V'$\,. Moreover, 
the map \,$\liem \to \liem,\;v \mapsto v \circ J$\, defines an isomorphism between \,$L_+(A)$\, and \,$L_-(A)$\, which respects the inner product on \,$\liem$\,.
In the sequel we will also write \,$J(v)$\, for \,$v \circ J$\, when \,$v \in \liem$\,. 

If \,$A \in \A$\,, \,$J \in \frakJ_A$\, and a unit vector \,$e_+ \in V_+'(A)$\, is given,
we put \,$e_- := J(e_+) \in V_-'(A)$\,, then \,$(e_+,e_-)$\, is a symplectic basis of \,$V'$\,. We call any such symplectic basis 
\emph{adapted to \,$(A,J)$\,} or simply \emph{adapted to \,$A$\,}.

\paragraph{Cartan subalgebras.}
\begin{Prop}
\label{P:curv:cartan}
Let \,$\liea \subset \liem$\, be a $2$-dimensional real subspace. Then \,$\liea$\, is a Cartan subalgebra of \,$\liem$\, if and only if there exists 
\,$A \in \A$\, and an orthonormal basis \,$(H_+,H_-)$\, of \,$\liea$\, with \,$H_\pm \in L_\pm(A)$\, and \,$H_+(V'_+(A)) \perp H_-(V'_-(A))$\,. 
\end{Prop}

\beweis
First we suppose that there exists a basis \,$(H_+,H_-)$\, of \,$\liea$\, as in the proposition. Using Equation~\eqref{eq:curv:adjointexplicit} and the fact
that we have \,$H_\pm \in L_\pm(A)$\,, one sees that \,$H_-H_+^* = H_+H_-^* = 0$\, holds; using the same equation and the property \,$H_+(V'_+(A)) \perp H_-(V'_-(A))$\,,
one further finds \,$H_+^* H_- = H_-^* H_+ = 0$\,. From these
equations it follows via Equation~\eqref{eq:curv:R} that \,$R(H_+,H_-)H_\pm = 0$\, holds. Therefore \,$\liea$\, is flat, and hence a Cartan subalgebra of \,$\liem$\,
because its dimension \,$2$\, equals the rank of \,$G_2(V'\oplus V)$\,. 

Conversely, let \,$\liea$\, be any Cartan subalgebra of \,$\liem$\,. Because any two Cartan algebras of \,$\liem$\, are conjugate under the isotropy action
on \,$\liem$\, (see \cite{Helgason:1978}, Theorem~V.6.2, p.~246), the Cartan subalgebra \,$\liea$\, is conjugate to a Cartan subalgebra of the type
described in the proposition, and therefore itself of that type.
\beweisende

\paragraph{Roots and root spaces.}
Let \,$\liea = \R H_+ \operp \R H_-$\, be a Cartan subalgebra of \,$\liem$\, as in Proposition~\ref{P:curv:cartan}. 

For the purpose of describing the roots and root spaces of \,$\liem$\,, we 
let \,$(\alpha_+,\alpha_-)$\,
be the basis of the space \,$\liea^*$\, of \,$\R$-linear forms on \,$\liea$\, which is dual to \,$(H_+,H_-)$\,. 
Moreover we fix a basis \,$(e_+,e_-)$\, of \,$V'$\, adapted to \,$A$\,,
and define for any \,$c \in \HH$\, and \,$\eps \in \{\pm 1\}$\, the \,$\HH$-linear map \,$M_{c,\eps}^{(e_+,e_-)} \in L(V',V)$\, by 
\begin{equation}
\label{eq:curv:Mdef}
M_{c,\eps}^{(e_+,e_-)}(e_+) = \tfrac{1}{\sqrt{2}}\,H_-(e_-)\cdot c \qmq{and} M_{c,\eps}^{(e_+,e_-)}(e_-) = \tfrac{1}{\sqrt{2}}\,\eps\,H_+(e_+)\cdot \overline{c} \; .
\end{equation}

We consider the real form \,$V_{\R}' := \spn_{\R}\{e_+,e_-\}$\, of \,$V'$\, and fix a real form \,$V_{\R}$\, of \,$V$\, with \,$H_+(e_+), H_-(e_-) \in V_{\R}$\,.
We define left multiplications on \,$V'$\, and on \,$V$\, with respect to these real forms, and consider \,$L(V',V)$\, as a \,$\HH$-left- and \,$\HH$-right-linear space
via these left multiplications. Note that we then have \,$c\cdot H_\pm = H_\pm \cdot c$\, for any \,$c \in \HH$\,. 

The following table gives
the roots \,$\lambda$\, of \,$G_2(V' \oplus V)$\, with respect to the Cartan subalgebra \,$\liea$\,, the corresponding root spaces
\,$\liem_\lambda$\, and the multiplicities \,$n_\lambda$\,; these data are easily derived from Equation~\eqref{eq:curv:R}:
\begin{center}
\begin{tabular}{|c|c|c|}
\hline
\,$\lambda$\, & \,$\liem_\lambda$\, & \,$n_\lambda$ \\
\hline
\,$\lambda_1 := \alpha_+$\, & \,$\Menge{v \in L_+(A)}{v(e_+) \perp H_+(e_+), H_-(e_-)}$\, & \,$4n-8$\, \\
\,$\lambda_2 := \alpha_-$\, & \,$\Menge{v \in L_-(A)}{v(e_-) \perp H_+(e_+), H_-(e_-)}$\, & \,$4n-8$\, \\
\,$\lambda_3 := \alpha_+ + \alpha_-$\, & \,$\Menge{M_{c,-1}^{(e_+,e_-)}}{c \in \HH}$\, & \,$4$\, \\
\,$\lambda_4 := \alpha_+ - \alpha_-$\, & \,$\Menge{M_{c,1}^{(e_+,e_-)}}{c \in \HH}$\, & \,$4$\, \\
\,$2\,\lambda_1$\, & \,$\IM(\HH) \cdot H_+$\, & \,$3$\, \\
\,$2\,\lambda_2$\, & \,$\IM(\HH) \cdot H_-$\, & \,$3$\, \\
\hline
\end{tabular}
\end{center}
Note that the root spaces (of course) do not in fact depend on the choice of the adapted basis \,$(e_+,e_-)$\, which appears in their description.
However, the individual vectors \,$M_{c,\eps}^{(e_+,e_-)}$\, \emph{do} depend on the choice of the basis; in fact we have for any \,$c,q_+,q_- \in \HH$\,
with \,$|q_\pm|=1$\,:
\begin{equation}
\label{eq:curv:Mtrafo}
M_{c,\eps}^{(e_+\,q_+, e_-\,q_-)} = M_{q_-\,c\,\overline{q_+}\,,\,\eps}^{(e_+,e_-)} \; . 
\end{equation}

From the table of roots we obtain the following, well-known root diagram for \,$G_2(V'\oplus V)$\,:
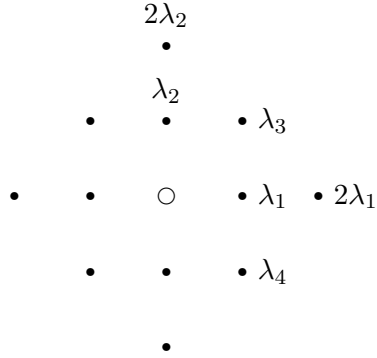
\begin{center}
\begin{minipage}{5cm}
\begin{center}
\strut \\[.3cm]
\setlength{\unitlength}{1.0cm}
\begin{picture}(2,6)
\put(1,3){\circle{0.2}}
\put(2,3){\circle*{0.1}}        
\put(3,3){\circle*{0.1}}        
\put(2,2){\circle*{0.1}}        
\put(1,2){\circle*{0.1}}        
\put(1,1){\circle*{0.1}}        
\put(0,2){\circle*{0.1}}        
\put(0,3){\circle*{0.1}}        
\put(-1,3){\circle*{0.1}}        
\put(0,4){\circle*{0.1}}        
\put(1,4){\circle*{0.1}}        
\put(1,5){\circle*{0.1}}        
\put(2,4){\circle*{0.1}}        
\put(2.2,1.9){{$\lambda_4$}}
\put(2.2,2.9){{$\lambda_1$}}
\put(3.2,2.9){{$2\lambda_1$}}
\put(2.2,3.9){{$\lambda_3$}}
\put(0.8,4.3){{$\lambda_2$}}
\put(0.7,5.3){{$2\lambda_2$}}
\end{picture}
\strut \\[1cm]
\end{center}
\end{minipage}
\end{center}

\begin{Remark}
Of course, the descriptions in the present section can easily be generalized to general quaternionic Grassmannians \,$G_r(\HH^{r+s})$\,,
as well as to Grassmannians over \,$\R$\, or \,$\C$\,.
In particular, Equation~\eqref{eq:curv:R} for the curvature tensor is valid in the general setting.
\end{Remark}

\section{The orbits of the isotropy action on \,$\boldsymbol{G_2(\HH^{n+2})}$\,}
\label{Se:orbits}

\begin{Prop}
\label{P:orbits:orbits}
Suppose \,$v \in \liem \setminus \{0\}$\,. 
\begin{enumerate}
\item
The endomorphism \,$v^*v \in \End(V')$\, is self-adjoint and positive semi-definite. Consequently \,$v^*v$\, is real diagonalizable with
real eigenvalues \,$t_1, t_2 \geq 0$\,. They satisfy the relation \,$t_1 + t_2 = \|v\|^2$\,, hence there is a unique angle \,$\vi(v) \in [0,\tfrac{\pi}{4}]$\,
so that \,$\{t_1,t_2\} = \{\|v\|^2\,\cos(\vi(v))^2, \|v\|^2\,\sin(\vi(v))^2\}$\, holds. We call \,$\vi(v)$\, the \emph{characteristic angle} of \,$v$\,. 

\item
If \,$\liea$\, is a Cartan subalgebra of \,$\liem$\, with \,$v \in \liea$\, (such exist, see \cite{Helgason:1978}, Theorem~V.6.2, p.~246), 
there exists an orthonormal basis \,$(H_+,H_-)$\, of \,$\liea$\, of the kind
described in Proposition~\ref{P:curv:cartan} so that
\begin{equation}
\label{eq:orbits:canrep}
v = \|v\| \cdot (\cos(\vi(v))\,H_+ + \sin(\vi(v))\,H_-)
\end{equation}
holds. We call any such presentation a \emph{canonical representation} of \,$v$\,. 

\begin{center}
\begin{minipage}{5cm}
\begin{center}
\strut \\[.3cm]
\setlength{\unitlength}{0.8cm}
\begin{picture}(3,2.5)
\put(1,1){\circle{0.2}}
\put(2,1){\circle*{0.1}}        
\put(2,0){\circle*{0.1}}        
\put(1,0){\circle*{0.1}}        
\put(0,0){\circle*{0.1}}        
\put(0,1){\circle*{0.1}}        
\put(0,2){\circle*{0.1}}        
\put(1,2){\circle*{0.1}}        
\put(2,2){\circle*{0.1}}        
\put(3,1){\circle*{0.1}}
\put(-1,1){\circle*{0.1}}
\put(1,3){\circle*{0.1}}
\put(1,-1){\circle*{0.1}}
\put(2.1,0.90){{\small $\lambda_1^\sharp$}}
\put(2.1,2.1){{\small $\lambda_3^\sharp$}}
\thicklines
\put(1,1){\vector(2,1){1}}
\dashline[100]{0.1}(1,1)(2,1)
\dashline[100]{0.1}(1,1)(2,2)
\put(1,1){\arc(0.6,0){26}}  
\put(1.32,1.05){{\scriptsize $\vi$}}
\put(2.05,1.4){$v$}
\end{picture}
\strut \\[1cm]
\end{center}
\end{minipage}
\end{center}
\vspace{1cm}

\item
The sets
$$ C_t := \Menge{v \in \Sph(\liem)}{\vi(v)=t} \qmq{with \,$t \in [0,\tfrac{\pi}{4}]$\,} $$
are the orbits of the isotropy action of \,$\Sp(V') \times \Sp(V)$\, on the unit sphere \,$\Sph(\liem)$\,. 
\end{enumerate}
\end{Prop}

{\footnotesize
\beweis
\emph{For (a).}
It is obvious that \,$v^*v$\, is self-adjoint, and for any \,$e \in V'$\, we have \,$\g{e}{v^*v\,e} = \g{v\,e}{v\,e} \geq 0$\,, hence \,$v^*v$\, is positive
semi-definite. It follows that \,$v^*v \in \End(V')$\, is real diagonalizable, and that its eigenvalues \,$t_1,t_2$\, are \,$\geq 0$\,. Let \,$e_1,e_2$\, be
unit eigenvectors corresponding to \,$t_1$\, resp.~\,$t_2$\, which are (in case \,$t_1=t_2$) \,$\HH$-orthogonal to each other. Then
\,$(e_1,e_2)$\, is a symplectic basis of \,$V'$\,, and therefore we have
$$ \|v\|^2 = \g{v}{v} = \g{v\,e_1}{v\,e_1} + \g{v\,e_2}{v\,e_2} = \g{e_1}{v^*v\,e_1} + \g{e_2}{v^*v\,e_2} = t_1 + t_2 \; . $$
Thus \,$\vi(v)$\, can be defined as in the proposition.

\emph{For (b).}
Let \,$\liea$\, be a Cartan subalgebra of \,$\liem$\, with \,$v \in \liea$\,. By Proposition~\ref{P:curv:cartan} there exists \,$A \in \A$\, and an
orthonormal basis \,$(H_+,H_-)$\, of \,$\liea$\, with \,$H_\pm \in L_\pm(A)$\, and \,$H_+(V'_+(A)) \perp H_-(V'_-(A))$\,. Then there exists \,$s \in [0,2\pi)$\,
so that
\begin{equation}
\label{eq:orbits:orbits:vproof}
v = \|v\| \cdot (\cos(s)\,H_+ + \sin(s)\,H_-) 
\end{equation}
holds. By changing the signs of \,$A$\,, \,$H_+$\, and \,$H_-$\, where necessary, we can ensure that \,$v$\, lies in the closed Weyl chamber 
delimited by the root vectors \,$\lambda_1^\sharp = H_+$\, and \,$\lambda_3^\sharp = H_+ + H_-$\,; then we have \,$s \in [0,\tfrac\pi4]$\,. 

It follows from Equation~\eqref{eq:orbits:orbits:vproof} that \,$v^*v$\, is described by the matrix 
\,$\left( \begin{smallmatrix} \|v\|^2\,\cos(s)^2 & 0 \\ 0 & \|v\|^2 \sin(s)^2  \end{smallmatrix} \right)$\, 
with respect to any basis \,$(e_+,e_-)$\, adapted to \,$A$\,; this shows that \,$s = \vi(v)$\, holds. Therefore 
\eqref{eq:orbits:orbits:vproof} is a canonical representation for \,$v$\,.

\emph{For (c).}
First, let \,$v \in \Sph(\liem)$\, and \,$(B_1,B_2) \in \Sp(V') \times \Sp(V)$\, be given, and let \,$\wt{v} := B_2\,v\,B_1^*$\, be the result of the
isotropy action of \,$(B_1,B_2)$\, on \,$v$\,. Then we have by Equation~\eqref{eq:curv:iso}
$$ \wt{v}^* \wt{v} = B_1\,v^*\,B_2^* \; B_2\,v\,B_1^* = B_1\,(v^* v)\,B_1^* \;. $$
Thus \,$\wt{v}^* \wt{v}$\, and \,$v^* v$\, are \,$\Sp(V')$-conjugate to each other, and therefore have the same eigenvalues. This shows that \,$\vi(\wt{v}) = \vi(v)$\,
holds.

Conversely, let \,$v,\wt{v} \in \Sph(\liem)$\, be given with \,$\vi(v) = \vi(\wt{v}) =: s$\,. Then let
\begin{equation}
\label{eq:orbits:orbits:vwtvproof}
v = \cos(s)\,H_+ + \sin(s)\,H_- \qmq{and} \wt{v} = \cos(s)\,\wt{H}_+ + \sin(s)\,\wt{H}_-
\end{equation}
be canonical representations of \,$v$\, resp.~\,$\wt{v}$\, (with \,$H_\pm \in L_\pm(A)$\, and \,$\wt{H}_\pm \in L_\pm(\wt{A})$\,, where \,$A, \wt{A} \in \A$\,,
and \,$H_+(V_+'(A)) \perp H_-(V_-'(A))$\,, \,$\wt{H}_+(V_+'(\wt{A})) \perp \wt{H}_-(V_-'(\wt{A}))$\,). There exists \,$(B_1,B_2) \in \Sp(V') \times \Sp(V)$\,
so that
$$ B_1(e_\pm) = \wt{e}_\pm \qmq{and} B_2(H_\pm(e_\pm)) = \wt{H}_\pm(\wt{e}_\pm) $$
holds, where \,$(e_+,e_-)$\, resp.~\,$(\wt{e}_+,\wt{e}_-)$\, is any basis adapted to \,$A$\, resp.~to \,$\wt{A}$\,. By Equations~\eqref{eq:orbits:orbits:vwtvproof}
we then have \,$B_2\,v\,B_1^* = \wt{v}$\,, so \,$v$\, and \,$\wt{v}$\, are members of the same orbit of the isotropy action on \,$\Sph(\liem)$\,. 
\beweisende

}

\begin{Remarks}
\label{R:orbits:orbits}
\begin{enumerate}
\item
A vector \,$v \in \liem \setminus \{0\}$\, is singular, i.e.~contained in more than one Cartan subalgebra of \,$\liem$\,, if and only if \,$v$\, lies on the
boundary of a Weyl chamber, i.e.~if and only if \,$\vi(v) \in \{0,\tfrac\pi4\}$\, holds. It follows that for \,$0<\vi(v)<\tfrac\pi4$\,, the canonical representation
\eqref{eq:orbits:canrep} is unique.
\item
By Proposition~\ref{P:orbits:orbits}(c), the submanifolds \,$C_t$\, are the orbits of the isotropy action of a rank~\,$2$\, Riemannian symmetric space
on the unit sphere. It follows from results of \textsc{Takagi} and \textsc{Takahashi}, see \cite{Takagi/Takahashi:homhyp-1972}, that 
\,$(C_t)_{0 < t < \tfrac\pi4}$\, is 
a family of isoparametric hypersurfaces on the sphere; the submanifolds \,$C_0$\, and \,$C_{\pi/4}$\, are the focal manifolds of this family.
\end{enumerate}
\end{Remarks}

The characteristic angle \,$\tfrac\pi4$\, plays a special role in many circumstances. This is, for example, evidenced by the fact that
for \,$v \in \liem\setminus\{0\}$\, with \,$\vi(v)=\tfrac\pi4$\,, the canonical representation \eqref{eq:orbits:canrep} of \,$v$\, can be obtained with respect
to any \,$A \in \A$\,. Another specialty of that characteristic angle is exhibited in the following lemma and proposition. 

\begin{Lemma}
\label{L:orbits:pi4}
Let \,$v \in C_{\pi/4}$\, and let \,$\gamma_v: \R \to G_2(V' \oplus V)$\, be the geodesic of \,$G_2(V' \oplus V)$\,
with \,$\gamma_v(0) = V'$\, and \,$\dot{\gamma}_v(0) = v$\,. Then for any \,$t \in \R$\,
either \,$\gamma_v(t)=V'$\, or \,$\gamma_v(t) \cap V' = \{0\}$\, holds (where in the last equation we regard \,$\gamma_v(t)$\, and \,$V'$\, as quaternionic linear subspaces
of \,$V' \oplus V$\,).
\end{Lemma}

\beweis
Let \,$v \in C_{\pi/4}$\, be given. 
As explained in Section~\ref{Se:curv}, \,$v \in \liem = L(V',V)$\, can also be (and originally was)
regarded as an element of \,$\liesp(V' \oplus V)$\,; in block matrix notation
with respect to \,$V' \oplus V$\,, \,$v$\, corresponds to \,$X := \left(\begin{smallmatrix} 0 & -v^* \\ v & 0 \end{smallmatrix} \right) \in \liesp(V' \oplus V)$\,.
Then the geodesic \,$\gamma_v$\, is given by \,$ \gamma_v(t) = \exp(tX)\cdot V'$\,, where \,$\exp: \liesp(V' \oplus V) \to \Sp(V' \oplus V)$\, is
the exponential map of \,$\Sp(V' \oplus V)$\,.

Because of \,$\vi(v)=\tfrac\pi4$\, Proposition~\ref{P:orbits:orbits}(a) shows that \,$v^*v = \tfrac12\, \id_{V'}$\, holds. Using this fact one easily calculates 
for \,$\nu \geq 1$\,
$$ X^{2\nu} = \left(-\tfrac12\right)^\nu \, \left( \begin{smallmatrix} \id & 0 \\ 0 & 2\,vv^* \end{smallmatrix} \right) $$
and for \,$\nu \geq 0$\,
$$ X^{2\nu+1} = \left(-\tfrac12\right)^\nu\,X \; . $$
Therefrom we obtain
\begin{align*}
\exp(tX) & = \sum_{\nu\geq 0} \tfrac{1}{\nu!}(tX)^\nu = \id + \sum_{\nu \geq 1} \tfrac{1}{(2\nu)!}(tX)^{2\nu} + \sum_{\nu \geq 0} \tfrac{1}{(2\nu+1)!}(tX)^{2\nu+1} \\
& = \id + \sum_{\nu \geq 1} \tfrac{(-1)^\nu}{(2\nu)!}\,\left(\tfrac{t}{\sqrt{2}}\right)^{2\nu}\,\left(\begin{smallmatrix} \id & 0 \\ 0 & 2\,vv^* \end{smallmatrix}\right)
        + \sqrt{2}\,\sum_{\nu \geq 0} \tfrac{(-1)^\nu}{(2\nu+1)!}\,\left(\tfrac{t}{\sqrt{2}}\right)^{2\nu+1}\,X \\
& = \id + (\cos(t/\sqrt{2})-1)\,\left( \begin{smallmatrix} \id & 0 \\ 0 & 2\,vv^* \end{smallmatrix} \right) 
        + \sqrt{2}\,\sin(t/\sqrt{2})\,\left( \begin{smallmatrix} 0 & -v^* \\ v & 0 \end{smallmatrix} \right) \; . 
\end{align*}
Denoting by \,$(e_1,e_2)$\, any symplectic basis of \,$V'$\,, we thus have
$$ \gamma_v(t) = \exp(tX)\cdot V' = \spn_{\HH}\{ \; \cos(t/\sqrt{2})\,e_1 + \sqrt{2}\,\sin(t/\sqrt{2})\,v(e_1) \;,\; \cos(t/\sqrt{2})\,e_2 + \sqrt{2}\,\sin(t/\sqrt{2})\,v(e_2) \;\} \; . $$
Because of \,$\vi(v)=\tfrac\pi4$\, the vectors \,$e_1$\,, \,$e_2$\,, \,$v(e_1)$\, and \,$v(e_2)$\, are all non-zero, and pairwise \,$\HH$-orthogonal. Therefore we see that
that \,$\gamma_v(t) = V'$\, holds if and only if we have \,$\sin(t/\sqrt{2})=0$\,, i.e.~\,$t \in \Z(\pi\sqrt{2})$\,, and that otherwise \,$\gamma_v(t) \cap V' = \{0\}$\, holds.
\beweisende

\begin{Prop}
\label{P:orbits:pi4}
Suppose that \,$M$\, is a connected, totally geodesic submanifold of \,$G_2(V' \oplus V)$\,, such that there exists \,$U \in M$\, so that 
 \,$\vi(v)=\tfrac\pi4$\, holds for every \,$v \in T_UM \setminus \{0\}$\,.
Then \,$U_1 \cap U_2 = \{0\}$\, holds for every \,$U_1,U_2 \in M$\, with \,$U_1 \neq U_2$\,. 
\end{Prop}

\beweis
Without loss of generality, we may suppose \,$M$\, to be complete, and therefore a symmetric subspace of \,$G_2(V' \oplus V)$\,. Then some subgroup of the isotropy
group \,$K$\, of \,$G_2(V' \oplus V)$\, acts transitively on \,$M$\,, and therefore the condition \,$\vi(v) = \tfrac\pi4$\, for every \,$v \in T_UM \setminus \{0\}$\, 
holds with respect to every point \,$U \in M$\,. We may also assume without loss of generality that \,$U_1 = V'$\, holds.

Because \,$M$\, is a connected and complete, totally geodesic submanifold of \,$G_2(V' \oplus V)$\,, there exist \,$v \in C_{\pi/4}$\, and \,$t >0$\, so that
the geodesic \,$\gamma_v: \R \to G_2(V' \oplus V)$\, with \,$\gamma_v(0) = V' = U_1$\, and \,$\dot{\gamma}_v(0) = v$\, satisfies
\,$\gamma_v(t) = U_2$\,. We have \,$U_1 \neq U_2$\,, and therefore Lemma~\ref{L:orbits:pi4} shows that \,$U_1 \cap U_2 = \{0\}$\, holds.
\beweisende

The totally geodesic submanifolds \,$M$\, of Grassmannian manifolds with the property that \,$U_1 \cap U_2 = \{0\}$\, holds 
for every \,$U_1,U_2 \in M$\, with \,$U_1 \neq U_2$\, have been classified by \textsc{Wolf} in \cite{Wolf:1963-spheres} and \cite{Wolf:1963-elliptic}.
Via Proposition~\ref{P:orbits:pi4} we will be able to apply Wolf's results to obtain a classification of the totally geodesic submanifolds
of \,$G_2(V'\oplus V)$\, whose tangent spaces have constant characteristic angle \,$\tfrac\pi4$\, as part of the proof of the classification of all totally geodesic
submanifolds in \,$G_2(V'\oplus V)$\, in Section~\ref{Se:claH} of the present paper.

\section{The classification of the Lie triple systems in \,$\boldsymbol{G_2(\HH^{n+2})}$\,}
\label{Se:claH}

\begin{Def}
\label{De:claH:HPtype}
Let \,$V$\, be a symplectic space. An \,$\R$-linear subspace \,$U \subset V$\, is called
\begin{enumerate}
\item \emph{quaternionic} or \emph{of \,$\HP$-type \,$(\HH,\dim_{\HH}(U))$\,}, if \,$U\cdot c \subset U$\, holds for every \,$c \in \HH$\,;
\item \emph{totally complex (with respect to some \,$i \in \Sph(\IM(\HH))$\,)} or \emph{of \,$\HP$-type \,$(\C,\dim_{\C}(U))$\,}, 
if \,$U\cdot i \subset U$\, and \,$U\cdot j \perp_{\R} U$\, holds for every \,$j \in (\R\,i)^{\perp,\IM(\HH)}$\,;
\item \emph{totally real} or \emph{of \,$\HP$-type \,$(\R,\dim_{\R}(U))$\,}, if \,$U\cdot c \perp_{\R} U$\, holds for every \,$c \in \IM(\HH)$\,;
\item \emph{of \,$\HP$-type \,$(\Sph^3)$\,}, if \,$U$\, is a real-$3$-dimensional subspace of a space \,$U' \subset V$\, of \,$\HP$-type \,$(\HH,1)$\,.
\end{enumerate}
Clearly, no subspace of \,$V$\, can be of more than one \,$\HP$-type. 

For any name \,$\tau$\, of a \,$\HP$-type, we define the \emph{dimension} \,$\dim(\tau)$\, and the \emph{width} \,$w(\tau)$\, of \,$\tau$\,:
If \,$\tau$\, is of the form \,$(\K,k)$\, with \,$\K \in \{\R,\C,\HH\}$\,, we put \,$\dim(\tau) := k$\, and \,$w(\tau) := \dim_{\R}(\K) \in \{1,2,4\}$\,.
For \,$\tau = (\Sph^3)$\, we put \,$\dim(\tau)=1$\, and \,$w(\tau) = 3$\,. Then in any case the spaces of \,$\HP$-type \,$\tau$\, have real
dimension \,$w(\tau) \cdot \dim(\tau)$\,. 
\end{Def}

\begin{Remark}
As is well-known (see \cite{Wolf:1963-elliptic}, \S 3), 
the Lie triple systems in a tangent space \,$T_p\HP^n$\, of the quaternionic projective space \,$\HP^n$\, are exactly those \,$\R$-linear subspaces
which have an \,$\HP$-type of dimension \,$\leq n$\,; 
two Lie triple systems are conjugate under the isotropy action of \,$\Sp(1) \times \Sp(n)$\, on \,$T_p \HP^n$\, if and only if
they are of the same \,$\HP$-type. The totally geodesic submanifolds of \,$\HP^n$\, corresponding to Lie triple systems of type \,$(\HH,\ell)$\,, \,$(\C,\ell)$\,,
\,$(\R,\ell)$\, and \,$(\Sph^3)$\, are isometric to \,$\HP^\ell$\,, \,$\CP^\ell$\,, \,$\RP^\ell$\, and \,$\Sph^3 \subset \Sph^4 \cong \HP^1$\,, respectively.
\end{Remark}

As was explained in the Introduction, the pivotal point of the classification of the totally geodesic submanifolds in a Riemannian symmetric
space is the classification of its Lie triple systems. In the present section, we solve the latter problem for the
quaternionic 2-Grassmannians \,$G_2(\HH^{n+2})$\,. 

We remain in the situation of the preceding two sections. In particular \,$\liem \cong L(V',V)$\, is the tangent space of the 
quaternionic 2-Grassmannian \,$G_2(V' \oplus V)$\,.

\newcommand{\Geo}{\ensuremath{\mathrm{Geo}}}
\newcommand{\Gtwo}{\ensuremath{\mathrm{G}_2}}
\newcommand{\PtP}{\ensuremath{\PP \times \PP}}
\newcommand{\StSsqrt}{\ensuremath{\Sph^1 \times \Sph^5}}
\newcommand{\Sptwo}{\ensuremath{\Sp_2}}
\newcommand{\Qthree}{\ensuremath{Q_3}}

\begin{Theorem}
\label{T:claH:claH}
Let \,$\{0\} \neq \liem'$\, be an \,$\R$-linear subspace of \,$\liem$\,. Then \,$\liem'$\, 
is a Lie triple system if and only if \,$\liem'$\, is of one of the following types:

\begin{itemize}
\item \,$\boldsymbol{(\Geo,\vi=t)}$\,, where \,$t \in [0,\tfrac{\pi}{4}]$\,. 

We have \,$\liem' = \R v$\, with some \,$v \in \Sph(\liem')$\, of characteristic angle (see Proposition~\ref{P:orbits:orbits}) \,$\vi(v)=t$\,. 

\item \,$\boldsymbol{(\PP,\vi=0,\tau)}$\,, where \,$\tau$\, is the name of a \,$\HP$-type with \,$\dim(\tau) \leq n$\,.

There exists \,$A \in \A$\, so that \,$\liem'$\, is a subspace of \,$L_+(A)$\, of \,$\HP$-type \,$\tau$\,.

\item \,$\boldsymbol{(\Sph,\vi=\arctan(\tfrac{1}{3}),\ell)}$\,, where \,$\ell \in \{2,3\}$\,. 

There exist \,$A \in \A$\,, a basis \,$(e_+,e_-)$\, adapted to \,$A$\,, \,$H_+ \in L_+(A)$\,, \,$H_- \in L_-(A)$\,
with \,$\|H_\pm(e_\pm)\|=1$\, and \,$\g{H_+(e_+)}{H_-(e_-)}=0$\,, and a canonical basis \,$(i,j,k)$\, of \,$\IM(\HH)$\,
so that
$$ \tfrac{1}{\sqrt{10}}\,(3H_+ + H_-) \;,\; \sqrt{\tfrac35}\,M_{1,1}^{(e_+,e_-)} + \sqrt{\tfrac{2}{5}}\,i\,H_- \;,\;
\underbrace{\sqrt{\tfrac35}\,M_{j,1}^{(e_+,e_-)} + \sqrt{\tfrac{2}{5}}\,k\,H_- }_{\text{only for \,$\ell=3$\,}}  $$
is an orthonormal basis of \,$\liem'$\,; here \,$M_{c,1}^{(e_+,e_-)}$\, is defined by \eqref{eq:curv:Mdef}.

\item \,$\boldsymbol{(\PP,\vi=\arctan(\tfrac12),\tau)}$\,, where \,$\tau$\, is the name of a \,$\HP$-type with \,$\dim(\tau) = 1$\, and
\,$w(\tau) \leq n-1$\,. 

There exist \,$A \in \A$\,, \,$J \in \frakJ_A$\,, a basis \,$(e_+,e_-)$\, adapted to \,$(A,J)$\,, \,$H_- \in \Sph(L_-(A))$\,, a \,$w(\tau)$-dimensional
real subspace \,$Q\subset \HH$\, with \,$1 \in Q$\,, a totally real, \,$w(\tau)$-dimensional subspace \,$U \subset L_+(A)$\, with
\,$U \perp J(H_-)$\,, and an \,$\R$-linear isometry \,$\Theta: Q \to U$\, so that
$$ \liem' = \Menge{q\,H_- + 2\,\Theta(q)}{q \in Q} $$
holds.

\item \,$\boldsymbol{(\PP,\vi=\arctan(\tfrac12),(\R,2))}$\,, if \,$n \geq 3$\,.

There exist \,$A \in \A$\,, \,$J \in \frakJ_A$\,, a basis \,$(e_+,e_-)$\, adapted to \,$(A,J)$\,, \,$H_+ \in L_+(A)$\, and \,$H_-,w_0 \in L_-(A)$\,,
where \,$H_+(e_+), H_-(e_-), w_0(e_-)$\, are pairwise \,$\HH$-orthogonal, so that 
$$ \tfrac{1}{\sqrt{5}}\,(2\,H_++H_-) \;,\; \sqrt{\tfrac25} M_{1,1}^{(e_+,e_-)} + \sqrt{\tfrac35} w_0 $$
is an orthonormal basis of \,$\liem'$\,. 

\item \,$\boldsymbol{(\PP,\vi=\arctan(\tfrac12),(\C,2))}$\,, if \,$n \geq 4$\,.

There exists a Lie triple system \,$\liem''$\, of \,$\liem$\, of type \,$(\PP,\vi=\arctan(\tfrac12),(\C,1))$\,
 -- with the data mentioned in the description of that type -- and \,$w_0 \in \Sph(L_-(A))$\, with \,$H_-$\,, \,$w_0$\, and \,$J(U)$\, pairwise \,$\HH$-orthogonal,
so that if we put \,$H_+ := \Theta(1) \in \Sph(L_+(A))$\, and let \,$i \in \Sph(\IM(\HH))$\, be so that \,$Q = \R \oplus \R\,i$\, holds, 
we have \,$\liem' = \liem'' \oplus \liem''_\perp$\,, where \,$\liem''_\perp$\, is the \,$\R$-linear space spanned by the orthonormal system
$$ \sqrt{\tfrac25}\,M_{1,1}^{(e_+,e_-)} - \tfrac{1}{\sqrt{5}}\,J(\Theta(i))\,i + \sqrt{\tfrac25} \,w_0 \;,\;
\sqrt{\tfrac25}\,M_{i,1}^{(e_+,e_-)} - \tfrac{1}{\sqrt{5}}\,J(\Theta(i)) + \sqrt{\tfrac25} \,w_0\,i \; . $$

\item \,$\boldsymbol{(\PP,\vi=\arctan(\tfrac12),(\HH,2))}$\,, if \,$n \geq 5$\,.

There exists a Lie triple system \,$\liem''$\, of \,$\liem$\, of type \,$(\PP,\vi=\arctan(\tfrac12),(\HH,1))$\,
 -- with the data mentioned in the description of that type -- so that with \,$H_+ := \Theta(1) \in \Sph(L_+(A))$\, we have \,$\liem' = \liem'' \oplus \liem''_\perp$\,,
where
$$ \liem''_\perp = \Mengegr{ \; M_{c,1}^{(e_+,e_-)} - \tfrac{1}{\sqrt{2}}\,\bigr(\,J(\Theta(i))\,\overline{c}\,i 
+ J(\Theta(j))\,\overline{c}\,j + J(\Theta(k))\,\overline{c}\,k \,\bigr) \; }{\;c \in \HH\;} $$
and \,$(i,j,k)$\, is any canonical basis of \,$\IM(\HH)$\,. 

\item \,$\boldsymbol{(\PP,\vi=\tfrac\pi4,\tau)}$\,, where \,$\tau$\, is the name of a \,$\HP$-type with \,$\dim(\tau) \leq \tfrac{n}{2}$\,.

There exists \,$A \in \A$\,, \,$J \in \frakJ_A$\,, two subspaces \,$W_1,W_2$\, of \,$L_+(A)$\, of \,$\HP$-type \,$\tau$\, (if \,$\tau = (\C,\ell)$\,: with respect
to the same \,$i \in \Sph(\IM(\HH))$\,) with \,$W_1 \perp W_2$\,, and an
anti-linear isometry \,$\Theta: W_1 \to W_2$\, so that
$$ \liem' = \Menge{x + J(\Theta(x))}{x \in W_1} $$
holds.

\item \,$\boldsymbol{(\Sph^5,\vi=\tfrac\pi4)}$\,.

There exists \,$A \in \A$\,, a basis \,$(e_+,e_-)$\, adapted to \,$A$\,, and \,$H_\pm \in \Sph(L_\pm(A))$\, with \,$\g{H_+(e_+)}{H_-(e_-)}=0$\, so that 
\,$\liem' = \R\,(H_+ + H_-) \oplus \Menge{M_{c,1}^{(e_+,e_-)}}{c \in \HH}$\, holds.

\item \,$\boldsymbol{(\Gtwo,\tau)}$\,, where \,$\tau$\, is the name of a \,$\HP$-type with \,$\dim(\tau) \leq n$\, and \,$\tau \neq (\Sph^3)$\,.

There exist \,$A \in \A$\,, \,$J \in \frakJ_A$\, and a subspace \,$U \subset L_+(A)$\, of \,$\HP$-type \,$\tau$\, so that \,$\liem' = U \operp J(U)$\, holds.

\item \,$\boldsymbol{(\PtP,\tau_1,\tau_2)}$\,, where \,$\tau_1$\, and \,$\tau_2$\, are names of \,$\HP$-types with \,$\dim(\tau_1) + \dim(\tau_2) \leq n$\,. 

There exist \,$A \in \A$\,, \,$J \in \frakJ_A$\, and subspaces \,$U_1,U_2 \subset L_+(A)$\, of \,$\HP$-type \,$\tau_1$\, resp.~\,$\tau_2$\, with 
\,$U_1 \perp U_2$\,, so that we have \,$\liem' = U_1 \operp J(U_2)$\,. 

\item \,$\boldsymbol{(\StSsqrt,\ell)}$\,, where \,$2 \leq \ell \leq 5$\, holds. 

There exists \,$A \in \A$\,, a basis \,$(e_+,e_-)$\, adapted to \,$A$\,, \,$H_\pm \in \Sph(L_\pm(A))$\, with \,$\g{H_+(e_+)}{H_-(e_-)}=0$\,
and an \,$\R$-linear subspace \,$C \subset \HH$\, of real dimension \,$\ell-1$\,, so that
\,$\liem' = \R\,H_+ \oplus \R\,H_- \oplus \Menge{M_{c,1}^{(e_+,e_-)}}{c \in C}$\, holds.

\item \,$\boldsymbol{(\Sptwo)}$\,. 

There exists an \,$\HH$-linear isometry \,$\Phi: V' \to V$\, so that \,$\liem' = \Menge{\Phi\circ X}{X \in \liesp(V')}$\, holds.

\item \,$\boldsymbol{(\Qthree)}$\,.

There exists a Lie triple system \,$\wh{\liem}' = U\oplus J(U)$\, of type \,$(\Gtwo,(\C,2))$\, (where \,$U$\, is totally complex with respect to \,$i \in \IM(\HH)$\,),
and a vector \,$v \in \Sph(\wh{\liem}')$\, with \,$\vi(v)=\tfrac\pi4$\,, so that \,$\liem'$\, is the \,$\R$-ortho-complement of \,$\R\,v \oplus \R\,v\,i$\,
in \,$\wh{\liem}'$\,.

\end{itemize}

We call the full name \,$(\Geo,\vi=t)$\,, \,$(\PP,\vi=0,\tau)$\, etc.~corresponding to a Lie triple system the \emph{type} of that Lie triple system.
If we identify 
\begin{gather*}
(\PtP,\tau_1,\tau_2) \cong (\PtP,\tau_2,\tau_1) \\
(\Geo,\vi=0) \cong (\PP,\vi=0,(\R,1)) \\
(\Geo,\vi=\arctan(\tfrac12)) \cong (\PP,\vi=\arctan(\tfrac12),(\R,1)) \\
(\Geo,\vi=\tfrac\pi4) \cong (\PP,\vi=\tfrac\pi4,(\R,1)) \;, 
\end{gather*}
then no Lie triple system is of more than one type, and two Lie triple systems are of the same type
if and only if they are congruent under the isotropy action of \,$\Sp(2)\times \Sp(n)$\, acting on \,$\liem$\,.

Moreover, the various types of Lie triple systems have the following properties:
\begin{center}
\begin{tabular}{|c|c|c|c|}
\hline
type of \,$\liem'$\, & \,$\dim_{\R}(\liem')$\, & \,$\rk(\liem')$\, & \,$\liem'$\, maximal \\
\hline
\,$(\Geo,\vi=t)$\, & $1$ & $1$ & no \\
\,$(\PP,\vi=0,\tau)$\, & \,$w(\tau)\,\dim(\tau)$\, & $1$ & for \,$\tau=(\HH,n)$\, \\
\,$(\Sph,\vi=\arctan(\tfrac13),\ell)$\, & \,$\ell$\, & $1$ & no \\
\,$(\PP,\vi=\arctan(\tfrac12),\tau)$\, & \,$w(\tau)\,\dim(\tau)$\, & $1$ & \tiny{if \,$n=4$\,: for \,$\tau=(\Sph^3)$\,; if \,$n=5$\,: for \,$\tau = (\HH,2)$\,} \\
\,$(\PP,\vi=\tfrac\pi4,\tau)$\, & \,$w(\tau)\,\dim(\tau)$\, & $1$ & no \\
\,$(\Sph^5,\vi=\tfrac\pi4)$\, & \,$5$\, & $1$ & no \\
\,$(\Gtwo,(\K,1))$\, & $2\,\dim_{\R}(\K)$ & $1$ & if \,$n=2$\,: for \,$\K=\HH$\, \\
\hline
\,$(\Gtwo,\tau)$\,, \,$\dim(\tau) \geq 2$\, & $2\,w(\tau)\, \dim(\tau)$ & $2$ &
        for \,$\tau=(\HH,n-1)$\, or \,$\tau=(\C,n)$\, \\
\,$(\PtP,\tau_1,\tau_2)$\, & {\footnotesize \,$w(\tau_1)\,\dim(\tau_1) + w(\tau_2)\,\dim(\tau_2)$\, } & $2$ & 
        for \,$\tau_\nu = (\HH,\ell_\nu)$\, with \,$\ell_1 + \ell_2 = n$\, \\
\,$(\Sph^1 \times \Sph^5,\ell)$\, & $1+\ell$ & $2$ & if \,$n=2$\,: for \,$\ell=5$\, \\
\,$(\Sptwo)$\, & $10$ & $2$ & if \,$n=2$\, \\
\,$(\Qthree)$\, & $6$ & $2$ & no \\
\hline
\end{tabular}
\end{center}
\end{Theorem}

\bigskip

\begin{Remark}
\label{R:claH:missing}
The following types of totally geodesic submanifolds are missing from the entry for \,$G_2(\HH^{n+2})$\, in the claimed classification of maximal totally
geodesic submanifolds of rank \,$2$\, Riemannian symmetric spaces in Table~VIII of \cite{Chen/Nagano:totges2-1978}:
\begin{enumerate}
\item
\,$(\Sph^1 \times \Sph^5,5)$\, (isometric to \,$(\Sph^5_{r=1/\sqrt{2}} \times \Sph^1_{r=1/\sqrt{2}})/\{\pm \id\}$\,, maximal in \,$G_2(\HH^4)$\,) 
\item 
\,$(\PP,\vi=\arctan(\tfrac12),(\Sph^3))$\, (isometric to \,$\Sph^3$\,, maximal in \,$G_2(\HH^6)$\,),
\item
\,$(\PP,\vi=\arctan(\tfrac12),(\HH,2))$\, (isometric to \,$\HP^2$\,, maximal in \,$G_2(\HH^7)$\,).
\end{enumerate}
While the totally geodesic submanifolds of type \,$(\Sph^1 \times \Sph^5,5)$\, are reflective in \,$G_2(\HH^4)$\, (the complementary type is \,$(\Sp_2)$\,)
and therefore can, for example, be found in the classification of reflective submanifolds of Riemannian symmetric spaces due to \textsc{Leung},
see \cite{Leung:reflective-1979}, the submanifolds of type \,$(\PP,\vi=\arctan(\tfrac12),\tau)$\, have, to my knowledge, never been described before.

Moreover, the existence of totally geodesic submanifolds of type \,$(\PP,\vi=\arctan(\tfrac13),\ell)$\, 
(isometric to \,$\Sph^\ell$\,), while not maximal in any \,$G_2(\HH^{n+2})$\, (rather, the type \,$(\PP,\vi=\arctan(\tfrac13),3)$\, is maximal in \,$\Sp(2)$\,), 
can not be deduced from Table~VIII of \cite{Chen/Nagano:totges2-1978} either. Also this type of totally geodesic submanifold has, as far as I know,
not been described before.
\end{Remark}

\begin{Remark}
\label{R:claH:pi4alternative}
There is an alternative description of the Lie triple systems of types \,$(\PP,\vi=\tfrac\pi4,(\C,\ell))$\, and \,$(\PP,\vi=\tfrac\pi4,(\HH,\ell))$\,:

\,$\liem' \subset \liem$\, is of type \,$(\PP,\vi=\tfrac\pi4,(\C,\ell))$\, (resp.~of type \,$(\PP,\vi=\tfrac\pi4,(\HH,\ell))$\,)
if and only if there exists \,$A \in \A$\,, \,$J \in \frakJ_A$\,, a totally real (resp.~totally complex), \,$2\ell$-dimensional subspace \,$W \subset L_+(A)$\,
and an orthogonal (resp.~anti-linear and orthogonal) map \,$\Xi: W \to W$\, with \,$\Xi^2 = -\id_W$\, so that
$$ \liem' = \Menge{x + J(\Xi(x))}{x \in W} $$
holds.
\end{Remark}

\bigskip

The remainder of the present section in concerned with the proof of Theorem~\ref{T:claH:claH}.

It is straightforward to see that the spaces mentioned in the theorem are indeed Lie triple systems, either by explicit calculations via Equation~\eqref{eq:curv:R}
or by checking that the embeddings described in Section~\ref{Se:Hembed} are indeed totally geodesic and correspond to spaces of the mentioned types. It is also easily
seen that two spaces of the same type are congruent under the isotropy action, and that the dimensions and ranks given in the table are correct.

To show that no two Lie triple systems of different type are congruent, it mostly suffices to note that two Lie triple systems for which the 
corresponding totally geodesic submanifolds have different isometry types (as given in Section~\ref{Se:Hembed}) cannot be congruent. This only leaves
the distinction between \,$(\Gtwo,(\K,1))$\, and \,$(\PP,\vi=0,(\K,2))$\, (for \,$\K \in \{\R,\C,\HH\}$\,). To show that these types indeed are not congruent,
we consider the normalizer group of Lie triple systems of the types involved under the action of the \,$\Sp(V')$-factor of the isotropy action on \,$\liem$\,.
The isomorphism type of these normalizer groups in dependence on the type \,$\tau$\, of the Lie triple system \,$\liem'$\, is given in the following table:

\smallskip
{
\scriptsize
\begin{center}
\begin{tabular}{|c||c|c||c|c||c|c|}
\hline
\,$\tau$\, & \,$(\Gtwo,(\R,1))$\, & \,$(\PP,\vi=0,(\R,2))$\, & \,$(\Gtwo,(\C,1))$\, & \,$(\PP,\vi=0,(\C,2))$\, & \,$(\Gtwo,(\HH,1))$\, & \,$(\PP,\vi=0,(\HH,2))$\, \\
\hline
normalizer & \,$\SO(2)$\, & \,$\Og(1)\times\Og(1)$\, & \,$\SU(2)$\, & \,$\Ug(1) \times \Ug(1)$\, & \,$\Sp(2)$\, & \,$\Sp(1)\times \Sp(1)$\, \\
\hline
\end{tabular} \; . 
\end{center}

}
\smallskip

Because the normalizers corresponding to \,$(\Gtwo,(\K,1))$\, and \,$(\PP,\vi=0,(\K,2))$\, are of different isomorphy type in each case, 
the Lie triple systems of these types cannot be congruent under the isotropy action.

For the data in the table on the maximality of the various Lie triple systems, we presume that the list of Lie triple systems given in the theorem
is complete. We first show why the Lie triple systems claimed to be maximal indeed are. For the Lie triple systems of rank \,$1$\,, note that 
such a Lie triple system can only be contained in another Lie triple system of rank \,$1$\, with the same characteristic angle \,$\vi$\,
(see Section~\ref{Se:orbits}), or in one of rank \,$2$\,. 
Using this fact, we see by inspection of the various types that Lie triple systems of the type \,$(\PP,\vi=0,(\HH,n))$\, are maximal, and that \,$(\Gtwo,(\HH,1))$\,
is maximal if \,$n=2$\, holds. 
Also, if \,$n=5$\, holds,
then \,$(\PP,\vi=\arctan(\tfrac12),(\HH,2))$\, is maximal: Assume to the contrary that there exists a Lie triple system \,$\liem''$\, with \,$\liem' \subsetneq \liem''
\subsetneq \liem$\,. Then \,$\liem''$\, must be of rank \,$2$\,, and the root system of \,$\liem''$\, must contain \,$\lambda_2$\, and \,$\lambda_4$\, with 
a multiplicity of at least \,$4$\,, also \,$\lambda_1$\, and \,$2\lambda_2$\, with a multiplicity of at least \,$3$\,. This shows that \,$\liem''$\, is of type
\,$(\Gtwo,(\HH,5))$\, and therefore equals \,$\liem$\, in contradiction to our assumption. A similar consideration shows that \,$(\PP,\vi=\arctan(\tfrac12),(\Sph^3))$\,
is maximal if \,$n=4$\, holds. 

For the Lie triple systems \,$\liem'$\, of rank \,$2$\,, note that such a Lie triple system can only be contained in another Lie triple system \,$\liem''$\, of rank \,$2$\,,
and if \,$\Delta(\liem',\liea)$\, and \,$\Delta(\liem'',\liea)$\, are the root systems of these two Lie triple systems with respect to a Cartan subalgebra 
\,$\liea \subset \liem'$\,, we have \,$\Delta(\liem',\liea) \subset \Delta(\liem'',\liea)$\, and \,$n_\alpha(\liem'') \geq n_\alpha(\liem')$\, for every
\,$\alpha \in \Delta(\liem',\liea)$\,. Using this fact, we see that Lie triple systems of the types \,$(\Gtwo,(\HH,n-1))$\, and \,$(\Gtwo,(\C,n))$\, are maximal, and
also that if \,$n=2$\, holds, then Lie triple systems of type \,$(\Sph^1 \times \Sph^5,5)$\, and \,$(\Sp_2)$\, are maximal. Finally, one sees by a
consideration of the explicit description of the types of rank \,$2$\, that Lie triple systems of type \,$(\PtP, (\HH,\ell_1), (\HH,\ell_2))$\, with
\,$\ell_1+\ell_2=n$\, are maximal.

That no Lie triple systems are maximal besides those mentioned above follows from the facts in the following table:
\begin{center}
\begin{longtable}{|c|c|}
\hline
Every Lie triple system of type ... & is contained in a Lie triple system of type ... . \\
\hline
\endhead
\hline
\endfoot
\,$(\Geo,\vi=t)$\, & \,$(\PtP,(\R,1),(\R,1))$\, \\
\,$(\PP,\vi=0,\tau)$\, with \,$\dim(\tau)\leq n-1$\, & \,$(\PP,\vi=0,(\HH,n))$\, \\
\,$(\PP,\vi=0,(\K,n))$\, with \,$\K \in \{\R,\C\}$\, & \,$(\PP,\vi=0,(\HH,n))$\, \\
\,$(\Sph,\vi=\arctan(\tfrac13),\ell)$\, & \,$(\Sp_2)$\, \\
\,$(\PP,\vi=\arctan(\tfrac12),(\K,\ell))$\, with \,$\K \in \{\R,\C\}$\, & \,$(\Gtwo,(\K,\ell+\dim_{\R}(\K)))$\, \\
\,$(\PP,\vi=\arctan(\tfrac12),(\Sph^3))$\, if \,$n \geq 5$\, & \,$(\PP,\vi=\arctan(\tfrac12),(\HH,1))$\, \\
\,$(\PP,\vi=\arctan(\tfrac12),(\HH,1))$\, if \,$n \geq 6$\, & \,$(\Gtwo,(\HH,5))$\, \\
\,$(\PP,\vi=\tfrac\pi4,\tau)$\, & \,$(\PtP,\tau,\tau)$\, \\
\,$(\Sph^5,\vi=\tfrac\pi4)$\, & \,$(\Sph^1 \times \Sph^5,5)$\, \\
\,$(\Gtwo,(\K,\ell))$\, with \,$\K \in \{\R,\C\}$\, and \,$\ell \leq n-1$\, & \,$(\Gtwo,(\K,n))$\, \\
\,$(\Gtwo,(\HH,\ell))$\, with \,$\ell \leq n-2$\, & \,$(\Gtwo,(\HH,n-1))$\, \\
\,$(\Gtwo,(\R,n))$\, & \,$(\Gtwo,(\C,n))$\, \\
\,$(\PtP,\tau_1,\tau_2)$\, with \,$\tau_\nu$\, not both \,$(\HH,\ell_\nu)$\, & \,$(\PtP,(\HH,\dim(\tau_1)),(\HH,\dim(\tau_2)))$\, \\
\,$(\PtP,(\HH,\ell_1),(\HH,\ell_2))$\, with \,$\ell_1+\ell_2 \leq n-1$\, & \,$(\PtP,(\HH,\ell_1),(\HH,n-\ell_1))$\, \\
\,$(\Sph^1 \times \Sph^5,\ell)$\, with \,$\ell \leq 4$\, & \,$(\Sph^1 \times \Sph^5,5)$\, \\
\,$(\Sph^1 \times \Sph^5,5)$\, if \,$n\geq 3$\, & \,$(\Gtwo,(\HH,2))$\, \\
\,$(\Sp_2)$\, if \,$n \geq 3$\, & \,$(\Gtwo,(\HH,2))$\, \\
\,$(Q_3)$\, & \,$(\Gtwo,(\C,2))$\, \\
\hline
\end{longtable}
\end{center}

We now focus on the main problem, namely the proof that the list of Lie triple systems given in the theorem is indeed complete.

As was emphasized before, for a classification of the Lie triple systems, it is not sufficient to know the root system
(with multiplicities) of the symmetric space under investigation; rather one has to know the structure of the curvature
tensor in all three variables in order to understand the actual transformations it induces between the various root spaces.

The required description of the curvature tensor \,$R$\, is in our situation, i.e.~for the symmetric space \,$G_2(V' \oplus V)$\,, essentially provided by
Equation~\eqref{eq:curv:R}. However, the more explicit description of \,$R$\, we now give, and which easily follows from Equation~\eqref{eq:curv:R},
is often more handy and will be used frequently throughout the classification.
For this, we let \,$A \in \A$\,, \,$u_\pm,v_\pm \in L_\pm(A)$\, and \,$w \in \liem$\, be given. Then we describe the linear map 
\,$R(u_\pm,v_\pm)w \in \liem \cong L(V',V)$\, by stating the image of an arbitrarily fixed basis \,$(e_+,e_-)$\, adapted to \,$A$\,:
\begin{align}
\label{eq:claH:R11}
R(u_+,v_+)w & = 
\begin{cases}
e_+ & \mapsto u_+(e_+)\,\g{v_+(e_+)}{w(e_+)} - v_+(e_+)\,\g{u_+(e_+)}{w(e_+)} + 2\,w(e_+)\,\IM(\g{v_+(e_+)}{u_+(e_+)}) \\
e_- & \mapsto u_+(e_+)\,\g{v_+(e_+)}{w(e_-)} - v_+(e_+)\,\g{u_+(e_+)}{w(e_-)} 
\end{cases} \\
\label{eq:claH:R12}
R(u_+,v_-)w & = 
\begin{cases}
e_+ & \mapsto w(e_-)\,\g{v_-(e_-)}{u_+(e_+)} \\
e_- & \mapsto -w(e_+) \, \g{u_+(e_+)}{v_-(e_-)}
\end{cases}  \\
\label{eq:claH:R22}
R(u_-,v_-)w & = 
\begin{cases}
e_+ & \mapsto u_-(e_-)\,\g{v_-(e_-)}{w(e_+)} - v_-(e_-)\,\g{u_-(e_-)}{w(e_+)} \\
e_- & \mapsto u_-(e_-)\,\g{v_-(e_-)}{w(e_-)} - v_-(e_-)\,\g{u_-(e_-)}{w(e_-)} + 2\,w(e_-)\,\IM(\g{v_-(e_-)}{u_-(e_-)}) \\
\end{cases}
\end{align}

Because the Riemannian symmetric space \,$G_2(V' \oplus V)$\, is of rank \,$2$\,, any Lie triple system \,$\liem'$\, of \,$\liem$\, has either rank \,$2$\,
or rank \,$1$\,. We will handle the classification for the two different ranks separately, in the following two subsections.

\subsection{The case of rank 2}
\label{SSe:claH:rk2}

In this subsection we let a Lie triple system \,$\liem'$\, of \,$\liem$\, with  \,$\rk(\liem')=2=\rk(\liem)$\, be given. Let us fix a Cartan subalgebra \,$\liea$\,
of \,$\liem'$\,, then \,$\liea$\, also is a Cartan subalgebra of \,$\liem$\,. 
Then we have the root systems \,$\Delta := \Delta(\liem,\liea)$\, and \,$\Delta' := \Delta(\liem',\liea')$\, of \,$\liem$\, resp.~\,$\liem'$\, with regard 
to \,$\liea$\,, and consequently the root space decompositions
\begin{equation}
\label{eq:claH:rk2:decomps}
\liem = \liea \oplus \bigoplus_{\lambda \in \Delta_+} \liem_\lambda \qmq{and} \liem' = \liea' \oplus \bigoplus_{\alpha \in \Delta_+'} \liem_\alpha' \; .
\end{equation}
Proposition~\ref{P:cla:subroots:subroots-neu}(b) shows that we have
\begin{equation}
\label{eq:claH:rk2:m'subspaces}
\liea' = \liea,\quad \Delta' \subset \Delta,\quad \forall \alpha \in \Delta' : \liem_\alpha' = \liem_\alpha \cap \liem' \; . 
\end{equation}

We now write down the action of the curvature tensor between certain root spaces explicitly; these formulas 
(whose derivation again requires the understanding of the full curvature tensor) will play
an important role in the classification: Because \,$\liem'$\, is a Lie triple system, \,$u,v,w \in \liem'$\, implies \,$R(u,v)w \in \liem'$\,;
the formulas therefore permit to derive from the presence of certain vectors in \,$\liem'$\, the presence of certain other vectors in \,$\liem'$\,. 

In relation to this, we note that 
by Proposition~\ref{P:curv:cartan}, there exist \,$A \in \A$\, and an orthonormal basis \,$(H_+,H_-)$\, of \,$\liea$\, with \,$H_\pm \in L_\pm(A)$\,
and \,$H_+(V_+'(A)) \perp H_-(V_-'(A))$\,.

\begin{Lemma}
\label{L:claH:Rroot}
Let \,$J \in \frakJ_A$\, and \,$(e_+,e_-)$\, be any basis adapted to \,$(A,J)$\,. We define \,$M_{c,\eps} := M_{c,\eps}^{(e_+,e_-)}$\, as in
Equation~\eqref{eq:curv:Mdef}.
\begin{enumerate}
\item
\textbf{Action on \,$\boldsymbol{\liem_{\lambda_1}}$\, and \,$\boldsymbol{\liem_{\lambda_2}}$\,.}
Let \,$u_+,v_+ \in \liem_{\lambda_1}$\, and \,$u_-,v_- \in \liem_{\lambda_2}$\, be given. Then we have
$$ R(H_+,u_+)v_+ = \g{u_+(e_+)}{v_+(e_+)}\,H_+ \qmq{and} R(H_-,u_-)v_- = \g{u_-(e_-)}{v_-(e_-)}\,H_- \; . $$
\item
\textbf{From \,$\boldsymbol{\liem_{\lambda_1}}$\, to \,$\boldsymbol{\liem_{\lambda_2}}$\, and vice versa, 
via \,$\boldsymbol{\liem_{\lambda_3}}$\, or \,$\boldsymbol{\liem_{\lambda_4}}$\,.} 
Let \,$c \in \HH$\,, \,$\eps \in \{\pm 1\}$\,, \,$w_+ \in \liem_{\lambda_1}$\, and \,$w_- \in \liem_{\lambda_2}$\, be given. 
Then we have
$$ R(H_-,M_{c,\eps})w_+ = \tfrac{1}{\sqrt{2}} \, J(w_+)\cdot \overline{c} \qmq{and} R(H_+,M_{c,\eps})w_- = -\tfrac{\eps}{\sqrt{2}}\, J(w_-)\cdot c \; . $$
\item
\textbf{From \,$\boldsymbol{\liem_{\lambda_3}}$\, to \,$\boldsymbol{\liem_{\lambda_4}}$\, and vice versa, 
via \,$\boldsymbol{\liem_{\lambda_1}}$\, or \,$\boldsymbol{\liem_{\lambda_2}}$\,.} 
Let \,$u_+,v_+ \in \liem_{\lambda_1}$\,, \,$u_-,v_- \in \liem_{\lambda_2}$\,, \,$c \in \HH$\, and \,$\eps \in \{\pm 1\}$\, be given.
Put \,$d_+ := \IM(\g{v_+(e_+)}{u_+(e_+)})$\, and \,$d_- := \IM(\g{v_-(e_-)}{u_-(e_-)})$\,.
Then we have
$$ R(u_+,v_+)M_{c,\eps} = M_{cd_+,-\eps} + M_{cd_+,\eps} \qmq{and} 
R(u_-,v_-)M_{c,\eps} = M_{d_-c,-\eps} - M_{d_-c,\eps} \; . $$
\item
\textbf{From \,$\boldsymbol{\liem_{\lambda_3}}$\, to \,$\boldsymbol{\liem_{\lambda_4}}$\, and vice versa,
via \,$\boldsymbol{\liem_{2\lambda_1}}$\, or \,$\boldsymbol{\liem_{2\lambda_2}}$\,.} 
Let \,$d \in \IM(\HH)$\,, \,$c \in \HH$\, and \,$\eps \in \{\pm 1\}$\, be given. Then we have
$$ R(H_+,d\,H_+)M_{c,\eps} = M_{-2cd,-\eps} \qmq{and} R(H_-,d\,H_-)M_{c,\eps} = M_{-2dc,-\eps} \; . $$
\item
\textbf{From \,$\boldsymbol{\liem_{2\lambda_1}}$\, to \,$\boldsymbol{\liem_{2\lambda_2}}$\, and vice versa,
via \,$\boldsymbol{\liem_{2\lambda_1} \oplus \liem_{2\lambda_2}}$\,.} 
Let \,$c,\wt{c} \in \HH$\,, \,$d \in \IM(\HH)$\, and \,$\eps \in \{\pm 1\}$\, be given. Then we have
\begin{align*}
R(dH_+,M_{c,\eps})M_{\wt{c},-\eps} & = \RE(d\overline{c}\wt{c})\,H_+ + \eps\,\IM(cd\overline{\wt{c}})\,H_- \\
\qmq{and} R(dH_-,M_{c,\eps})M_{\wt{c},-\eps} & = -\eps\,\IM(\overline{c}d\wt{c})\,H_+ + \RE(\wt{c}\overline{c}d)\,H_- \; .  
\end{align*}
\item
\textbf{From \,$\boldsymbol{\liem_{\lambda_1} \oplus \liem_{\lambda_2}}$\, to \,$\boldsymbol{\liem_{\lambda_3} \oplus \liem_{\lambda_4}}$\,.}
Let \,$u_+ \in \liem_{\lambda_1}$\,, \,$v_- \in \liem_{\lambda_2}$\, and \,$\eps \in \{\pm 1\}$\, be given. 
Then we have
$$ R(u_+,v_-)(-\eps\,H_++H_-) = \sqrt{2}\,M_{\g{v_-(e_-)}{u_+(e_+)}\,,\,\eps} \; . $$
\item
\textbf{From \,$\boldsymbol{\liem_{\lambda_3} \oplus \liem_{\lambda_4}}$\, to \,$\boldsymbol{\liem_{2\lambda_1} \oplus \liem_{2\lambda_2}}$\,.}
Let \,$c,\wt{c} \in \HH$\, be given. Then we have
$$ R(H_+,M_{c,1})M_{\wt{c},-1} = H_+ \cdot \IM(\overline{c}\,\wt{c}) - H_- \cdot \IM(\wt{c}\,\overline{c}) \; . $$
\end{enumerate}
\end{Lemma}

\beweis
All these formulas are easily derived from Equations~\eqref{eq:claH:R11}, \eqref{eq:claH:R12} and \eqref{eq:claH:R22}.
\beweisende

We further prepare the classification by the following lemma, which exhibits how \,$\liem_{2\lambda_\nu}'$\, controls the structure of \,$\liem_{\lambda_\nu}'$\,
(for \,$\nu \in \{1,2\}$\,).

\begin{Lemma}
\label{L:claH:mstructure}
Let \,$\nu \in \{1,2\}$\, and suppose \,$\lambda_\nu \in \Delta'$\,; put \,$\pm := +$\, for \,$\nu=1$\,, \,$\pm := -$\, for \,$\nu=2$\,. 
Then \,$K := \Menge{c \in \HH}{H_\pm\,c \in \liem'}$\, is a sub-field of \,$\HH$\, with
\,$K \cdot \liem_{\lambda_\nu}' \subset \liem_{\lambda_\nu}'$\, and \,$K^{\perp,\HH} \cdot \liem_{\lambda_\nu}' \perp_{\R} \liem_{\lambda_\nu}'$\,.
Consequently, we have \,$n_{2\lambda_\nu}' \in \{0,1,3\}$\, and:
\begin{enumerate}
\item If \,$n_{2\lambda_\nu}'=0$\, holds, then \,$\liem_{\lambda_\nu}'$\, is a totally real subspace of \,$\liem_{\lambda_\nu}$\,.
\item If \,$n_{2\lambda_\nu}'=1$\, holds, say \,$\liem_{2\lambda_\nu}' = \R\,H_\pm\,i$\, with \,$i \in \Sph(\IM(\HH))$\,,
then \,$\liem_{\lambda_\nu}'$\, is a totally complex subspace of \,$\liem_{\lambda_\nu}$\, with respect to \,$i$\,. 
\item If \,$n_{2\lambda_\nu}'=3$\, holds, then \,$\liem_{\lambda_\nu}'$\, is a quaternionic subspace of \,$\liem_{\lambda_\nu}$\,.
\end{enumerate}
\end{Lemma}

\beweis
We consider the case \,$\nu=1$\,; the case \,$\nu=2$\, is proved analogously. Let us fix a basis \,$(e_+,e_-)$\, adapted to \,$A$\, 
and let \,$L_+' := (L_+(A) \cap \liem')$\,; note that because of Equations~\eqref{eq:claH:rk2:decomps} and \eqref{eq:claH:rk2:m'subspaces},
we have the splitting 
\begin{equation}
\label{eq:claH:mstructure:L+split}
L_+' = K\,H_+ \oplus L_+'' \qmq{with} L_+'' := (\HH\,H_+)^\perp \cap L_+' \; . 
\end{equation}
For any \,$v \in L_+' \setminus \{0\}$\,, we consider the \,$\R$-linear subspace
$$ K_v := \Menge{ c\in \HH }{v\,c \in L_+'} $$
of \,$\HH$\,; clearly we have \,$K_{H_+} = K$\,. Below we will show that for any \,$v,w \in L_+' \setminus \{0\}$\,
\begin{equation}
\label{eq:claH:mstructure:Kv}
\g{v(e_+)}{w(e_+)}=0 \;\;\;\Longrightarrow\;\;\; K_v = K_w
\end{equation}
holds; because of the splitting \eqref{eq:claH:mstructure:L+split} of \,$L_+'$\, into two non-zero, orthogonal summands
it follows that  \,$K_v$\, does not depend on \,$v \in L_+' \setminus \{0\}$\,, and therefore we then have
\begin{equation}
\label{eq:claH:mstructure:Kinvar}
\forall v \in L_+'\setminus\{0\} \; : \; K_v = K \; . 
\end{equation}

For the proof of \eqref{eq:claH:mstructure:Kv}, we let \,$v,w \in L_+' \setminus \{0\}$\, with \,$\g{v(e_+)}{w(e_+)}=0$\,
and \,$c \in K_v$\, be given. Then we also have \,$v \, c \in \liem'$\, by the definition of \,$K_v$\,, and therefore, owing to
the fact that \,$\liem'$\, is a Lie triple system, \,$R(w,v)(v\,c) \in \liem'$\,. Via Equation~\eqref{eq:claH:R11} we calculate
\begin{equation*}
\liem' \ni R(w,v)(v\,c) = w\cdot \underbrace{\|v(e_+)\|^2}_{\neq 0}\,c \; ;
\end{equation*}
this shows that \,$w\,c \in \liem'$\, and therefore \,$c \in K_w$\, holds. This proves the inclusion \,$K_v \subset K_w$\,;
the opposite inclusion is shown in the same way.

Next we show that for any \,$v,w \in \liem_{\lambda_1}'$\,,
\begin{equation}
\label{eq:claH:mstructure:skp}
\g{v(e_+)}{w(e_+)} \in K
\end{equation}
holds. In fact, we have \,$H_+,v,w \in \liem'$\, and therefore \,$R(H_+,v)w \in \liem'$\,; using Equation~\eqref{eq:claH:R11} one calculates that
\,$ R(H_+,v)w = H_+ \cdot \g{v(e_+)}{w(e_+)} $\,
holds. Thus we have \,$H_+ \cdot \g{v(e_+)}{w(e_+)} \in \liem'$\, and therefore \eqref{eq:claH:mstructure:skp} holds.

It is a consequence of \eqref{eq:claH:mstructure:Kinvar} that the \,$\R$-linear subspace \,$K$\,
of \,$\HH$\, with \,$1 \in K$\, is closed under multiplication, and is therefore a sub-field of \,$\HH$\,. Hence we have either \,$K = \R$\,,
or \,$K = \spn\{1,i\} \cong \C$\, with some \,$i \in \Sph(\IM(\HH))$\,, or \,$K = \HH$\,. 

If \,$K=\R$\, holds, we have \,$n_{2\lambda_1}' = 0$\,, and \,$\liem_{\lambda_1}'$\, is totally real by \eqref{eq:claH:mstructure:skp}.

If \,$K=\spn\{1,i\}$\, holds with some \,$i \in \Sph(\IM(\HH))$\,, we have \,$n_{2\lambda_1}'=1$\,, and \,$\liem_{\lambda_1}'$\, is totally complex
with respect to \,$i$\, by \eqref{eq:claH:mstructure:Kinvar} and \eqref{eq:claH:mstructure:skp}.

If \,$K=\HH$\, holds, we have \,$n_{2\lambda_1}'=3$\,, and \,$\liem_{\lambda_1}'$\, is quaternionic by \eqref{eq:claH:mstructure:Kinvar}.
\beweisende

In the sequel, we divide the classification into the following cases, depending on the configuration of~\,$\Delta'$\,:
\begin{itemize}
\item[(1)] \,$\lambda_1,\lambda_2,\lambda_3,\lambda_4 \in \Delta'$\,,
\item[(2)] \,$\lambda_1,\lambda_2 \in \Delta'$\,, but \,$\lambda_3,\lambda_4 \not\in \Delta'$\,,
\item[(3)] either \,$\lambda_1$\, or \,$\lambda_2$\,, but not both, are elements of \,$\Delta'$\,,
\item[(4)] \,$\lambda_1,\lambda_2 \not\in \Delta'$\,, \,$\lambda_3,\lambda_4 \in \Delta'$\,,
\item[(5)] \,$\lambda_1,\lambda_2 \not\in \Delta'$\,, either \,$\lambda_3$\, or \,$\lambda_4$\,, but not both, are elements of \,$\Delta'$\,,
\item[(6)] \,$\lambda_1,\lambda_2, \lambda_3,\lambda_4 \not\in \Delta'$\,.
\end{itemize}
It follows from the invariance of the root system \,$\Delta'$\, under the action of its Weyl group, that whenever we have \,$\lambda_1 \in \Delta'$\,
or \,$\lambda_2\in \Delta'$\,, the presence of either of the two roots \,$\lambda_3$\, and \,$\lambda_4$\, in \,$\Delta'$\, implies the presence of the other.
Therefore these six cases exhaust all possibilities for \,$\Delta'$\,. 
We will now handle these cases separately.

\bigskip

\textbf{Case (1).} We suppose \,$\lambda_1,\lambda_2,\lambda_3,\lambda_4 \in \Delta'$\,. In particular, we have \,$\liem_{\lambda_3}' \neq \{0\}$\,, and therefore
Equation~\eqref{eq:curv:Mtrafo} shows that there exists \,$J \in \frakJ_A$\, and a basis \,$(e_+,e_-)$\, adapted to \,$(A,J)$\, so that 
\,$M_{1,-1}^{(e_+,e_-)} \in \liem_{\lambda_3}'$\, holds. We show that the following relations hold between root spaces of \,$\liem'$\, corresponding to roots
from the same orbit of the action of the Weyl group of \,$\Delta'$\,:
\begin{align}
\label{eq:claH:rk2:m1m2}
\liem_{\lambda_2}' & = J(\liem_{\lambda_1}') \\
\label{eq:claH:rk2:m3m4}
\liem_{\lambda_4}' & = \Menge{M_{c,1}^{(e_+,e_-)}}{c \in \HH,\; M_{c,-1}^{(e_+,e_-)} \in \liem_{\lambda_3}'} \\
\label{eq:claH:rk2:m21m22}
\liem_{2\lambda_2} & = \Menge{d\,H_-}{d \in \IM(\HH),\; d\,H_+ \in \liem_{2\lambda_1}'} \; . 
\end{align}

\emph{For \eqref{eq:claH:rk2:m1m2}.} Let \,$w \in \liem_{\lambda_1}'$\, be given. We have \,$H_-,M_{1,-1}^{(e_+,e_-)},w \in \liem'$\,, and therefore
also \,$\liem' \ni R(H_-,M_{1,-1}^{(e_+,e_-)})w \overset{(*)}{=} \tfrac{1}{\sqrt{2}}\,J(w_+)$\, (the equals sign marked \,$(*)$\, follows from Lemma~\ref{L:claH:Rroot}(b));
in fact we have \,$J(w_+) \in \liem' \cap \liem_{\lambda_2} = \liem_{\lambda_2}'$\,. 
This shows the inclusion ``$\supset$'' in \eqref{eq:claH:rk2:m1m2}, and the converse inclusion is shown in the same way.

\emph{For \eqref{eq:claH:rk2:m3m4}.} Let \,$c \in \HH$\, with 
\,$M_{c,-1}^{(e_+,e_-)} \in \liem_{\lambda_3}'$\, be given, and let us fix \,$w \in \liem_{\lambda_1}'$\, with \,$\|w(e_+)\|=1$\,.
We have \,$H_-,M_{c,-1}^{(e_+,e_-)},(-w) \in \liem'$\, and therefore also \,$\liem' \ni R(H_-,M_{c,-1}^{(e_+,e_-)}(-w) \overset{(*)}{=} -\tfrac{1}{\sqrt{2}}\,J(w) \cdot \overline{c}$\,
($(*)$ again follows from Lemma~\ref{L:claH:Rroot}(b)), hence \,$v_- := -J(w) \cdot \overline{c} \in \liem_{\lambda_2} \cap \liem' = \liem_{\lambda_2}'$\,. 
With \,$u_+ := w \in \liem_{\lambda_1}'$\, we have \,$\g{v_-(e_-)}{u_+(e_+)} = \g{w(e_+)\,\overline{c}}{w(e_+)} = c$\,,
and therefore -- because of \,$u_+,v_-,H_++H_- \in \liem'$\, -- \,$\liem' \ni R(u_+,v_-)(H_++H_-) \overset{(*)}{=} \sqrt{2}\,M_{c,1}^{(e_+,e_-)}$\,
($(*)$ by Lemma~\ref{L:claH:Rroot}(f)), hence \,$M_{c,1}^{(e_+,e_-)} \in \liem_{\lambda_4} \cap \liem' = \liem_{\lambda_4}'$\,. 
This shows the inclusion ``$\supset$'' in \eqref{eq:claH:rk2:m3m4}; the converse inclusion is again shown analogously.

\emph{For \eqref{eq:claH:rk2:m21m22}.} We have \,$M_{1,-1}^{(e_+,e_-)} \in \liem_{\lambda_3}'$\, and therefore, by Equation~\eqref{eq:claH:rk2:m3m4}, 
also \,$M_{1,1}^{(e_+,e_-)} \in \liem_{\lambda_4}'$\,. Therefore Lemma~\ref{L:claH:Rroot}(e) shows that for any \,$d \in \IM(\HH)$\, with 
\,$d\,H_+ \in \liem_{2\lambda_1}'$\,, we also have \,$d\,H_- \in \liem_{2\lambda_2}'$\,, and vice versa. 

Let us now consider \,$K := \Menge{c \in \HH}{c\,H_+ \in \liem'} = \R \oplus \Menge{d \in \IM(\HH)}{d\,H_+ \in \liem_{2\lambda_1}'}$\, of \,$\HH$\,;
we saw in Lemma~\ref{L:claH:mstructure} that \,$K$\, is a sub-field of \,$\HH$\,, and how \,$K$\, controls the structure of \,$\liem_{\lambda_1}$\,. 
Because of Equation~\eqref{eq:claH:rk2:m21m22} we in fact have
\begin{equation}
\label{eq:claH:rk2:Hpm}
\R H_\pm \oplus \liem_{2\lambda_\nu} = K\cdot H_\pm \qmq{for \,$\nu \in \{1,2\}$\,}
\end{equation}
hence \,$K$\, also determines the structure of \,$\liem_{\lambda_2}$\, in the way described in Lemma~\ref{L:claH:mstructure}. 

We will now show that \,$K$\, moreover ``controls'' the root spaces \,$\liem_{\lambda_3}'$\, and \,$\liem_{\lambda_4}'$\,, more specifically, that we have
\begin{equation}
\label{eq:claH:rk2:m3m4'}
\liem_{\lambda_\nu}' = \Menge{M_{c,\mp 1}^{(e_+,e_-)}}{c \in K} \qmq{for \,$\nu \in \{3,4\}$\,.} 
\end{equation}
We prove this equation for \,$\nu=3$\,; for \,$\nu=4$\, the proof runs analogously. So let \,$c \in \HH$\, with \,$M_{c,-1}^{(e_+,e_-)} \in \liem_{\lambda_3}'$\, 
be given. Again we fix \,$u_+ \in \liem_{\lambda_1}'$\, with \,$\|u_+(e_+)\| = 1$\,. Then we have on one hand \,$J(u_+) \in \liem_{\lambda_2}'$\, by
Equation~\eqref{eq:claH:rk2:m1m2}, on the other hand \,$J(u_+)\cdot \overline{c} \in \liem_{\lambda_2}'$\, by Lemma~\ref{L:claH:Rroot}(b). Because
of \,$H_-, J(u_+)\cdot\overline{c}, J(u_+) \in \liem'$\, we conclude by Lemma~\ref{L:claH:Rroot}(a) that \,$H_-\,c \in \liem'$\,, and therefore by
Equation~\eqref{eq:claH:rk2:Hpm} \,$c \in K$\, holds.

Conversely, let \,$c \in K$\, be given. With \,$u_+ \in \liem_{\lambda_1}'$\, being as before, we then also have \,$u_+ \cdot c \in \liem_{\lambda_1}'$\,
by Lemma~\ref{L:claH:mstructure}, and \,$J(u_+) \in \liem_{\lambda_2}'$\, by Equation~\eqref{eq:claH:rk2:m1m2}. Because of \,$u_+\cdot c, J(u_+), H_++H_- \in \liem'$\,
it follows from Lemma~\ref{L:claH:Rroot}(f) that \,$M_{c,- 1}^{(e_+,e_-)} \in \liem_{\lambda_3}'$\, holds. This concludes the proof
of Equation~\eqref{eq:claH:rk2:m3m4'}.

Now we have
\begin{align*}
\liem' & = \underbrace{\liea \oplus \liem_{2\lambda_1}' \oplus \liem_{2\lambda_2}'}_{\underset{\eqref{eq:claH:rk2:Hpm}}{=} K\cdot H_+ \oplus K\cdot H_-}
\; \oplus \;  \underbrace{\liem_{\lambda_3}' \oplus \liem_{\lambda_4}'}_{\underset{\eqref{eq:claH:rk2:m3m4'}}{=} K\cdot J(H_+) \oplus K \cdot J(H_-)} 
\; \oplus \;  \liem_{\lambda_1}' \oplus \underbrace{\liem_{\lambda_2}'}_{\underset{\eqref{eq:claH:rk2:m1m2}}{=} J(\liem_{\lambda_1}')} \\
& = U \oplus J(U)
\end{align*}
with \,$U := K\cdot H_+ \oplus K\cdot J(H_-) \oplus \liem_{\lambda_1}' \subset L_+(A)$\,. 

By Lemma~\ref{L:claH:mstructure} \,$U$\, is a either a quaternionic subspace (if \,$n_{2\lambda_k}' = 3$\,), a totally complex subspace (if \,$n_{2\lambda_k}' = 1$\,), 
or a totally real subspace (if \,$n_{2\lambda_k}' = 0$\,) of \,$L_+(A)$\,. This shows that -- depending on which of these three cases holds -- \,$\liem'$\, 
is either of the type \,$(\Gtwo, (\HH,2+\tfrac{1}{4} n_{\lambda_1}'))$\,, or of the type \,$(\Gtwo, (\C,2+\tfrac{1}{2} n_{\lambda_1}'))$\, 
or of the type \,$(\Gtwo, (\R, 2+n_{\lambda_1}'))$\,.

\bigskip

\textbf{Case (2).} We now suppose \,$\lambda_1,\lambda_2 \in \Delta'$\,, but \,$\lambda_3,\lambda_4 \not\in \Delta'$\,. As before, Lemma~\ref{L:claH:mstructure}
shows that \,$\liem_{\lambda_1}'$\, and \,$\liem_{\lambda_2}'$\, are (individually) either quaternionic, totally complex, or totally real. More specifically,
if we put
$$ K_\pm := \Menge{c \in \HH}{c\,H_\pm \in \liem'} \;,$$
then \,$\liem_{\lambda_k}'$\, is invariant under multiplication with \,$K_\pm$\, and satisfies \,$(K_\pm)^{\perp,\HH} \cdot \liem_{\lambda_k}' \perp_{\R} \liem_{\lambda_k}'$\,.
However, because of the
absence of the roots \,$\lambda_3$\, and \,$\lambda_4$\, in \,$\Delta'$\,, there is no binding between \,$K_+$\, and \,$K_-$\, anymore.

From the absence of the roots \,$\lambda_3$\, and \,$\lambda_4$\,, a condition on the relative positions of \,$\liem_{\lambda_1}'$\, and \,$\liem_{\lambda_2}'$\,
can be derived: Let \,$u_+ \in \liem_{\lambda_1}'$\, and \,$v_- \in \liem_{\lambda_2}'$\, be given. We then have
\,$R(u_+,v_-)(H_++H_-) \in (\liem_{\lambda_3} \oplus \liem_{\lambda_4})\cap \liem' = \{0\}$\, and therefore \,$0 = R(u_+,v_-)(H_++H_-) = \sqrt{2} \, M_{\g{v_-(e_-)}{u_+(e_+)}\,,\,-1}^{(e_+,e_-)}$\,
(where the last equality sign follows from Lemma~\ref{L:claH:Rroot}(f)), hence \,$\g{v_-(e_-)}{u_+(e_+)} = 0$\,. Thus we have shown
$$ J(\liem_{\lambda_2}') \perp \liem_{\lambda_1}' \; . $$

Therefore we have
$$
\liem' \;=\; \underbrace{\liea \oplus \liem_{2\lambda_1}' \oplus \liem_{2\lambda_2}'}_{= K_+\cdot H_+ \oplus K_-\cdot H_-} 
\oplus \, \liem_{\lambda_1}' \oplus \liem_{\lambda_2}' 
\;=\; U_1 \oplus J(U_2)
$$
with \,$U_1 := K_+ \cdot H_+ \oplus \liem_{\lambda_1}' \subset L_+(A)$\, and \,$U_2 := J(K_- \cdot H_- \oplus \liem_{\lambda_2}') \subset L_+(A)$\,. 
Moreover, we have \,$U_1 \perp U_2$\,, and \,$U_1$\, resp.~\,$U_2$\, is a quaternionic, totally complex, or totally real space according to the 
isomorphy type of \,$K_+$\, resp.~\,$K_-$\,. Therefore \,$\liem'$\, is of type \,$(\PtP, (\K_1,1+\tfrac{1}{r_1} n_{\lambda_1}'),  
(\K_2,1+\tfrac{1}{r_2} n_{\lambda_2}'))$\,, where \,$\K_1$\, resp.~\,$\K_2$\, is \,$\R$\,, \,$\C$\, or \,$\HH$\,, according to the isomorphy type of \,$K_+$\,
resp.~\,$K_-$\,, and we put \,$r_\nu := \dim_{\R}(\K_\nu) \in \{1,2,4\}$\,.  

\bigskip

\textbf{Case (3).}
We suppose that one, but not both, of the roots \,$\lambda_1$\, and \,$\lambda_2$\, is a member of \,$\Delta'$\,; without loss of generality we may suppose
\,$\lambda_1 \in \Delta'$\, and \,$\lambda_2 \not\in \Delta'$\,. Because \,$\Delta'$\, is invariant under its Weyl transformation group,
we then also have \,$\lambda_3,\lambda_4 \not\in \Delta'$\,. 

Let us again consider the sub-field \,$K := \Menge{c \in \HH}{c\,H_+ \in \liem'}$\, of \,$\HH$\, and choose \,$\K \in \{\R,\C,\HH\}$\, according to the
isomorphy type of \,$K$\,, then according to Lemma~\ref{L:claH:mstructure}, \,$\liem_{\lambda_1}'$\, is of \,$\HP$-type \,$(\K,\tfrac{1}{r}\,n_{\lambda_1}')$\, with
\,$r := \dim_{\R}(\K)$\,, hence \,$K\cdot H_+ \oplus \liem_{\lambda_1}'$\, is of \,$\HP$-type \,$\tau_1 := (\K,1+\tfrac{1}{q}\,n_{\lambda_1}')$\,.
Moreover \,$\R\,H_- \oplus \liem_{2\lambda_2}'$\, is of \,$\HP$-type \,$\tau_2$\,, which is \,$(\R,1)$\,, \,$(\C,2)$\,, \,$(\Sph^3)$\,
or \,$(\HH,1)$\,, according to whether \,$n_{2\lambda_2}'$\, is \,$0$\,, \,$1$\,, \,$2$\, or \,$3$\,, respectively. 
Thus we have
$$ \liem' = \underbrace{\liea \oplus \liem_{2\lambda_1}'}_{= K\cdot H_+ \oplus \R\cdot H_-} \oplus \liem_{2\lambda_2}' \oplus \liem_{\lambda_1}' = U_1 \oplus J(U_2) $$
with \,$U_1 := K\cdot H_+ \oplus \liem_{\lambda_1}' \subset L_+(A)$\, and \,$U_2 := J(\R\cdot H_- \oplus \liem_{2\lambda_2}') \subset L_+(A)$\,. 
We have \,$U_1 \perp J(U_2)$\,, and
\,$U_1$\, resp.~\,$U_2$\, is of \,$\HP$-type \,$\tau_1$\, resp.~\,$\tau_2$\,. Therefore \,$\liem'$\, is of type
\,$(\PtP, \tau_1,\tau_2)$\,.

\bigskip

\textbf{Case (4).} 
Let us next consider the case \,$\lambda_1,\lambda_2 \not\in \Delta'$\,, \,$\lambda_3,\lambda_4 \in \Delta'$\,. 
Similarly as in case (1), there exists \,$J \in \frakJ_A$\, and a basis
\,$(e_+,e_-)$\, adapted to \,$(A,J)$\, so that \,$M_{1,-1}^{(e_+,e_-)} \in \liem_{\lambda_3}'$\, holds. We now consider the \,$\R$-linear subspaces
\begin{gather*}
C_{\lambda_\nu} := \Menge{c \in \HH}{M_{c,\mp 1}^{(e_+,e_-)} \in \liem_{\lambda_\nu}'} \;\subset\; \HH \qmq{for \,$\nu \in \{3,4\}$\,} \\
\qmq{and} D_{2\lambda_\nu} := \Menge{d \in \IM(\HH)}{d\,H_\pm \in \liem_{2\lambda_\nu}'} \;\subset\; \IM(\HH) \qmq{for \,$\nu \in \{1,2\}$\,.} 
\end{gather*}
(Unlike the similarly defined subspaces in previous cases, they are not generally sub-fields of \,$\HH$\,.)
We have \,$1 \in C_{\lambda_3}$\,. 

From the fact that \,$\liem'$\, is a Lie triple system, we derive the following inclusions via Lemma~\ref{L:claH:Rroot}:
\begin{align}
\label{eq:claH:rk2:D1}
C_{\lambda_3} \cdot D_{2\lambda_1},\;D_{2\lambda_2}\cdot C_{\lambda_3} \subset C_{\lambda_4}
& \qmq{and}
C_{\lambda_4} \cdot D_{2\lambda_1},\;D_{2\lambda_2}\cdot C_{\lambda_4} \subset C_{\lambda_3} \; ; \\
\label{eq:claH:rk2:D2}
\IM(C_{\lambda_3}\cdot C_{\lambda_4}) \subset D_{2\lambda_1} & \qmq{and} \IM(C_{\lambda_4}\cdot C_{\lambda_3}) \subset D_{2\lambda_2} \; . 
\end{align}
Indeed, for the first inclusion of \eqref{eq:claH:rk2:D1}, let \,$c \in C_{\lambda_3}$\, and \,$d \in D_{2\lambda_1}$\, be given. By definition, we then
have \,$M_{c,-1}^{(e_+,e_-)},\;d\,H_+ \in \liem'$\,, and therefore also \,$\liem' \ni R(H_+,d\,H_+)M_{c,-1}^{(e_+,e_-)} \overset{(*)}{=} M_{-2cd,1}^{(e_+,e_-)}$\,
(where the equals sign marked $(*)$ follows from Lemma~\ref{L:claH:Rroot}(d)),
hence \,$M_{cd,1}^{(e_+,e_-)} \in \liem_{\lambda_4} \cap \liem' = \liem_{\lambda_4}'$\, and therefore \,$cd \in C_{\lambda_4}$\,. The other inclusions of 
\eqref{eq:claH:rk2:D1} are shown similarly.

For the first inclusion of \eqref{eq:claH:rk2:D2}, let \,$c \in C_{\lambda_4}$\, and \,$\wt{c} \in C_{\lambda_3}$\, be given. 
Because \,$C_{\lambda_3}$\, is an \,$\R$-linear subspace of \,$\HH$\, with \,$1\in C_{\lambda_3}$\,,
it is invariant under quaternionic conjugation, and hence we also have \,$\overline{\wt{c}} \in C_{\lambda_3}$\,. It follows from Lemma~\ref{L:claH:Rroot}(g)
that \,$D_{2\lambda_1} \ni \IM(\overline{c}\cdot\overline{\wt{c}}) =-\IM(\wt{c}\cdot c)$\, and hence \,$\IM(\wt{c}\cdot c) \in D_{2\lambda_1}$\, holds.
This shows the first inclusion of \eqref{eq:claH:rk2:D2}; the second inclusion is shown in the same way.

Because of \,$1 \in C_{\lambda_3}$\, we derive from \eqref{eq:claH:rk2:D1}: \,$D_{2\lambda_\nu} \subset C_{\lambda_4}$\,,
and from \eqref{eq:claH:rk2:D2}: \,$\IM(C_{\lambda_4}) \subset D_{2\lambda_\nu}$\,. Because of these inclusions and \,$D_{2\lambda_\nu} \subset \IM(\HH)$\,,
we have \,$D_{2\lambda_\nu} \subset C_{\lambda_4} \cap \IM(\HH) \subset \IM(C_{\lambda_4}) \subset D_{2\lambda_\nu}$\,. In this way we conclude
\begin{equation}
\label{eq:claH:rk2:Did}
D_{2\lambda_1} = D_{2\lambda_2} = C_{\lambda_4} \cap \IM(\HH) = \IM(C_{\lambda_4}) \; . 
\end{equation}
It follows that we have
\begin{equation}
\label{eq:claH:rk2:D4alt}
\text{either \quad \,$C_{\lambda_4} = D_{2\lambda_\nu}$\, \quad or \quad \,$C_{\lambda_4} = \R \oplus D_{2\lambda_\nu}$\,.}
\end{equation}

We have \,$n_{2\lambda_\nu}' = \dim(D_{2\lambda_\nu})$\,, and therefore by \eqref{eq:claH:rk2:Did}, \,$n_{2\lambda_1}' = n_{2\lambda_2}' =: n_2' \in \{0,\dotsc,3\}$\,. 
To finish off the present case of the classification, 
we now consider the four possible values for \,$n_2'$\, separately.

First suppose \,$n_2' = 0$\,, and hence \,$D_{2\lambda_\nu} = \{0\}$\,. By \eqref{eq:claH:rk2:D4alt} and the fact that \,$\lambda_4\in\Delta'$\, holds, we therefore have
\,$C_{\lambda_4} = \R$\,. Because we thus have in particular \,$1 \in C_{\lambda_4}$\,, \eqref{eq:claH:rk2:D2} shows \,$\IM(C_{\lambda_3}) \subset D_{2\lambda_1} = \{0\}$\,,
hence \,$C_{\lambda_3} = \R$\,. Thus we have
$$ \liem' = \liea \oplus \liem_{\lambda_3}' \oplus \liem_{\lambda_4}' 
= \R H_+ \oplus \R H_- \oplus \underbrace{\R\,M_{1,-1}^{(e_+,e_-)} \oplus \R\,M_{1,1}^{(e_+,e_-)}}_{= \R\,J(H_+) \oplus \R\,J(H_-)}
= U \oplus J(U) $$
with the subspace \,$U := \R\,H_+ \oplus \R\,J(H_-)$\, of \,$L_+(A)$\, of \,$\HP$-type \,$(\R,2)$\,. Thus \,$\liem'$\, is of type \,$(\Gtwo,(\R,2))$\,.

Now suppose \,$n_2' = 1$\,. Then there exists \,$i \in \Sph(\IM(\HH))$\, with \,$D_{2\lambda_1} \overset{\eqref{eq:claH:rk2:Did}}{=} D_{2\lambda_2} = \R\,i$\,,
and by \eqref{eq:claH:rk2:D4alt} we have either \,$C_{\lambda_4} = \R\,i$\, or \,$C_{\lambda_4} = \R \oplus \R\,i$\,. In either case we have
\,$\IM(C_{\lambda_3} \cdot C_{\lambda_4}) \subset D_{2\lambda_1} = \R\,i$\, by \eqref{eq:claH:rk2:D2}, and therefore \,$C_{\lambda_3} \cdot C_{\lambda_4} \subset
\R \oplus \R\,i$\,. Because of \,$i \in C_{\lambda_4}$\,, the latter inclusion implies \,$C_{\lambda_3} \subset \R \oplus \R\,i$\,; because of \,$1 \in C_{\lambda_3}$\,,
we in fact have either \,$C_{\lambda_3} = \R$\, or \,$C_{\lambda_3} = \R \oplus \R\,i$\,. 
Because the dimensions of \,$C_{\lambda_3}$\, and \,$C_{\lambda_4}$\, are equal (the roots \,$\lambda_3$\, and \,$\lambda_4$\, correspond to each other 
under the Weyl transformation induced by \,$2\lambda_1 \in \Delta'$\,, and have therefore the same multiplicity), 
we have either \,$C_{\lambda_3} = C_{\lambda_4} = \R \oplus \R\,i$\,, or \,$C_{\lambda_3} = \R$\, and \,$C_{\lambda_4} = \R\,i$\,.

If the former case holds, we have
$$ \liem' = \spn_{\R}\{H_+, H_-, i\,H_+, i\,H_-, M_{1,-1}^{(e_+,e_-)}, M_{i,-1}^{(e_+,e_-)}, M_{1,1}^{(e_+,e_-)}, M_{i,1}^{(e_+,e_-)}\} = U \oplus J(U) $$
with the subspace \,$U:=\C\,H_+ \oplus \C\,J(H_-)$\, of \,$L_+(A)$\, of \,$\HP$-type \,$(\C,2)$\,. Hence \,$\liem'$\, then is of type \,$(\Gtwo, (\C,2))$\,.

Otherwise we have \,$C_{\lambda_3} = \R$\, and \,$C_{\lambda_4} = \R\,i$\,, and therefore
$$ \liem' = \spn_{\R}\{H_+, H_-, i\,H_+, i\,H_-, M_{1,-1}^{(e_+,e_-)}, M_{i,1}^{(e_+,e_-)}\} \;, $$
hence \,$\liem'$\, is the \,$\R$-ortho-complement of \,$\spn_{\R}\{M_{i,-1}^{(e_+,e_-)},M_{1,1}^{(e_+,e_-)}\} = \spn_{\R}\{v,v\,i\}$\, with \,$v := M_{i,-1}^{(e_+,e_-)}$\,
in \,$U \oplus J(U)$\,. We have \,$\vi(v)=\tfrac\pi4$\,, and therefore \,$\liem'$\, is of type \,$(\Qthree)$\,. 

Next suppose \,$n_2' = 2$\,. Then there exists a canonical basis \,$(i,j,k)$\, of \,$\IM(\HH)$\, 
so that \,$D_{2\lambda_1} = D_{2\lambda_2} = \R\,i \oplus \R\,j$\, holds. By \eqref{eq:claH:rk2:D4alt} we have either \,$C_{\lambda_4} = \R\,i \oplus \R\,j$\,
or \,$C_{\lambda_4} = \R \oplus \R\,i \oplus \R\,j$\,. 
 In either case, we have \,$C_{\lambda_3} \supset C_{\lambda_4} \cdot D_{2\lambda_1} \supset \R \oplus \R\,k$\, by \eqref{eq:claH:rk2:D1}.
Assume that we had \,$C_{\lambda_4} = \R \oplus \R\,i \oplus \R\,j$\,. Then we would have \,$k = 1\cdot k \in \IM(C_{\lambda_3} \cdot C_{\lambda_4})$\,, and therefore
\,$k \in D_{2\lambda_1}$\, by \eqref{eq:claH:rk2:D2}, a contradiction. Therefore we have \,$C_{\lambda_4} = \R\,i \oplus \R\,j$\,, and thus
because of \,$\dim C_{\lambda_3} = \dim C_{\lambda_4}$\,, \,$C_{\lambda_3} = \R \oplus \R\,k$\,. Thus we have
$$ \liem' = (\R \oplus \R\,i \oplus \R\,j)H_+ \,\oplus\,(\R \oplus \R\,i \oplus \R\,j)H_- \,\oplus\, 
\spn_{\R}\{M_{1,-1}^{(e_+,e_-)}, M_{i,1}^{(e_+,e_-)}, M_{j,1}^{(e_+,e_-)}, M_{k,-1}^{(e_+,e_-)} \} \; . $$
This shows that we have \,$\liem' = \Menge{\Phi \circ X}{X \in \liesp(V')}$\,, where the \,$\HH$-linear isometry \,$\Phi: V' \to V$\, is given by
\,$\Phi(e_+) = k\cdot H_+(e_+)$\, and \,$\Phi(e_-) = k\cdot H_-(e_-)$\,. 
Therefore \,$\liem'$\, is of type \,$(\Sptwo)$\,. 

Finally suppose \,$n_2' = 3$\,. Then we have \,$D_{2\lambda_\nu} = \IM(\HH)$\,, hence \,$\IM(\HH) \subset C_{\lambda_4}$\, by \eqref{eq:claH:rk2:D4alt}, 
and therefore by \eqref{eq:claH:rk2:D1}:
\,$C_{\lambda_3} \supset C_{\lambda_4} \cdot D_{2\lambda_1} \supset \IM(\HH) \cdot \IM(\HH) = \HH$\,, i.e.~\,$C_{\lambda_3} = \HH$\,, 
and thus also \,$C_{\lambda_4} \supset C_{\lambda_3}
\cdot D_{2\lambda_1} = \HH \cdot \IM(\HH) = \HH$\,, i.e.~\,$C_{\lambda_4} = \HH$\,. Thus we have \,$\liem_{\lambda}' = \liem_{\lambda}$\, for \,$\lambda \in \{\lambda_3,
\lambda_4,2\lambda_1,2\lambda_2\}$\,, and therefore
$$ \liem' \;=\; \underbrace{\liea \oplus \liem_{2\lambda_1}' \oplus \liem_{2\lambda_2}'}_{= \HH\,H_+ \oplus \HH\,H_-}
\oplus \underbrace{\liem_{\lambda_3}' \oplus \liem_{\lambda_4}'}_{= \HH\,J(H_+) \oplus \HH\,J(H_-)} \;=\; U \oplus J(U) $$
with the subspace \,$U := \HH\,H_+ \oplus \HH \, J(H_-)$\, of \,$L_+(A)$\, of \,$\HP$-type \,$(\HH,2)$\,. Thus \,$\liem'$\, is of type \,$(\Gtwo,(\HH,2))$\,. 

\bigskip

\textbf{Case (5).}
Here we suppose \,$\lambda_1,\lambda_2 \not\in \Delta'$\, and that either, but not both, of \,$\lambda_3$\, and \,$\lambda_4$\, is a member of \,$\Delta'$\,;
we may suppose without loss of generality that \,$\lambda_4 \in \Delta'$\,, \,$\lambda_3 \not\in \Delta'$\, holds. Because of the invariance of \,$\Delta'$\,
under its Weyl transformation group, we have \,$2\lambda_1, 2\lambda_2 \not\in \Delta'$\,, and therefore \,$\liem' = \liea \oplus \liem_{\lambda_4}'$\,. It
follows that \,$\liem'$\, is of type \,$(\StSsqrt,1+n_{\lambda_4}')$\,. 

\bigskip

\textbf{Case (6).}
Finally, we suppose \,$\lambda_1,\lambda_2, \lambda_3,\lambda_4 \not\in \Delta'$\,, and therefore 
we have \,$\liem' = \liea \oplus \liem_{2\lambda_1}' \oplus \liem_{2\lambda_2}' = U_1 \operp J(U_2)$\, 
with \,$U_1 := \R H_+ \oplus \liem_{2\lambda_1}' \subset L_+(A)$\, and \,$U_2 := J(\R H_- \oplus \liem_{2\lambda_2}') \subset L_+(A)$\,. We have 
\,$U_1 \perp U_2$\, and \,$U_\nu$\, is of \,$\HP$-type \,$\tau_\nu$\,, which is \,$(\R,1)$\,, \,$(\C,1)$\,, \,$(\Sph^3)$\, or \,$(\HH,1)$\, according to whether
\,$n_{2\lambda_\nu}'$\, is \,$0$\,, \,$1$\,, \,$2$\, or \,$3$\,, respectively. Consequently \,$\liem'$\, is of type \,$(\PtP, \tau_1, \tau_2)$\,. 

\bigskip

The concludes the classification for the case of rank 2.

\subsection{The case of rank 1}
\label{SSe:claH:rk1}

We now let a Lie triple system \,$\liem'$\, of \,$\liem$\, with \,$\rk(\liem')=1$\, be given. If \,$\dim(\liem')=1$\, holds, then we have \,$\liem' = \R\,H$\,
with \,$H \in \Sph(\liem')$\,, and therefore \,$\liem'$\, then is of type \,$(\Geo,\vi=\vi(H))$\,. Thus we may suppose \,$\dim(\liem') \geq 2$\, in the sequel.

Because \,$\liem'$\, is of rank 1, any two unit vectors of \,$\liem'$\, are congruent
under the isotropy action of \,$G_2(V'\oplus V)$\,. Therefore all vectors in \,$\liem' \setminus \{0\}$\, have the same characteristic angle
(see Proposition~\ref{P:orbits:orbits}) \,$\vi_0 \in [0,\tfrac\pi4]$\,. 

We fix a unit vector \,$H \in \Sph(\liem')$\,, then \,$\liea' := \R\,H$\, is a Cartan subalgebra of \,$\liem'$\,. 
We choose a Cartan subalgebra \,$\liea$\, of \,$\liem$\, so that \,$\liea' = \liea \cap \liem'$\, holds. (Such Cartan subalgebras exist, as was discussed
in Section~\ref{Se:generallts}. In fact, for \,$\vi_0 \not\in \{0,\tfrac\pi4\}$\,, \,$\liea$\, is unique.)
Because of
\,$\dim(\liem') > \rk(\liem')$\,, we have \,$\Delta' := \Delta(\liem',\liea') \neq \varnothing$\,. Consider \,$\alpha \in \Delta'$\,. \,$\alpha$\,
is either elementary or composite in the sense of Definition~\ref{D:cla:subroots:Elemcomp}. Proposition~\ref{P:cla:subroots:Comp} shows that
if \,$\alpha$\, is elementary, there exists \,$\lambda \in \Delta := \Delta(\liem,\liea)$\, with \,$\lambda^\sharp \in \liea'$\, and therefore 
\,$\vi(H) = \vi(\lambda^\sharp)$\,; if \,$\alpha$\, is composite, there exist \,$\lambda,\mu \in \Delta$\, with \,$\lambda \neq \mu$\,
and \,$\lambda^\sharp - \mu^\sharp$\, being orthogonal to \,$\liea'$\,, and therefore \,$\vi(H) = \vi(\lambda^\sharp - \mu^\sharp)$\,
(note that if \,$v,v' \in \liea \setminus \{0\}$\, are orthogonal to each other, we have \,$\vi(v) = \vi(v')$\,).
We have
$$ \vi(\lambda_1^\sharp) = \vi(\lambda_2^\sharp) = \vi(2\lambda_1^\sharp) = \vi(2\lambda_2^\sharp) = 0, \quad \vi(\lambda_3^\sharp) = \vi(\lambda_4^\sharp) = \tfrac\pi4 $$
and
$$ \forall \lambda,\mu \in \Delta,\;\lambda \neq \mu \; : \; \vi(\lambda^\sharp - \mu^\sharp) \in \{0,\arctan(\tfrac13),\arctan(\tfrac12),\tfrac\pi4\}\; . $$
Thus we see that in any case
$$ \vi_0 = \vi(H) \in \{0,\arctan(\tfrac13),\arctan(\tfrac12),\tfrac\pi4\}\; $$
holds. Below we will handle the cases corresponding to these four possible values of \,$\vi_0$\, separately.

In preparation to their treatment, we note that by Proposition~\ref{P:cla:subroots:subroots-neu}(a) we have
\begin{equation}
\label{eq:claH:rk1:Delta'}
\Delta' \;\subset\; \;\bigr\{\;\lambda(H)\, \alpha_0\; \bigr| \;\lambda \in \Delta,\; \lambda(H) \neq 0\; \bigr\} 
\end{equation}
with the linear form \,$\alpha_0: \liea' \to \R,\;tH \mapsto t$\,, and
\begin{equation}
\label{eq:claH:rk1:m'decomp}
\liem' = \liea' \oplus \bigoplus_{\alpha \in \Delta'_+} \liem_\alpha'
\end{equation}
with
\begin{equation}
\label{eq:claH:rk1:malpha'}
\forall \alpha \in \Delta'\;:\; \liem_\alpha' = \left( \bigoplus_{\substack{\lambda \in \Delta \\ \lambda|\liea' = \alpha}} \liem_\lambda \right) \;\cap\; \liem' \; .
\end{equation}

We also note that by Proposition~\ref{P:orbits:orbits}(b) there exists \,$A \in \A$\,
and an orthonormal basis \,$(H_+,H_-)$\, of \,$\liea$\, with \,$H_\pm \in L_\pm(A)$\, and \,$H_+(V'_+(A)) \perp H_-(V'_-(A))$\, so that
\begin{equation}
\label{eq:claH:rk1:H}
H = \cos(\vi_0)\,H_+ + \sin(\vi_0)\,H_- 
\end{equation}
holds. 

\bigskip

\textbf{The case \,$\boldsymbol{\vi_0=0}$\,.}
In this case we have \,$H = H_+$\, by Equation~\eqref{eq:claH:rk1:H}, and therefore 
$$ \lambda_1(H) = \lambda_3(H) = \lambda_4(H) = 1,\quad 2\lambda_1(H) = 2,\quad \lambda_2(H) = 2\lambda_2(H) = 0 \; . $$
If we consider the linear form \,$ \alpha: \liea' \to \R,\; tH \mapsto t$\,,
we therefore have \,$\Delta' \subset \{\pm\alpha, \pm 2\alpha\}$\, by Equation~\eqref{eq:claH:rk1:Delta'}, and 
\begin{equation}
\label{eq:claH:rk1:0:m'decomp}
\liem' = \liea' \oplus \liem_{\alpha}' \oplus \liem_{2\alpha}'
\end{equation}
with 
$$ \liem_\alpha' := \bigr(\liem_{\lambda_1} \oplus \liem_{\lambda_3} \oplus \liem_{\lambda_4}\bigr) \,\cap\, \liem' 
\;\subset\; \Menge{v \in L_+(A)}{v(V_+'(A)) \perp H_+(V_+'(A))} \oplus \HH\,J(H_+)$$
(where we fix \,$J \in \frakJ_A$\,) and 
$$ \liem_{2\alpha}' := \liem_{2\lambda_1} \cap \liem' \;\subset\; \IM(\HH)\,H_+ $$
by Equations~\eqref{eq:claH:rk1:m'decomp} and \eqref{eq:claH:rk1:malpha'}.

Let \,$v \in \liem_{\alpha}'$\, be given, say \,$v = v_+ + c\,J(H_+)$\, with \,$v_+ \in L_+(A)$\,, \,$v_+(V_+'(A)) \perp H_+(V_+'(A))$\, and \,$c \in \HH$\,. 
We fix a basis \,$(e_+,e_-)$\, adapted to \,$(A,J)$\,, then we calculate by means of Equations~\eqref{eq:claH:R11} and \eqref{eq:claH:R12}:
$$ R(H_+,v)v = R(H_+,v_+)v + R(H_+,c\,J(H_+))v = 
H_+ \, (\|v(e_+)\|^2 + |c|^2) -2\,J(v)\cdot c \; . $$
Because \,$\liem'$\, is a Lie triple system, it follows that \,$J(v)\cdot c \in \liem'$\, holds. However Equation~\eqref{eq:claH:rk1:0:m'decomp} shows
that \,$J(v)\cdot c$\, is also orthogonal to \,$\liem'$\,. Hence we have \,$J(v)\cdot c = 0$\, and therefore
either \,$v=0$\, or \,$c=0$\,. Thus we have shown \,$\liem_{\alpha}' \subset \Menge{v \in L_+(A)}{\g{v(e_+)}{H_+(e_+)} =0} \cup \HH\,J(H_+)$\,; because
\,$\liem_{\alpha}'$\, is a linear space, it follows that
\begin{equation*}
\qmq{either} \liem_{\alpha}' \subset \Menge{v \in L_+(A)}{\g{v(e_+)}{H_+(e_+)} =0} \qmq{or} \liem_{\alpha}' \subset \HH\,J(H_+)
\end{equation*}
holds. We now handle these two possibilities separately.

First suppose that \,$\liem_{\alpha}' \subset \Menge{v \in L_+(A)}{\g{v(e_+)}{H_+(e_+)} =0}$\, and therefore \,$\liem' \subset L_+(A)$\, holds. \,$L_+(A)$\, itself
is a Lie triple system (of type \,$(\PP,\vi=0,(\HH,n))$\,) in \,$\liem$\,, 
and Equation~\eqref{eq:claH:R11} shows that the restriction of the curvature tensor \,$R$\, of \,$G_2(V'\oplus V)$\,
to \,$L_+(A)$\, is the curvature tensor of the quaternion projective space \,$\HP^n$\,. \,$\liem'$\, also is a Lie triple system if regarded as a subspace of 
\,$L_+(A)$\,, its position is therefore determined by the known classification of the totally geodesic submanifolds of \,$\HP^n$\,. Therefore there exists
a \,$\HP$-type \,$\tau$\, so that \,$\liem'$\, is of \,$\HP$-type \,$\tau$\, in \,$L_+(A)$\,, and therefore of type \,$(\PP,\vi=0,\tau)$\, in \,$\liem$\,. 

Now suppose that \,$\liem_{\alpha}' \subset \HH\,J(H_+)$\, and therefore \,$\liem' \subset \HH\,H_+ \oplus \HH J(H_+) =: \wh{\liem}'$\, holds. Again \,$\wh{\liem}'$\,
itself is a Lie triple system in \,$\liem$\,, namely of type \,$(\Gtwo,(\HH,1))$\,, and the restriction of \,$R$\, to \,$\wh{\liem}'$\, is the curvature
tensor of \,$G_2(\HH^3) \cong \HP^2$\,. By the same argument as in the preceding case, there exists a \,$\HP$-type \,$\tau$\, so that \,$\liem'$\, is of type
\,$\tau$\, in \,$\wh{\liem}$\,. If \,$\tau = (\K,2)$\, holds with \,$\K \in \{\R,\C,\HH\}$\,, then \,$\liem'$\, is of type \,$(\Gtwo,(\K,1))$\,;
if \,$\tau = (\K,1)$\, or \,$\tau = (\Sph^3)$\, holds, then \,$\liem'$\, is of type \,$(\PP,\vi=0,\tau)$\,.

\bigskip

\textbf{The case \,$\boldsymbol{\vi_0=\arctan(\tfrac13)}$\,.}
In this case, we have by Equation~\eqref{eq:claH:rk1:H}
\begin{equation}
\label{eq:claH:rk1:13:H}
H = \tfrac{3}{\sqrt{10}} \,H_+ + \tfrac{1}{\sqrt{10}}\,H_- 
\end{equation}
and therefore
$$ \lambda_1(H) = \tfrac{3}{\sqrt{10}} \;,\; \lambda_2(H) = \tfrac{1}{\sqrt{10}} \;,\; \lambda_3(H) = \tfrac{4}{\sqrt{10}} \;,\; 
\lambda_4(H) = \tfrac{2}{\sqrt{10}} \;,\; 2\lambda_1(H) = \tfrac{6}{\sqrt{10}} \;,\; 2\lambda_2(H) = \tfrac{2}{\sqrt{10}} \; . $$
Because there does not exist a root \,$\lambda \in \Delta$\, with \,$\lambda^\sharp \in \liea'$\,, Proposition~\ref{P:cla:subroots:Comp}(a) shows that
every root of \,$\Delta'$\, is composite in the sense of Definition~\ref{D:cla:subroots:Elemcomp}. Therefore we have \,$\Delta' = \{\pm \alpha\}$\,
with the linear form \,$\alpha: \liea' \to \R,\; tH \mapsto \tfrac{2}{\sqrt{10}}\,t$\,,
and by Equations~\eqref{eq:claH:rk1:m'decomp} and \eqref{eq:claH:rk1:malpha'} we have
$$ \liem' = \liea' \oplus \liem_\alpha' $$
with 
\begin{equation}
\label{eq:claH:rk1:13:malpha'}
\liem_\alpha' := (\liem_{\lambda_4} \oplus \liem_{2\lambda_2})\,\cap\, \liem' \; .
\end{equation}
We have \,$\alpha^\sharp = \tfrac{2}{\sqrt{10}}\cdot H = \tfrac35\,H_+ + \tfrac15\,H_- = \tfrac35\,\lambda_4^\sharp + \tfrac25\,(2\lambda_2)^\sharp$\,,
and therefore Proposition~\ref{P:cla:skew} shows that there exists a \,$\R$-linear subspace \,$C \subset \HH$\, and an \,$\R$-linear isometry
\,$\Phi: C \to \IM(\HH)$\, so that
\begin{equation}
\label{eq:claH:rk1:13:malpha'2}
\liem_\alpha' = \Menge{M_{c,1}^{(e_+,e_-)} + \sqrt{\tfrac23}\,\Phi(c)\,H_-}{c \in C}
\end{equation}
holds, herein \,$(e_+,e_-)$\, is any basis adapted to \,$A$\,; because of Equation~\eqref{eq:curv:Mtrafo} the basis \,$(e_+,e_-)$\, can be chosen in such
a way that we have \,$1 \in C$\,.

In the case \,$n_\alpha' = 1$\, this already shows that \,$\liem' = \spn_{\R}\{H,\;M_{1,1}^{(e_+,e_-)} + \sqrt{\tfrac23}\,\Phi(1)\,H_-\}$\, is of type
\,$(\PP,\vi=\arctan(\tfrac13),2)$\,.

So we now suppose \,$\dim C = n_\alpha' \geq 2$\, and let \,$c_1,c_2 \in C$\, be given. By Equation~\eqref{eq:claH:rk1:13:malpha'2} we have
 \,$v_\nu := M_{c_\nu,1}^{(e_+,e_-)} + \sqrt{\tfrac23}\,\Phi(c)\,H_- \in \liem_\alpha'$\, for \,$\nu \in \{1,2\}$\,.
Using Equations~\eqref{eq:claH:R11}, \eqref{eq:claH:R12} and \eqref{eq:claH:R22}
we calculate \,$R(H,v_1)v_2$\,, which is again an element of \,$\liem'$\, because \,$\liem'$\, is a Lie triple system:
\begin{equation}
\label{eq:claH:rk1:13:R}
R(H,v_1)v_2 = \tfrac23\,\RE(\overline{c_1}\,c_2)\,H + \tfrac{2}{\sqrt{15}}\,M_{\Phi(c_2)\,c_1 - \Phi(c_1)\,c_2\,,\,-1}^{(e_+,e_-)} \; . 
\end{equation}
Because \,$\liem_{\lambda_3}$\, is orthogonal to \,$\liem'$\,, the \,$\liem_{\lambda_3}$-component of \eqref{eq:claH:rk1:13:R}, which is proportional to 
\,$M_{\Phi(c_2)\,c_1 - \Phi(c_1)\,c_2\,,\,-1}^{(e_+,e_-)}$\,, vanishes, and thus we have shown
\begin{equation}
\label{eq:claH:rk1:13:c1c2}
\forall c_1,c_2 \in C \; : \; \Phi(c_2)\cdot c_1 = \Phi(c_1) \cdot c_2 \; .
\end{equation}
By specializing \,$c_1=1$\, in this equation, we obtain in particular
\begin{equation}
\label{eq:claH:rk1:13:c}
\forall c \in C \; : \; \Phi(c) = \Phi(1) \cdot c \; .
\end{equation}

It follows that the case \,$n_\alpha' \geq 3$\, does not occur. In fact, under the assumption \,$n_\alpha' \geq 3$\, there would exist a canonical basis
\,$(i,j,k)$\, of \,$\IM(\HH)$\, so that \,$1,i,j \in C$\, holds. We would then have on one hand by Equation~\eqref{eq:claH:rk1:13:c}:
\,$\Phi(j) = \Phi(1)\cdot j$\,, on the other hand by Equation~\eqref{eq:claH:rk1:13:c1c2}: \,$\Phi(j) = \Phi(i)\cdot j\cdot i^{-1} = \Phi(i)\cdot k
\overset{\eqref{eq:claH:rk1:13:c}}{=} \Phi(1)\cdot i\cdot k = -\Phi(1)\cdot j$\,, and therefore \,$\Phi(j)=0$\,, in contradiction to the fact that 
\,$\Phi: C \to \IM(\HH)$\, is a linear isometry.

Therefore the only remaining case is that of \,$n_\alpha' = 2$\,. In this case, there exists \,$j \in \Sph(\IM(\HH))$\, so that \,$C = \R\oplus \R\,j$\, holds. 
Because \,$\IM(\HH) \ni \Phi(j) \overset{\eqref{eq:claH:rk1:13:c}}{=} \Phi(1)\cdot j$\, holds, \,$i := \Phi(1) \in \Sph(\IM(\HH))$\, must be orthogonal to \,$j$\,.
Therefore \,$(i,j,k)$\, with \,$k:=ij$\, is a canonical basis of \,$\IM(\HH)$\,, and we have \,$\Phi(j) \overset{\eqref{eq:claH:rk1:13:c}}{=} \Phi(1)\cdot j = k$\,.
Therefore \,$\liem' = \spn_{\R}\{H,\; M_{1,1}^{(e_+,e_-)} + \sqrt{\tfrac23}\,i\,H_-,\;M_{j,1}^{(e_+,e_-)} + \sqrt{\tfrac23}\,k\,H_-\}$\,
is of type \,$(\Sph,\vi=\arctan(\tfrac13),3)$\,. 

\bigskip

\textbf{The case \,$\boldsymbol{\vi_0=\arctan(\tfrac12)}$\,.}
In this case, we have by Equation~\eqref{eq:claH:rk1:H}
\begin{equation}
\label{eq:claH:rk1:12:H}
H = \tfrac{2}{\sqrt{5}} \,H_+ + \tfrac{1}{\sqrt{5}}\,H_- 
\end{equation}
and therefore
$$ \lambda_1(H) = \tfrac{2}{\sqrt{5}} \;,\; \lambda_2(H) = \tfrac{1}{\sqrt{5}} \;,\; \lambda_3(H) = \tfrac{3}{\sqrt{5}} \;,\; 
\lambda_4(H) = \tfrac{1}{\sqrt{5}} \;,\; 2\lambda_1(H) = \tfrac{4}{\sqrt{5}} \;,\; 2\lambda_2(H) = \tfrac{2}{\sqrt{5}} \; . $$
Because there does not exist a root \,$\lambda \in \Delta$\, with \,$\lambda^\sharp \in \liea'$\,, Proposition~\ref{P:cla:subroots:Comp}(a) shows that
every root of \,$\Delta'$\, is composite in the sense of Definition~\ref{D:cla:subroots:Elemcomp}. It follows that we have \,$\Delta' \subset \{\pm \alpha,\pm 2\alpha\}$\,
with the linear form \,$\alpha: \liea' \to \R,\; tH \mapsto \tfrac{1}{\sqrt{5}}\,t$\,
and by Equations~\eqref{eq:claH:rk1:m'decomp} and \eqref{eq:claH:rk1:malpha'} we have
\begin{equation}
\label{eq:claH:rk1:12:m'}
\liem' = \liea' \oplus \liem_\alpha' \oplus \liem_{2\alpha}' 
\end{equation}
with 
\begin{equation}
\label{eq:claH:rk1:12:malpha'def}
\liem_\alpha' := (\liem_{\lambda_2} \oplus \liem_{\lambda_4})\,\cap\, \liem' \qmq{and}
\liem_{2\alpha}' := (\liem_{\lambda_1} \oplus \liem_{2\lambda_2})\,\cap\, \liem' \; .
\end{equation}
We have \,$\alpha^\sharp = \tfrac{1}{\sqrt{5}}\,H = \tfrac25\,H_+ + \tfrac15\,H_- = \tfrac{2}{5}\,\lambda_4^\sharp + \tfrac{3}{5}\,\lambda_2^\sharp$\,. By 
Proposition~\ref{P:cla:skew} it follows that for any basis \,$(e_+,e_-)$\, adapted to \,$A$\,
there exists an \,$\R$-linear subspace \,$C \subset \HH$\, and an \,$\R$-linear isometry \,$\Phi: C \to \liem_{\lambda_2}$\,
so that
\begin{equation}
\label{eq:claH:rk1:12:malpha'}
\liem_\alpha' = \Menge{M_{c,1}^{(e_+,e_-)} + \sqrt{\tfrac32}\,\Phi(c)}{c \in C}
\end{equation}
holds, and that we have \,$n_\alpha' \leq 4$\,. If \,$C \neq \{0\}$\, holds, the basis \,$(e_+,e_-)$\, adapted to \,$A$\, can be chosen in such a way
that \,$1 \in C$\, holds.

Similarly, we have \,$(2\alpha)^\sharp = \tfrac{4}{5}\,H_+ + \tfrac25\,H_- = \tfrac{1}{5}\,(2\lambda_2)^\sharp + \tfrac{4}{5}\,\lambda_1^\sharp$\,. By
Proposition~\ref{P:cla:skew} it follows that there exists an \,$\R$-linear subspace \,$D \subset \IM(\HH)$\, and an \,$\R$-linear isometry
\,$\wt{\Phi}: D \to \liem_{\lambda_1}$\, so that 
\begin{equation}
\label{eq:claH:rk1:12:m2alpha'}
\liem_{2\alpha}' = \Menge{d\,H_- + 2\,\wt{\Phi}(d)}{d \in D}
\end{equation}
holds, and that we have \,$n_{2\alpha}' \leq 3$\,. 

Let \,$c_1,c_2 \in C$\, and \,$d_1,d_2 \in D$\, be given;
by Equation~\eqref{eq:claH:rk1:12:malpha'} resp.~\eqref{eq:claH:rk1:12:m2alpha'} we then have
$$ v_\nu := M_{c_\nu,1}^{(e_+,e_-)} + \sqrt{\tfrac32}\,\Phi(c_\nu) \in \liem_{\alpha}' \qmq{and} \wt{v}_\nu := d_\nu\,H_- + 2\,\wt{\Phi}(d_\nu) \in \liem_{2\alpha}' $$
for \,$\nu \in \{1,2\}$\,. 
Via Equations~\eqref{eq:claH:R11}, \eqref{eq:claH:R12} and \eqref{eq:claH:R22} we then calculate various instances of the curvature tensor;
all the resulting vectors are again elements of \,$\liem'$\, by the fact that \,$\liem'$\, is a Lie triple system: 
\begin{align}
\label{eq:claH:rk1:12:R-vvH}
\sqrt{5}\,R(v_1,v_2)H & = -3\,H_-\cdot \IM(\g{\Phi(c_1)}{\Phi(c_2)}) + \sqrt{3}\,J\bigr(\Phi(c_1)\,c_2 - \Phi(c_2)\,c_1\bigr) \\
\label{eq:claH:rk1:12:R-Hvvt}
\sqrt{5}\,R(H,v_1)\wt{v}_1 & = -\sqrt{\tfrac32}\,\Phi(c_1)\,d_1 + \sqrt{2}\,J(\wt{\Phi}(d_1))\,\overline{c_1} \notag \\
& \hspace{1cm} 
+ \sqrt{3}\,M_{\g{J\Phi(c_1)}{\wt{\Phi}(d_1)}\,,\,1}^{(e_+,e_-)} + \sqrt{3}\,M_{\g{J\Phi(c_1)}{\wt{\Phi}(d_1)}\,,\,-1}^{(e_+,e_-)} + M_{d_1\,c_1\,,\,-1}^{(e_+,e_-)} \\
\label{eq:claH:rk1:12:R-Hvtvt}
\sqrt{5}\,R(H,\wt{v}_1)\wt{v}_2 & = 8\,H_+ \cdot \g{\wt{\Phi}(d_1)}{\wt{\Phi}(d_2)} - 4\,H_- \cdot \RE(d_1\,d_2) \\
\label{eq:claH:rk1:12:R-vvvt}
R(v_1,v_2)\wt{v}_1 & = H_+ \cdot \sqrt{3}\,\bigr(\overline{c_1}\,\g{J\Phi(c_2)}{\wt{\Phi}(d_1)} - \overline{c_2}\,\g{J\Phi(c_1)}{\wt{\Phi}(d_1)} \bigr) \notag \\
& \hspace{1cm} 
+ H_- \cdot \bigr( \IM(c_1\,\overline{c_2})\,d_1 - d_1\,\IM(c_1\,\overline{c_2}) + 3\,d_1\,\IM(\g{\Phi(c_2)}{\Phi(c_1)}) \bigr) \notag \\
& \hspace{1cm}
+ 2 \, \wt{\Phi}(d)\,\IM(\overline{c_2}\,c_1) + 3\,J(\Phi(c_1))\,\g{J\Phi(c_2)}{\wt{\Phi}(d_1)} - 3\,J(\Phi(c_2))\,\g{J\Phi(c_1)}{\wt{\Phi}(d_1)} \; . 
\end{align}

We will now use these equations to derive results concerning the structure of the data \,$(C,\Phi)$\, describing \,$\liem_\alpha'$\,, the structure of the
data \,$(D,\wt{\Phi})$\, describing \,$\liem_{2\alpha}'$\,, and the relations between these two sets of data.

Because \,$\liem_{2\lambda_1}$\, is orthogonal to \,$\liem'$\, by Equations~\eqref{eq:claH:rk1:12:m'} and \eqref{eq:claH:rk1:12:malpha'def}, 
the \,$\liem_{2\lambda_1}$-component of the element \eqref{eq:claH:rk1:12:R-Hvtvt} of \,$\liem'$\,, which is proportional to 
\,$H_+\cdot \IM(\g{\wt{\Phi}(d_1)}{\wt{\Phi}(d_2)})$\,, vanishes. This shows that \,$\wt{\Phi}(D)$\, is a totally real subspace of \,$\liem_{\lambda_1}$\,. 
It follows that \,$\liea' \oplus \liem_{2\alpha}'$\, is a Lie triple system of type \,$(\PP,\vi=\arctan(\tfrac12),\tau)$\,, where
the \,$\HP$-type \,$\tau$\, is \,$(\R,1)$\,, \,$(\C,1)$\,, \,$(\Sph^3)$\, or \,$(\HH,1)$\,, according to whether \,$n_{2\alpha}'$\, 
is \,$0$\,, \,$1$\,, \,$2$\, or \,$3$\,, respectively. Therefore if \,$n_\alpha'=0$\, holds, then already \,$\liem'$\, is of that type.

In the sequel we thus suppose \,$n_\alpha' \geq 1$\,. Hence we have \,$C \neq \{0\}$\, and therefore, by our choice of the basis \,$(e_+,e_-)$\, adapted
to \,$A$\,, \,$1 \in C$\,. 

Because \,$\liem_{\lambda_3}$\, is orthogonal to \,$\liem'$\,, the \,$\liem_{\lambda_3}$-com\-po\-nent of \eqref{eq:claH:rk1:12:R-Hvvt}, which equals
\,$\sqrt{3}\,M_{\g{J\Phi(c_1)}{\wt{\Phi}(d_1)}\,,\,-1}^{(e_+,e_-)} + M_{d_1\,c_1\,,\,-1}^{(e_+,e_-)}$\,, vanishes. Hence we have
\,$\g{J\Phi(c_1)}{\wt{\Phi}(d_1)} = -\tfrac{1}{\sqrt{3}}\,d_1\cdot c_1$\,. 
Moreover the \,$(\liem_{\lambda_4} \oplus \liem_{\lambda_2})$-component of \eqref{eq:claH:rk1:12:R-Hvvt},
which equals \,$\sqrt{3}\,M_{\g{J\Phi(c_1)}{\wt{\Phi}(d_1)}\,,\,1}^{(e_+,e_-)} - \sqrt{3/2} \, \Phi(c_1)\,d_1 + \sqrt{2}\,J(\wt{\Phi}(d_1))\,\overline{c_1}$\,,
is a member of \,$(\liem_{\lambda_4} \oplus \liem_{\lambda_2}) \cap \liem' = \liem_\alpha'$\,; it follows by Equation~\eqref{eq:claH:rk1:12:malpha'} that 
\begin{equation}
\label{eq:claH:rk1:12:dc-Phi}
-\tfrac{1}{\sqrt{3}}\,d_1\cdot c_1 = \g{J\Phi(c_1)}{\wt{\Phi}(d_1)} \;\in\; C
\end{equation}
and 
\begin{equation}
\label{eq:claH:rk1:12:dc-Phi2}
\Phi(d_1\,c_1) = \Phi(c_1)\,d_1 - \tfrac{2}{\sqrt{3}}\,J(\wt{\Phi}(d_1))\,\overline{c_1}
\end{equation}
holds. \eqref{eq:claH:rk1:12:dc-Phi} shows in particular that we have
\begin{equation}
\label{eq:claH:rk1:12:DCC}
D \cdot C \subset C \;;
\end{equation}
because of \,$1 \in C$\, it follows that
\begin{equation}
\label{eq:claH:rk1:12:DC}
D \subset C 
\end{equation}
holds.

Next we note that \eqref{eq:claH:rk1:12:R-vvH} is a member of \,$(\liem_{2\lambda_2} \oplus \liem_{\lambda_1}) \cap \liem' = \liem_{2\alpha}'$\,;
it follows by Equation~\eqref{eq:claH:rk1:12:m2alpha'} that \,$(-3)\,\IM(\g{\Phi(c_1)}{\Phi(c_2)}) \in D$\, and that 
$$ \wt{\Phi}(\,(-3)\,\IM(\g{\Phi(c_1)}{\Phi(c_2)})\,) = \tfrac12\,\sqrt{3}\,J\bigr( \Phi(c_1)\,c_2 - \Phi(c_2)\,c_1 \bigr) \;, $$
hence
\begin{equation}
\label{eq:claH:rk1:12:ccd}
\IM(\g{\Phi(c_1)}{\Phi(c_2)}) \in D 
\end{equation}
and
\begin{equation}
\label{eq:claH:rk1:12:Phicc}
\Phi(c_2)\,c_1 = \Phi(c_1)\,c_2 - 2\,\sqrt{3}\,J\bigr(\,\wt{\Phi}(\IM \g{\Phi(c_1)}{\Phi(c_2)})\,\bigr)
\end{equation}
holds.

A final relation is obtained 
by considering the \,$\liem_{2\lambda_2}$-component of \eqref{eq:claH:rk1:12:R-vvvt}; we see that for any \,$c_1,c_2 \in C$\, and \,$d_1 \in D$\, we have
\begin{equation}
\label{eq:claH:rk1:12:IM-to-D}
\IM\bigr( \; \IM(c_1\,\overline{c_2})\,d_1 - d_1\,\IM(c_1\,\overline{c_2}) + 3\,d_1\,\IM(\g{\Phi(c_2)}{\Phi(c_1)}) \; \bigr) \;\in\; D \; . 
\end{equation}

Now we use these relations to show
that \,$\liem'$\, is indeed of a type of the form \,$(\PP,\vi=\arctan(\tfrac12),(\K,2))$\,. 

First suppose \,$n_\alpha' = 1$\,. Then we have \,$C = \R$\, and therefore \,$D = \{0\}$\, by \eqref{eq:claH:rk1:12:DC}. Hence Equations~\eqref{eq:claH:rk1:12:m'}
and \eqref{eq:claH:rk1:12:malpha'} show that \,$\liem'$\, is spanned
by \,$H$\, and \,$M_{1,1}^{(e_+,e_-)} + \sqrt{\tfrac32}\,\Phi(1)$\,. It follows that \,$\liem'$\, is of type \,$(\PP,\vi=\arctan(\tfrac12),(\R,2))$\,. 

Thus we now suppose \,$n_\alpha' \geq 2$\,. In this situation we have
\,$n_{2\alpha}' \geq 1$\,: Assume to the contrary that \,$n_{2\alpha}' = 0$\, and therefore \,$D = \{0\}$\, holds. Because of \,$n_\alpha' \geq 2$\,
there exists \,$i \in C \cap \Sph(\IM(\HH))$\,.
We would have \,$\IM(\g{\Phi(1)}{\Phi(i)}) = 0$\,
by \eqref{eq:claH:rk1:12:ccd}, and therefore \,$\Phi(i) = \Phi(1)\cdot i$\, by Equation~\eqref{eq:claH:rk1:12:Phicc}. But using the latter relation
we now calculate \,$\g{\Phi(1)}{\Phi(i)} = \g{\Phi(1)}{\Phi(1)\,i} = \|\Phi(1)\|^2 \cdot i = i$\,, hence \,$i \in D$\, by \eqref{eq:claH:rk1:12:ccd}
in contradiction to our assumption.

If \,$n_\alpha' = 2$\, holds, we have \,$C = \R \oplus \R\,i$\, with some \,$i \in \Sph(\IM(\HH))$\, and then \,$D = \R\,i$\, because of
\,$n_{2\alpha}' \geq 1$\, and \eqref{eq:claH:rk1:12:DC}.
We then have \,$\g{J\Phi(1)}{\wt{\Phi}(i)} = -\tfrac{1}{\sqrt{3}}\,i$\, by 
Equation~\eqref{eq:claH:rk1:12:dc-Phi}, and therefore there exists \,$w_0 \in \Sph(\liem_{\lambda_2})$\, with \,$\g{w_0}{\wt{\Phi}(i)} = 0$\, so that
$$ \Phi(1) = -\tfrac{1}{\sqrt{3}}\,J(\wt{\Phi}(i))\,i + \sqrt{\tfrac23}\,w_0 $$
holds. Further, we have by Equation~\eqref{eq:claH:rk1:12:dc-Phi2} (applied with \,$c_1=1$\,, \,$d_1=i$\,)
$$ \Phi(i) = \Phi(1)\,i - \tfrac{2}{\sqrt{3}}\,J(\wt{\Phi}(i)) = -\tfrac{1}{\sqrt{3}}\,J(\wt{\Phi}(i)) + \sqrt{\tfrac23}\,w_0\,i \; . $$
This shows that \,$\liem_\alpha' = \spn_{\R}\{M_{1,1}^{(e_+,e_-)}-\tfrac{1}{\sqrt{2}}J(\wt{\Phi})i+w_0 \;,\; M_{i,1}^{(e_+,e_-)}-\tfrac{1}{\sqrt{2}}J(\wt{\Phi})+w_0\,i\}$\,
is of the form of \,$\liem''_\perp$\, defined in the description of the type \,$(\PP,\vi=\arctan(\tfrac12),(\C,2))$\,. 
\,$\liea' \oplus \liem_{2\alpha}'$\,,
the ortho-complement of \,$\liem_{\alpha}'$\, in \,$\liem'$\,, is of type \,$(\PP,\vi=\arctan(\tfrac12),(\C,1))$\,
because of \,$n_{2\alpha}'=1$\,.
Hence \,$\liem'$\, is of type \,$(\PP,\vi=\arctan(\tfrac12),(\C,2))$\,.

Thus we now suppose \,$n_\alpha' \geq 3$\,. We will then show that the multiplicities of both \,$\alpha$\, and \,$2\alpha$\, are already maximal,
i.e.~that \,$n_\alpha'=4$\, and \,$n_{2\alpha'}=3$\, holds; it will follow therefrom that \,$\liem'$\, is of type \,$(\PP,\vi=\arctan(\tfrac12),(\HH,2))$\,. 

We know by the preceding arguments that \,$n_{2\alpha}' \geq 1$\, holds, therefore there exists \,$i \in \Sph(\IM(\HH))$\, so that \,$\R\,i \subset D$\, holds.
We have \,$D \subset C$\, by \eqref{eq:claH:rk1:12:DCC}, so we can extend \,$i$\, to a canonical basis \,$(i,j,k)$\, of \,$\IM(\HH)$\, so that
\,$\R \oplus \R\,i \oplus \R\,j \subset C$\, holds. By \eqref{eq:claH:rk1:12:DCC} we then also have \,$C \supset i\cdot(\R \oplus \R\,i \oplus \R\,j) 
= \R\,i \oplus \R \oplus \R\,k$\, and therefore \,$C = \HH$\,, hence \,$n_{\alpha}' = 4$\,. 

Assume that \,$n_{2\alpha'} < 3$\, holds. Then there would exist some \,$q \in \Sph(\IM(\HH))$\, with \,$q \perp_{\R} D$\,.
Because of \,$i \in D$\, we would in particular have \,$q \perp_{\R} i$\, and therefore \,$qi=iq$\,, moreover \,$D \subset \R\,i \oplus \R\,qi$\,. 
By applying \eqref{eq:claH:rk1:12:IM-to-D} with \,$c_1 := q \in C$\,, 
\,$c_2 := 1 \in C$\, and \,$d_1 := i \in D$\,, we see that
\begin{equation}
\label{eq:claH:rk1:12:43trick}
\IM \bigr( \; \IM(q\,\overline{1})\,i - i\,\IM(q\,\overline{1}) + 3\,i\,\IM(\g{\Phi(1)}{\Phi(q)})\;\bigr) 
\;=\; 2\,q\,i + \underbrace{3\,\IM\bigr(\,i\,\IM(\g{\Phi(1)}{\Phi(q)})\,\bigr)}_{(*)} 
\end{equation}
is a member of \,$D$\,. We have \,$\g{\Phi(1)}{\Phi(q)} \in D \subset \R\,i \oplus \R\,qi$\, by \eqref{eq:claH:rk1:12:ccd}, hence 
\,$i\,\g{\Phi(1)}{\Phi(q)} \in \R \oplus \R\,q$\, and therefore the term marked \,$(*)$\, in  \eqref{eq:claH:rk1:12:43trick} lies in \,$\R\,q$\,.
Because \eqref{eq:claH:rk1:12:43trick}, as a member of \,$D$\,, is \,$\R$-perpendicular to \,$q$\,, it follows that the term marked \,$(*)$\, in fact vanishes,
and thus we see that \,$q\,i \in D$\, holds. 
Now we apply \eqref{eq:claH:rk1:12:IM-to-D} again, this time with \,$c_1 := i \in C$\,, \,$c_2 := 1 \in C$\,
and \,$d := qi \in D$\,, to see that 
\begin{equation}
\label{eq:claH:rk1:12:43trick2}
\IM \bigr( \; \IM(i\,\overline{1})\,qi - qi\,\IM(i\,\overline{1}) + 3 \, qi \, \IM(\g{\Phi(1)}{\Phi(i)}) \; \bigr)
= 2\,q + 3\, \IM\bigr(\,qi \, \IM(\g{\Phi(1)}{\Phi(i)})\,\bigr) 
\end{equation}
is a member of \,$D$\,. We now have by Equation~\eqref{eq:claH:rk1:12:Phicc}: \,$\Phi(i) = \Phi(1)\,i - 2\,\sqrt{3}\,J\wt{\Phi}(\IM \g{\Phi(1)}{\Phi(i)})$\,
and therefore
\begin{align*}
\g{\Phi(1)}{\Phi(i)}
& \,\,=\,\, \g{\Phi(1)}{\Phi(1)\,i} - 2\,\sqrt{3}\,\bigr\langle\, \Phi(1) \,,\, J\wt{\Phi}(\IM\, \g{\Phi(1)}{\Phi(i)}) \,\bigr\rangle \\
& \,\,=\,\, i + 2\,\sqrt{3}\,\bigr\langle\, J\Phi(1) \,,\, \wt{\Phi}(\IM \g{\Phi(1)}{\Phi(i)}) \,\bigr\rangle \\
& \overset{\eqref{eq:claH:rk1:12:dc-Phi}} = i - 2\,\IM(\g{\Phi(1)}{\Phi(i)}) \; . 
\end{align*}
From this calculation, we first conclude \,$\g{\Phi(1)}{\Phi(i)} \in \IM(\HH)$\, and then \,$\g{\Phi(1)}{\Phi(i)} = \tfrac13\,i$\,. So we see
that the element \eqref{eq:claH:rk1:12:43trick2} of \,$D$\, in fact equals \,$q$\,. Thus we have shown \,$q \in D$\,, in contradiction to
our assumption. Therefore we have \,$n_{2\alpha}' = 3$\, and hence \,$D = \IM(\HH)$\,. 

Now let \,$c \in \HH$\, be given. Then we have for any \,$d \in \HH$\, by Equation~\eqref{eq:claH:rk1:12:dc-Phi}
$$ \g{J(\wt{\Phi}(d))}{\Phi(c)} = -\tfrac{1}{\sqrt{3}}\,\overline{c}\,d \; . $$
If \,$(i,j,k)$\, is any canonical basis of \,$\IM(\HH)$\,, we therefore have
$$ \Phi(c) \;=\;  \underbrace{-\tfrac{1}{\sqrt{3}}\,\bigr(\, J(\wt{\Phi}(i))\,\overline{c}\,i + J(\wt{\Phi}(j))\,\overline{c}\,j + J(\wt{\Phi}(k))\,\overline{c}\,k \,\bigr)}_{(*)}
\;+\; w' $$
with some vector \,$w' \in \liem_{\lambda_2}$\, which is \,$\HH$-perpendicular to \,$\wt{\Phi}(D)$\,. However, both \,$\Phi(c)$\, and the vector marked \,$(*)$\, above
are vectors of length \,$|c|$\,, whence \,$w'=0$\, follows. Hence we have
\begin{align*}
\liem_{\alpha}' & \overset{\eqref{eq:claH:rk1:12:malpha'}}{=} \Menge{M_{c,1}^{(e_+,e_-)} + \sqrt{\tfrac32}\,\Phi(c)}{c \in \HH} \\
& \;=\; \Mengegr{ \; M_{c,1}^{(e_+,e_-)} - \tfrac{1}{\sqrt{2}}\,\bigr(\,J(\wt{\Phi}(i))\,\overline{c}\,i 
+ J(\wt{\Phi}(j))\,\overline{c}\,j + J(\wt{\Phi}(k))\,\overline{c}\,k \,\bigr) \; }{\;c \in \HH\;} \; . 
\end{align*}
This shows that \,$\liem_\alpha'$\, is of the form of \,$\liem''_\perp$\, defined in the description of the type \,$(\PP,\vi=\arctan(\tfrac12),(\HH,2))$\,. 
\,$\liea' \oplus \liem_{2\alpha}'$\,,
the ortho-complement of \,$\liem_{\alpha}'$\, in \,$\liem'$\,, is of type \,$(\PP,\vi=\arctan(\tfrac12),(\HH,1))$\,
because of \,$n_{2\alpha}'=3$\,.
Hence \,$\liem'$\, is of type \,$(\PP,\vi=\arctan(\tfrac12),(\HH,2))$\,.

\bigskip

\textbf{The case \,$\boldsymbol{\vi_0=\tfrac\pi4}$\,.}
This case could be solved by similar calculations as in the other three cases. However, we here employ a different method, using the
classification result due to \textsc{Wolf} already mentioned in the Introduction.

Let \,$M$\, be the connected, complete, totally geodesic submanifold of \,$G_2(V' \oplus V)$\, with \,$V' \in M$\, and \,$T_{V'}M = \liem'$\,. By
Proposition~\ref{P:orbits:pi4} \,$M$\, has the property that  \,$U_1 \cap U_2 = \{0\}$\, holds for every \,$U_1,U_2 \in M$\, with \,$U_1 \neq U_2$\,. 
The totally geodesic submanifolds with this property have been classified in all Grassmannian manifolds by Wolf in \cite{Wolf:1963-spheres} 
and \cite{Wolf:1963-elliptic}, and we now apply his classification to our situation. When reading the theorems cited below, it should be noted
that the condition of isoclinicity occuring in them is for a totally geodesic submanifold \,$M$\, of \,$G_2(V' \oplus V)$\, equivalent to our property
(\,$U_1 \cap U_2 = \{0\}$\, for every \,$U_1,U_2 \in M$\, with \,$U_1 \neq U_2$\,) by \cite{Wolf:1963-elliptic}, Theorem~2$'$.

\,$M$\, is a symmetric space of rank \,$1$\,, and therefore either isometric to a sphere or to a projective space. 

If \,$M$\, is isometric to a sphere, then by \cite{Wolf:1963-spheres}, Theorem~8 the dimension of \,$M$\, is at most \,$5$\, and every such sphere
is contained in a \,$5$-dimensional one. By \cite{Wolf:1963-spheres}, Theorem~7 any two such spheres of the same dimension are congruent to each other. 
The totally geodesic submanifolds corresponding to the type \,$(\Sph^5,\vi=\tfrac\pi4)$\, are isometric to spheres (see the investigation in Section~\ref{Se:Hembed})
of the type of the present case and therefore \,$M$\, is contained in such a sphere, hence \,$\liem'$\, is contained in a Lie triple system of type
\,$(\Sph^5,\vi=\tfrac\pi4)$\,. Therefore \,$\liem'$\, is of one of the types \,$(\PP,\vi=\tfrac\pi4,(\R,1))$\,,
\,$(\PP,\vi=\tfrac\pi4,(\C,1))$\,, \,$(\PP,\vi=\tfrac\pi4,(\Sph^3))$\,, \,$(\PP,\vi=\tfrac\pi4,(\HH,1))$\, and \,$(\Sph^5,\vi=\tfrac\pi4)$\,, 
depending on whether its dimension is \,$1$\,, \,$2$\,, \,$3$\,, \,$4$\, or \,$5$\,. 

Otherwise \,$M$\, is isometric to a projective space. \cite{Wolf:1963-elliptic}, Lemma~7 shows that the Cayley projective plane does not occur. 
Therefore \,$M$\, is isometric to \,$\KP^\ell$\, for some \,$\K \in \{\R,\C,\HH\}$\, and \,$\ell \geq 2$\,. By \cite{Wolf:1963-elliptic}, Theorem~3
we have \,$\ell \leq \tfrac{n}{2}$\,, and two submanifolds of the present type which are isometric to the same projective space are already congruent in \,$G_2(V' \oplus V)$\,. 
The explicit description of the tangent space of
these totally geodesic submanifolds in \cite{Wolf:1963-elliptic}, Proposition~1 shows that the tangent space \,$\liem'$\,
of \,$M$\, is of type \,$(\PP,\vi=\tfrac\pi4,(\K,\ell))$\,.
\strut\hfilll $\Box$

\section{Totally geodesic submanifolds of \,$\boldsymbol{G_2(\HH^{n+2})}$\,}
\label{Se:Hembed}

Of course we want to know the geometry of the totally geodesic submanifolds which correspond to the various types of Lie triple systems of \,$G_2(\HH^{n+2})$\,.
Their local isometry type can be read off their root systems, which have been constructed in the proof of Theorem~\ref{T:claH:claH}. However, to determine
the global isometry type, and the position of the submanifolds in \,$G_2(\HH^{n+2})$\, (for example, by explicitly describing totally geodesic embeddings),
one has to investigate the global structure for each type separately.

The global isometry types of the various totally geodesic submanifolds, which result from the considerations of the present section, are given in the following
table. Herein, we ascribe the type of a Lie triple system also to the corresponding totally geodesic submanifold (or to a corresponding totally geodesic
embedding). For \,$\ell \in \N$\, and \,$r > 0$\, we denote by \,$\Sph^\ell_r$\, the \,$\ell$-dimensional sphere of radius \,$r$\,, and for \,$\vkap > 0$\,
we denote by \,$\RP^\ell_\vkap$\,, \,$\CP^\ell_\vkap$\, and \,$\HP^\ell_\vkap$\, the respective projective spaces, their metric being scaled in such a way that the \emph{minimal}
sectional curvature is \,$\vkap$\,. (\,$\RP^\ell_\vkap$\, is then of constant sectional curvature \,$\vkap$\, and \,$\CP^\ell_\vkap$\, is of constant holomorphic
sectional curvature \,$4\vkap$\,. The totally geodesic submanifolds of \,$\HP^\ell_\vkap$\, are \,$\HP^{\ell'}_\vkap$\,, \,$\CP^{\ell'}_\vkap$\,, \,$\RP^{\ell'}_\vkap$\, and
\,$\Sph^3_{r=1/2\sqrt{\vkap}}$\, with \,$\ell' \leq \ell$\,.)

\begin{center}
\begin{tabular}{|c|c|}
\hline
type & corresponding global isometry type\\
\hline
\,$(\Geo,\vi=t)$\, & \,$\R$\, or \,$\Sph^1$\, \\
\,$(\PP,\vi=0,(\K,\ell))$\, & \,$\KP^\ell_1$\, \\
\,$(\PP,\vi=0,(\Sph^3))$\, & \,$\Sph^3_{r=1/2}$\, \\
\,$(\Sph,\vi=\arctan(\tfrac13),\ell)$\, & \,$\Sph^\ell_{r=\tfrac12\,\sqrt{10}}$\, \\
\,$(\PP,\vi=\arctan(\tfrac12),(\K,\ell))$\, & \,$\KP^{\ell}_{1/5}$\, \\
\,$(\PP,\vi=\arctan(\tfrac12),(\Sph^3))$\, & \,$\Sph^3_{r=2\sqrt{5}}$\, \\
\,$(\PP,\vi=\tfrac\pi4,(\K,\ell))$\, & \,$\KP^{\ell}_{1/2}$\, \\ 
\,$(\PP,\vi=\tfrac\pi4,(\Sph^3))$\, & \,$\Sph^3_{r=1/\sqrt{2}}$\, \\
\,$(\Sph^5,\vi=\tfrac\pi4)$\, & \,$\Sph^5_{r=1/\sqrt{2}}$\, \\
\hline
\,$(\Gtwo,(\K,\ell))$\, & \,$G_2(\K^{\ell+2})$\, \\
\,$(\PtP,(\K_1,\ell_1),(\K_2,\ell_2))$\, & \,$\mathrm{\K_1P}^{\ell_1}_1 \times \mathrm{\K_2P}^{\ell_2}_1$\, \\
\,$(\PtP,(\K,\ell),(\Sph^3))$\, & \,$\mathrm{\K P}^{\ell}_1 \times \Sph^3_{r=1/2}$\, \\
\,$(\PtP,(\Sph^3),(\Sph^3))$\, & \,$\Sph^3_{r=1/2} \times \Sph^3_{r=1/2}$\, \\
\,$(\Sph^1 \times \Sph^5,\ell)$\, &  $(\Sph^1_{r=1/\sqrt{2}} \times \Sph^\ell_{r=1/\sqrt{2}})/\{\pm \id\}$ \\
\,$(\Sptwo)$\, & \,$\Sp(2)$\, \\
\,$(\Qthree)$\, & \,$G_2^+(\R^5)$\, \\
\hline
\end{tabular}
\end{center}

\begin{Remark}
\label{R:Hembed:diam}
It is an interesting observation that the geodesic diameter of certain of the totally geodesic submanifolds of \,$G_2(\HH^{n+2})$\, 
is strictly larger than the geodesic diameter
\,$\tfrac{\pi}{2}\,\sqrt{2}$\, of \,$G_2(\HH^{n+2})$\, itself: The geodesic diameter of the spheres corresponding to the type \,$(\Sph,\vi=\arctan(\tfrac13),\ell)$\,
equals \,$\tfrac{\pi}{2}\,\sqrt{10}$\,, and the geodesic diameter of the projective spaces corresponding to the types \,$(\PP,\vi=\arctan(\tfrac12),\tau)$\, equals
\,$\tfrac{\pi}{2}\,\sqrt{5}$\,. 
\end{Remark}

As mentioned above, the correctness of the local isometry types (in particular, the correctness of the given radii and curvatures) is easily seen by
inspecting the root systems of the various Lie triple systems. To justify the global isometry type, and also to gain more insight into the 
geometry of some of the types of totally geodesic submanifolds, we now study the individual types separately. Because any two Lie triple systems 
(and therefore also any two totally geodesic submanifolds) of the same type are congruent under the isotropy action of \,$G_2(\HH^{n+2})$\,,
it suffices to provide one example of totally geodesic embedding per type (where the ``type'' of a totally geodesic submanifold as used below
is the same as the type defined in Theorem~\ref{T:claH:claH}).

\paragraph{Type \,$\boldsymbol{(\Gtwo,\tau)}$\,.}
The canonical embedding \,$\HH^{\ell+2} \hookrightarrow \HH^{n+2}$\, induces a totally geodesic isometric embedding \,$G_2(\HH^{\ell+2}) \to G_2(\HH^{n+2})$\, 
of type \,$(\Gtwo,(\HH,\ell))$\,. Moreover for \,$\K \in \{\R,\C\}$\,, let \,$\iota: \K^{\ell+2} \hookrightarrow \HH^{n+2}$\, be the canonical embedding,
then \,$G_2(\K^{\ell+2}) \to G_2(\HH^{n+2}),\;U \mapsto \spn_{\HH}(\iota(U))$\, is a totally geodesic isometric embedding of type \,$(\Gtwo,(\K,\ell))$\,.

\paragraph{Types \,$\boldsymbol{(\PP,\vi=0,\tau)}$\, and \,$\boldsymbol{(\PtP,\tau_1,\tau_2)}$\,.} 
If \,$\HH^{n+2} = V_1 \oplus V_2$\, is an orthogonal splitting into symplectic subspaces \,$V_\nu$\,
of dimension \,$\ell_\nu+1$\, (where \,$\ell_\nu \geq 0$\,, \,$\ell_1 + \ell_2 = n+2$\,), then the map 
$$ f: \HP_1(V_1) \times \HP_1(V_2) \to G_2(\HH^{n+2}), \; (p_1, p_2) \;\mapsto\; p_1 \oplus p_2 $$
is an isometric, totally geodesic embedding. It is of type \,$(\PP,\vi=0,(\HH,n))$\, for \,$\ell_2 =0$\,, of type  \,$(\PtP,(\HH,\ell_1),(\HH,\ell_2))$\,
for \,$\ell_1,\ell_2 \geq 1$\,. 
If now \,$\tau_1,\tau_2$\, are \,$\HP$-types with \,$\dim(\tau_\nu) \leq \ell_\nu$\,,
then there exist totally geodesic submanifolds \,$M_\nu$\, of \,$\HP_1(V_\nu)$\,, and the restriction of \,$f$\, to \,$M_1 \times M_2$\, is a
totally geodesic embedding into \,$G_2(\HH^{n+2})$\, of type \,$(\PP,\vi=0,\tau_1)$\, resp.~\,$(\PtP,\tau_1,\tau_2)$\,.

Note that the Lie triple systems of type \,$(\PtP,(\R,1),(\R,1))$\, are the Cartan subalgebras of \,$\liem$\,, and therefore the totally geodesic submanifolds
of \,$G_2(\HH^{n+2})$\, of that type are the maximal tori of \,$G_2(\HH^{n+2})$\,. 
They are therefore isometric to \,$\RP^1_1 \times \RP^1_1 \cong \Sph^1_{r=1/2} \times \Sph^1_{r=1/2}$\,.

Because any two points of \,$G_2(\HH^{n+2})$\, are connected by a minimal geodesic, which runs in a maximal torus, the geodesic diameter of \,$G_2(\HH^{n+2})$\,
equals the geodesic diameter of its maximal tori, i.e.~\,$\pi/\sqrt{2}$\,.

\paragraph{Type \,$\boldsymbol{(\Geo,\vi=t)}$\,.} The totally geodesic submanifolds of this type are, of course, the traces of geodesics
\,$\gamma: \R \to G_2(\HH^{n+2})$\, with \,$\vi(\dot{\gamma}(0)) = t$\,. \,$\gamma$\, runs within a maximal torus of \,$G_2(\HH^{n+2})$\,,
the latter being isometric to \,$\Sph^1_{r=1/2} \times \Sph^1_{r=1/2}$\, by the preceding result.
Thus it follows from the well-known behavior of the geodesics
on the torus that if \,$\tan(t)$\, is irrational, then \,$\gamma$\, is injective, and \,$\gamma(\R)$\, is a non-embedded totally geodesic submanifold
which is isometric to \,$\R$\, and which is dense in the maximal torus to which \,$\dot{\gamma}(0)$\, is tangential. On the other hand, if \,$\tan(t)$\, 
is rational, say \,$\tan(t) = \tfrac{\ell_1}{\ell_2}$\, with \,$\ell_1,\ell_2 \in \N$\, relatively prime (in the case \,$t=0$\, we put \,$\ell_1 := 0$\,,
\,$\ell_2 := 1$\,), then \,$\gamma$\, is periodic with period \,$L := \pi\,\sqrt{\ell_1^2 + \ell_2^2}$\,, \,$\gamma|[0,L)$\, is injective, and
therefore \,$\gamma(\R)$\, is isometric to \,$\Sph^1_{r=\tfrac12\,\sqrt{\ell_1^2 + \ell_2^2}}$\,. 

\paragraph{Type \,$\boldsymbol{(\Sph,\vi=\arctan(\tfrac13),\ell)}$\,.} The totally geodesic submanifolds \,$M$\, of this type are of dimension \,$\ell$\, and
of constant curvature \,$\tfrac25 = \tfrac{1}{r^2}$\, with \,$r := \tfrac12\,\sqrt{10}$\,. Therefore they are isometric either
to \,$\Sph^\ell_r$\, or to \,$\RP^\ell_{1/r^2}$\,. To distinguish between these two cases, we calculate the length of a geodesic tangential to \,$M$\,:
The preceding discussion of the type \,$(\Geo,\vi=t)$\, shows that the the submanifolds of type \,$(\Geo,\vi=\arctan(\tfrac13))$\, are 
circles of radius \,$\tfrac12 \cdot \sqrt{1^2 + 3^2} = r$\,. Therefore \,$M$\, is isometric to \,$\Sph^\ell_{r=\tfrac12\,\sqrt{10}}$\,. 

The spheres of type \,$(\Sph,\vi=\arctan(\tfrac13),2)$\, are related to the submanifolds of type \,$(\PP,\vi=\arctan(\tfrac12),(\R,2))$\,,
see the discussion of this type in Section~\ref{Se:claC} below.

\paragraph{Type \,$\boldsymbol{(\PP,\vi=\arctan(\tfrac12),\tau)}$\,.}
We will describe the \,$\HP^2$\, which is a maximal totally geodesic submanifold in \,$G_2(\HH^7)$\, of type \,$(\PP,\vi=\arctan(\tfrac12),(\HH,2))$\,.
The other types \,$(\PP,\vi=\arctan(\tfrac12),\tau)$\, correspond to totally geodesic submanifolds of that \,$\HP^2$\, of \,$\HP$-type \,$\tau$\,;
for the types \,$\tau=(\C,2)$\, and \,$\tau=(\R,2)$\, also see Section~\ref{Se:claC}.

Let \,$W$\, be a complex-6-dimensional unitary space, and \,$\tau: W \to W$\, be an anti-unitary transformation
(i.e.~\,$\tau$\, is anti-linear and orthogonal with respect to the real inner product on \,$W$\,) with \,$\tau^2 = -\id_W$\,. Via \,$\tau$\,,
\,$W$\, becomes a symplectic space of quaternionic dimension \,$3$\,;  we also have the corresponding symplectic group
$$ \Sp(W,\tau) := \Menge{B \in \SU(W)}{B\circ \tau = \tau \circ B} \cong \Sp(3) \; . $$

Let us now consider the three-fold alternating product \,$\bigwedge^3 W$\, of the complex linear space \,$W$\,; this is a complex-$20$-dimensional unitary space.
Note that every endomorphism \,$f: W \to W$\, induces an endomorphism \,$f^{(3)}: \bigwedge^3 W \to \bigwedge^3 W$\, characterized by
\,$f^{(3)}(w_1 \wedge w_2 \wedge w_3) = f(w_1) \wedge f(w_2) \wedge f(w_3)$\, for all \,$w_1,w_2,w_3 \in W$\,. \,$\tau^{(3)}$\, is an anti-unitary transformation
on \,$\bigwedge^3 W$\, with \,$(\tau^{(3)})^2 = -\id_{\bigwedge^3 W}$\,, so \,$\bigwedge^3 W$\, becomes a quaternionic-10-dimensional symplectic space
via \,$\tau^{(3)}$\,. \,$\Sp(W,\tau)$\, acts on \,$\bigwedge^3 W$\, by symplectic transformations via \,$B \mapsto B^{(3)}$\,. 

Now let \,$(b_1,b_2,b_3)$\, be any symplectic basis of the symplectic space \,$W$\, and put
$$ \omega := b_1 \wedge \tau(b_1) \,+\, b_2 \wedge \tau(b_2) \,+\, b_3 \wedge \tau(b_3)  \;\in\; {\textstyle \bigwedge^2 W} \; , $$
\,$\omega$\, is non-zero and does not depend on the choice of the symplectic basis \,$(b_1,b_2,b_3)$\,. Because of the latter property, we have
\,$B^{(2)}(\omega) = \omega$\, for every \,$B \in \Sp(W,\tau)$\,, and therefore the symplectic, quaternionic-$7$-dimensional subspace \,$V := (W \wedge \omega)^\perp$\,
of \,$\bigwedge^3 W$\, is invariant under the action of \,$\Sp(W,\tau)$\,. 

We will construct the totally geodesic \,$\HP^2$\, of type \,$(\PP,\arctan(\tfrac12),(\HH,2))$\, in the quaternionic Grassmannian \,$G_2(V) \cong G_2(\HH^7)$\,. 
For this, note that the action of \,$\Sp(W,\tau) \cong \Sp(3)$\, on \,$V$\, induces an action on \,$G_2(V)$\, in the obvious way, 
and we will find an orbit \,$M$\, of this action which is totally geodesic and isomorphic to \,$\HP^2 \cong \Sp(3)/(\Sp(2) \times \Sp(1))$\,.

Fix an orthogonal splitting \,$W = W^2 \operp W^1$\, of \,$W$\, into symplectic subspaces \,$W^2$\, and \,$W^1$\, of quaternionic dimension \,$2$\, resp.~\,$1$\,. 
Let \,$(b_1,b_2)$\, resp.~\,$(b_3)$\, be any symplectic basis of \,$W^2$\, resp.~\,$W^1$\,, and put
$$ \eta := b_1 \wedge \tau(b_1) \,+\, b_2 \wedge \tau(b_2) \,-\, b_3 \wedge \tau(b_3)  \;\in\; {\textstyle \bigwedge^2 V} \; ; $$
\,$\eta$\, is non-zero and does not depend on the choice of the bases. We have \,$U := W^2 \wedge \eta \in G_2(V)$\,.

The isotropy group \,$K$\, of the \,$\Sp(W,\tau)$-action on \,$G_2(V)$\, at the point \,$U$\, equals \,$\Sp(W,\tau)_{2,1} := \Menge{B \in \Sp(W,\tau)}{B(W^2)=W^2}
\cong \Sp(2) \times \Sp(1)$\,. Therefore the orbit \,$M \subset G_2(V)$\, of the said action through \,$U$\, is isomorphic to \,$\HP^2$\,. Moreover,
\,$M$\, is a totally geodesic submanifold of \,$G_2(V)$\, of type \,$(\PP,\vi=\arctan(\tfrac12),(\HH,2))$\,.

{\footnotesize \emph{Proof of the statements on \,$K$\, and \,$M$\,.}
For every \,$B \in \Sp(W,\tau)_{2,1}$\, we have \,$B^{(2)}(\eta) = \eta$\, and therefore \,$\Sp(W,\tau)_{2,1} \subset K$\, holds. To show the converse inclusion,
and also that the orbit \,$M \subset G_2(V)$\, is totally geodesic, we first work on Lie algebra niveau. The Lie algebra of \,$\Sp(W,\tau)_{2,1}$\, is
\,$\liesp(W,\tau)_{2,1} := \Menge{X \in \liesp(W,\tau)}{X(W^2) \subset W^2}$\,. Let any \,$X \in \liesp(W,\tau)$\, which is orthogonal to
\,$\liesp(W,\tau)_{2,1}$\, with respect to the Killing form of \,$\liesp(W,\tau)$\, be given, then we have \,$X(W^2) \subset W^1$\, and \,$X(W^1) \subset W^2$\,. 
Let \,$\Phi: \Sp(W,\tau) \to \Sp(V,\tau^{(3)}),\;B \mapsto B^{(3)}|V$\, be the canonical embedding, and \,$\Phi_L: \liesp(W,\tau) \to \liesp(V,\tau^{(3)})$\, 
its linearization. Then one can calculate that \,$\Phi_L(X)U \subset U^{\perp,V}$\, and \,$\Phi_L(X)U^{\perp,V} \subset U$\, holds by using the explicit presentation
of \,$\Phi_L$\,:
$$ \forall X \in \liesp(W,\tau),\;w_1,w_2,w_3 \in W \; : \; (\Phi_L(X))(w_1 \wedge w_2 \wedge w_3) 
= Xw_1 \wedge w_2 \wedge w_3 + w_1 \wedge Xw_2 \wedge w_3 + w_1 \wedge w_2 \wedge Xw_3 $$
and the symplectic bases
$$ \tfrac{1}{\sqrt{2}}\bigr( b_1 \wedge b_2 \wedge \tau\,b_2 - b_1 \wedge b_3 \wedge \tau\,b_3 \bigr) \;,\;  
\tfrac{1}{\sqrt{2}}\bigr( b_2 \wedge b_1 \wedge \tau\,b_1 - b_2 \wedge b_3 \wedge \tau\,b_3 \bigr) $$
of \,$U$\, resp.~
$$ \tfrac{1}{\sqrt{2}}\bigr( b_3 \wedge b_1 \wedge \tau\,b_1 - b_3 \wedge b_2 \wedge \tau\,b_2 \bigr) \;,\;  
b_1 \wedge b_2 \wedge b_3 \;,\; \tau\,b_1 \wedge b_2 \wedge b_3 \;,\; b_1 \wedge \tau\,b_2 \wedge b_3 \;,\; b_1 \wedge b_2 \wedge \tau\,b_3 $$
of \,$U^{\perp,V}$\,. This shows that \,$\Phi_L$\, maps the ortho-complement \,$\liem$\, of \,$\liesp(W,\tau)_{2,1}$\, in \,$\liesp(W,\tau)$\, 
(with respect to the Killing form)
into the ortho-complement \,$\liep$\, of \,$\liesp(U,\tau^{(3)}) \oplus \liesp(U^{\perp,V},\tau^{(3)})$\, in \,$\liesp(V,\tau^{(3)})$\, (with respect to the Killing form).

The decomposition \,$\liesp(V,\tau^{(3)}) = (\liesp(U,\tau^{(3)}) \oplus \liesp(U^{\perp,V},\tau^{(3)})) \,\oplus\, \liep$\, is the Cartan decomposition
of \,$\liesp(V,\tau^{(3)})$\, induced by the symmetric structure of \,$G_2(V)$\,; therefore the fact \,$\Phi_L(\liem) \subset \liep$\, has several consequences:
First, the Lie algebra \,$\liek$\, of the isotropy group \,$K$\, is contained in \,$\liesp(W,\tau)_{2,1}$\,;
because we have already seen \,$\Sp(W,\tau)_{2,1} \subset K$\,, we in fact have \,$\liek = \liesp(W,\tau)_{2,1}$\,. Thus the neutral component of \,$K$\,
equals \,$\Sp(W,\tau)_{2,1}$\,. Therefore, second, the decomposition \,$\liesp(W,\tau) = \liek \oplus \liem$\, is a Cartan decomposition of \,$\liesp(W,\tau)$\,,
corresponding to the symmetric structure of \,$\HP^2$\,. Therefore the \,$\Sp(W,\tau)$-orbit \,$M$\, is a locally symmetric space, locally isometric to \,$\HP^2$\,. Third, 
\,$\Phi_L(\liem) \subset \liep$\, shows that \,$M$\, is a totally geodesic submanifold of the symmetric space \,$G_2(V)$\,. Fourth, because \,$M$\, is therefore
a globally symmetric space which is locally isometric to \,$\HP^2$\,, \,$M$\, is in fact globally isometric to \,$\HP^2$\, (there do not exist any non-trivial
symmetric covering maps below \,$\HP^2$\,). This finally shows the isotropy group \,$K$\, to be connected, and therefore to be equal to \,$\Sp(W,\tau)_{2,1}$\,.
\strut\hfill $\Box$

}

\paragraph{Type \,$\boldsymbol{(\PP,\vi=\tfrac\pi4,\tau)}$\,.}
Let \,$\tau$\, be a \,$\HP$-type with \,$\dim(\tau) \leq \tfrac{n}{2}$\, and 
let \,$f: M \times M \to G_2(\HH^{n+2})$\, be the totally geodesic embedding of type \,$(\PtP,\tau,\tau)$\, described above (so \,$M$\, is either a \,$\KP^\ell_1$\,
or \,$\Sph^3_{r=1/2}$\, according to \,$\tau$\,). Then the ``diagonal map'' \,$M \to G_2(\HH^{n+2}),\;p \mapsto f(p,p)$\, is a totally geodesic, isometric
embedding of type \,$(\PP,\vi=\tfrac\pi4,\tau)$\,.

\paragraph{Type \,$\boldsymbol{(\Sph^5,\vi=\tfrac\pi4)}$\,.}
The totally geodesic  submanifolds of this type are \,$5$-dimensional, have sectional curvature \,$2=\tfrac{1}{r^2}$\, with \,$r := \tfrac{1}{\sqrt2}$\,, and the
geodesics starting tangential to it (they are of type \,$(\Geo,t=\tfrac\pi4)$\,) have length \,$\pi \cdot \sqrt{1^2+ 1^2} = 2\pi r$\,. It follows
by the analogous argument as for the type \,$(\PP,\vi=\arctan(\tfrac13),\ell)$\, that the submanifolds of type \,$(\Sph^5,\vi=\tfrac\pi4)$\, are 
isometric to \,$\Sph^5_{r=1/\sqrt{2}}$\,.

\paragraph{Type \,$\boldsymbol{(\Sp_2)}$\,.}
The totally geodesic submanifolds of type \,$(\Sp_2)$\, and \,$(\Sph^1 \times \Sph^5,5)$\, are reflective and complementary to each other in \,$G_2(\HH^4)$\,,
see the classification of reflective submanifolds in \cite{Leung:reflective-1979}. There it is also stated that 
the submanifolds of type \,$(\Sp_2)$\, are globally isometric to \,$\Sp(2)$\,. 

We now give another geometric construction of the totally geodesic submanifolds of type \,$(\Sp_2)$\, in \,$G_2(\HH^4)$\,; this construction will also show
again that these submanifolds are globally isometric to \,$\Sp(2)$\,.

For this we need some important
geometric concepts introduced by \textsc{Chen} and \textsc{Nagano}
(see \cite{Chen/Nagano:totges2-1978} and \cite{Chen:1987}), which are applicable to any Riemannian symmetric space \,$M$\, of compact type:
Let \,$p \in M$\, and \,$s_p: M \to M$\, be the geodesic symmetry of \,$M$\, at \,$p$\,, then
the connected components \,$\neq \{p\}$\, of the fixed point set of \,$s_p$\, are called \emph{polars} of \,$M$\, (with respect to \,$p$\,), they
are totally geodesic submanifolds of \,$M$\,.
A \emph{pole} of \,$M$\, is a polar which is a singleton. 
For \,$p_1,p_2 \in M$\,, a point \,$q\in M$\, is called a \emph{center point} between \,$p_1$\, and \,$p_2$\,, if there exists a geodesic
joining \,$p_1$\, to \,$p_2$\, so that \,$q$\, is the middle point on that geodesic.
If \,$p_2$\, is a pole with respect to \,$p_1$\,, we call the set
\,$C(p_1,p_2)$\, of center points between \,$p_1$\, and \,$p_2$\, the \emph{centrosome} of \,$p_1$\, and \,$p_2$\,. It is easy to see that
it is invariant under the action of the isotropy group at \,$p_1$\, (or \,$p_2$\,) and that its connected
components are totally geodesic submanifolds of \,$M$\, (see \cite{Chen:1987}, Proposition~5.1). 

Now let \,$(e_1,\dotsc,e_4)$\, be a symplectic basis of \,$\HH^4$\,, \,$U:= \spn_{\HH}\{e_1,e_2\}
\in G_2(\HH^4)$\,, and \,$K := \Menge{B \in \Sp(4)}{B(U)=U} \cong \Sp(2) \times \Sp(2)$\, be the isotropy group of \,$G_2(\HH^4)$\, at \,$U$\,.
Then \,$U^\perp = \spn_{\HH}\{e_3,e_4\}$\, is a pole to \,$U$\, (in fact, the only pole), 
\,$U' := \spn_{\HH}\{e_1+e_3,e_2+e_4\} = \Menge{(c_1,\dotsc,c_4)\in \HH^4}{c_1=c_3,c_2=c_4}$\,
is a center point between \,$U$\, and \,$U^\perp$\,, and the \,$K$-orbit \,$N$\, through \,$U'$\, is a connected component of the centrosome \,$C(U,U^\perp)$\,,
and therefore a totally geodesic submanifold of \,$G_2(\HH^4)$\,. By explicitly calculating a tangent space of \,$N$\,, one checks that \,$N$\, is of type \,$(\Sp_2)$\,. 
The isotropy group of the action of \,$K$\, at \,$U'$\, equals the diagonal \,$\Delta(K) := \Menge{(B,B)}{B \in \Sp(2)}$\,
of \,$K \cong \Sp(2) \times \Sp(2)$\,, and therefore \,$N$\, is isomorphic as homogeneous space to \,$K/\Delta(K) \cong \Sp(2)$\,.
This confirms that the totally geodesic submanifolds of type \,$(\Sp_2)$\, are globally isometric to \,$\Sp(2)$\,. 

\paragraph{Type \,$\boldsymbol{(\Sph^1 \times \Sph^5,\ell)}$\,.}
Let \,$\liem'$\, be a Lie triple system in \,$\liem$\, of type \,$(\Sph^1 \times \Sph^5,5)$\,. By inspection of its root system, we see that the
totally geodesic submanifold corresponding to \,$\liem'$\, is locally isometric to \,$\Sph^5_{r=1/\sqrt{2}} \times \R$\,. 
There exists a unit vector \,$v_0 \in \liem'$\, so that
\,$\liem'' := (\R\,v_0)^{\perp,\liem'}$\, is a Lie triple system of type \,$(\Sph^5,\vi=\tfrac\pi4)$\,; \,$\liem''$\, corresponds to a totally geodesic submanifold 
of isometry type \,$\Sph^5_{r=1/\sqrt{2}}$\, by the preceding consideration of that type.
Moreover, if \,$v_1 \in \liem''$\, is any unit vector, then \,$\R\,v_0 \oplus \R\,v_1$\, is a Cartan subalgebra of \,$\liem$\,; the totally geodesic submanifold
corresponding to it is therefore a maximal torus, and hence isometric to \,$\Sph^1_{r=1/2} \times \Sph^1_{r=1/2}$\, by the corresponding investigation above.
It follows from these observations that the totally geodesic submanifold corresponding to \,$\liem'$\, is isometric to \,$(\Sph^5_{r=1/\sqrt{2}} \times \Sph^1_{r=1/\sqrt{2}})
/ \{\pm \id\}$\,. 

\paragraph{Type \,$\boldsymbol{(Q_3)}$\,.}
Let \,$M$\, be a totally geodesic submanifold of \,$G_2(\HH^4)$\, of type \,$(\Gtwo,(\C,2))$\,, then \,$M$\, is as Riemannian symmetric space isomorphic
to \,$G_2(\C^4)$\,, as we saw in the treatment of the type \,$(\Gtwo,\tau)$\,. But \,$G_2(\C^4)$\, is as Riemannian symmetric space isomorphic to
the Grassmannian of \emph{oriented} planes \,$G_2^+(\R^6)$\,, as can be seen for example by comparing the Dynkin diagrams of the two spaces,
and noting that both are simply connected (the Lie group isomorphism underlying this isomorphy is \,$\SU(4) \cong \Spin(6)$\,). 
The canonical embedding \,$\R^5 \hookrightarrow \R^6$\,
therefore induces a totally geodesic embedding \,$G_2^+(\R^5) \to M$\,; seen as an embedding into \,$G_2(\HH^4)$\, it is of type \,$(Q_3)$\,.

\section{Totally geodesic submanifolds in the complex 2-Grassmannian}
\label{Se:claC}

Because \,$G_2(\C^{n+2})$\, is a totally geodesic submanifold of \,$G_2(\HH^{n+2})$\, (of type \,$(\Gtwo,(\C,n))$\,) we can now easily determine
the Lie triple systems resp.~the totally geodesic submanifolds of \,$G_2(\C^{n+2})$\,. All we need to do is to see which of the congruence classes
of Lie triple systems in \,$G_2(\HH^{n+2})$\, have members which are contained in a Lie triple system of type \,$(\Gtwo,(\C,n))$\,. In this way
we obtain the result of the following theorem. Here we call a \,$\HP$-type \,$\tau$\, a \emph{\,$\CP$-type}, if it is either of the form
\,$\tau=(\C,\ell)$\, or of the form \,$\tau=(\R,\ell)$\,. 

\,$G_2(\C^{n+2})$\, carries both a \,$\SU(n+2)$-invariant K\"ahler structure \,$J$\, and a \,$\SU(n+2)$-invariant quaternionic K\"ahler structure \,$\frakJ$\,. 
For a description
of these structures, see \cite{Berndt:grassmann-1997}. In the following theorem we also describe the position of the totally geodesic submanifolds
of \,$G_2(\C^{n+2})$\, with respect to these structures.

\begin{Theorem}
\label{T:claC:claC}
Let \,$\liem_1$\, be a Lie triple system of \,$\liem$\, of type \,$(\Gtwo,(\C,n))$\,, i.e.~\,$\liem_1$\, corresponds to \,$G_2(\C^{n+2})$\,. 
Let \,$\{0\} \neq \liem'$\, be an \,$\R$-linear subspace of \,$\liem_1$\,. Then \,$\liem'$\, 
is a Lie triple system of \,$\liem_1$\, if and only if \,$\liem'$\, is of one of the following types:

\begin{itemize}
\item \,$(\Geo,\vi=t)$\,, where \,$t \in [0,\tfrac{\pi}{4}]$\,. 
\item \,$(\PP,\vi=0,\tau)$\,, where \,$\tau$\, is the name of a \,$\CP$-type with \,$\dim(\tau) \leq n$\,.
\item \,$(\Sph,\vi=\arctan(\tfrac13),2)$\,
\item \,$(\PP,\vi=\arctan(\tfrac12),\tau)$\, where \,$\tau$\, is the name of a \,$\CP$-type with \,$\dim(\tau) \leq 2$\, and \,$n \geq \dim(\tau)+w(\tau)$\,
\item \,$(\PP,\vi=\tfrac\pi4,\tau)$\,, where \,$\tau$\, is the name of a \,$\HP$-type (sic!) with \,$\dim(\tau) \leq \tfrac{n}{2}$\,
\item \,$(\Gtwo,\tau)$\,, where \,$\tau$\, is the name of a \,$\CP$-type with \,$\dim(\tau) \leq n$\,
\item \,$(\PtP,\tau_1,\tau_2)$\,, where \,$\tau_1$\, and \,$\tau_2$\, are names of \,$\CP$-types with \,$\dim(\tau_1) + \dim(\tau_2) \leq n$\,
\item \,$(\Sph^1\times\Sph^5,\ell)$\,, where \,$\ell \leq 3$\,
\item \,$(Q_3)$\,
\end{itemize}

In the following table, we give for each type of Lie triple system in \,$G_2(\C^{n+2})$\, the (global) isometry type of the corresponding totally geodesic
submanifolds, their position with respect to the complex structure \,$J$\, (complex, totally real or neither) and with respect to the quaternionic
structure \,$\frakJ$\, (quaternionic, totally complex, totally real, or neither), and state if the Lie triple systems are maximal.
\begin{center}
\begin{longtable}{|c|c|c|c|c|}
\hline
type & isometry type & \,$J$-position & \,$\frakJ$-position & maximal \\
\hline
\endhead
\hline
\endfoot
\,$(\Geo,\vi=t)$\, & \,$\R$\, or \,$\Sph^1$\, & totally real & totally real & no \\
\,$(\PP,\vi=0,(\R,\ell))$\, & \,$\RP^\ell_1$\, & totally real & totally real & no \\
\,$(\PP,\vi=0,(\C,\ell))$\, & \,$\CP^\ell_1$\, & complex & totally complex & for \,$\ell=n$\, \\
\,$(\Sph,\vi=\arctan(\tfrac13),2)$\, & \,$\Sph^2_{r=\sqrt{10}/2}$\, & neither & neither & no \\
\,$(\PP,\vi=\arctan(\tfrac12),(\R,\ell))$\, & \,$\RP^{\ell}_{1/5}$\, & totally real & totally real & no \\
\,$(\PP,\vi=\arctan(\tfrac12),(\C,\ell))$\, & \,$\CP^{\ell}_{1/5}$\, & neither & neither & if \,$n=4$\,: for \,$\ell=2$\, \\
\,$(\PP,\vi=\tfrac\pi4,(\R,\ell))$\, & \,$\RP^\ell_{1/2}$\, & totally real & totally real & no \\
\,$(\PP,\vi=\tfrac\pi4,(\C,\ell))$\, & \,$\CP^\ell_{1/2}$\, & totally real & totally complex & no \\
\,$(\PP,\vi=\tfrac\pi4,(\Sph^3))$\, & \,$\Sph^3_{r=1/\sqrt{2}}$\, & totally real & neither & no \\
\,$(\PP,\vi=\tfrac\pi4,(\HH,\ell))$\, & \,$\HP^\ell_{1/2}$\, & totally real & quaternionic & for \,$2\ell=n$\, \\
\hline
\,$(\Gtwo,(\R,\ell))$\, & \,$G_2(\R^{\ell+2})$\, & totally real & totally complex & for \,$\ell = n$\, \\
\,$(\Gtwo,(\C,\ell))$\, & \,$G_2(\C^{\ell+2})$\, & complex & quaternionic & for \,$\ell = n-1$\, \\
\,$(\PtP,(\R,\ell_1),(\R,\ell_2))$\, & \,$\RP^{\ell_1}_1 \times \RP^{\ell_2}_1$\, & totally real & totally real & no \\
\,$(\PtP,(\R,\ell_1),(\C,\ell_2))$\, & \,$\RP^{\ell_1}_1 \times \CP^{\ell_2}_1$\, & neither & neither & no \\
\,$(\PtP,(\C,\ell_1),(\C,\ell_2))$\, & \,$\CP^{\ell_1}_1 \times \CP^{\ell_2}_1$\, & complex & totally complex & for \,$\ell_1+\ell_2=n$\, \\
\,$(\Sph^1\times\Sph^5,\ell)$\,& {\footnotesize \,$(\Sph^1_{r=1/\sqrt{2}} \times \Sph^\ell_{r=1/\sqrt{2}})/\{\pm \id\}$\,} & totally real & neither & if \,$n=2$\,: for \,$\ell=3$\, \\
\,$(Q_3)$\, & \,$G_2^+(\R^5)$\, & complex & neither & if \,$n=2$\, \\
\hline
\end{longtable}
\end{center}
\end{Theorem}

\beweis
The only thing that is not immediately obvious from the classification of totally geodesic submanifolds in \,$G_2(\HH^{n+2})$\, given in Theorem~\ref{T:claH:claH}
is the fact that \,$\liem_1$\, contains Lie triple systems of type \,$(\PP,\vi=\tfrac\pi4,(\HH,\ell))$\, for every \,$\ell \leq \tfrac{n}{2}$\,. To see this,
one should use the alternative description of that type given in Remark~\ref{R:claH:pi4alternative}.
\beweisende

\begin{Remarks}
\label{Rs:claC:claC}
\begin{enumerate}
\item
The type \,$(\PP,\vi=\arctan(\tfrac12),(\C,2))$\, (isomorphic to \,$\CP^2_{1/5}$\,, maximal in \,$G_2(\C^6)$\,) is not found in 
the ``classification'' of totally geodesic submanifolds of \,$G_2(\C^{n+2})$\, in Table~VIII of \cite{Chen/Nagano:totges2-1978}. It has never been described
before, as far as I know. 

Moreover, the existence of the type \,$(\PP,\vi=\arctan(\tfrac12),2)$\,, although not maximal in any \,$G_2(\C^{n+2})$\, (but rather
in \,$G_2(\R^5)$\,), can not be deduced from that table; it corresponds to the type \,$(\mathrm{A})$\, of my classification of totally
geodesic submanifolds in the complex quadric \,$Q^n \cong G_2^+(\R^{n+2})$\, in \cite{Klein:2007-claQ}.
\item
We again encounter the phenomenon (compare Remark~\ref{R:Hembed:diam}) that the geodesic diameter of certain  totally geodesic submanifolds
is strictly larger than the geodesic diameter of the ambient space \,$G_2(\C^{n+2})$\,, which equals \,$\tfrac{\pi}{2}\,\sqrt{2}$\,. 
As in the case of \,$G_2(\HH^{n+2})$\,, this is true of the totally geodesic submanifolds of type \,$(\Sph,\vi=\arctan(\tfrac13),2)$\,
(diameter \,$\tfrac{\pi}{2}\,\sqrt{10}$\,), and of those of the types \,$(\PP,\vi=\arctan(\tfrac12),\tau)$\, (diameter
\,$\tfrac{\pi}{2}\,\sqrt{5}$\,).
\end{enumerate}
\end{Remarks}

Of course, one would also like to understand totally geodesic, isometric embeddings into \,$G_2(\C^{n+2})$\, for the various types of its Lie triple systems.
For the most part they can be easily obtained by appropriately restricting embeddings into \,$G_2(\HH^{n+2})$\, described in Section~\ref{Se:Hembed}.
So we here discuss only the ``unexpected'' types \,$(\Sph,\vi=\arctan(\tfrac13),2)$\, and \,$(\PP,\vi=\arctan(\tfrac12),\tau)$\,.

\paragraph{Type \,$\boldsymbol{(\PP,\vi=\arctan(\tfrac12),\tau)}$\,.}
In particular we wish to construct a totally geodesic submanifold of type \,$(\PP,\vi=\arctan(\tfrac12),(\C,2))$\,, which is maximal in \,$G_2(\C^6)$\,.
For this, we continue the discussion of the type \,$(\PP,\vi=\arctan(\tfrac12),(\HH,2))$\, in Section~\ref{Se:Hembed}, and use the objects introduced there.

We fix a complex form \,$W^{\C}$\, of \,$W$\, (i.e.~a complex-3-dimensional, totally complex subspace of \,$(W,\tau)$\,).
Under the action of \,$\SU(W^{\C}) \cong \SU(3)$\, on \,$\bigwedge^3 {W}$\, (where \,$B \in \SU(W^{\C})$\,
acts via \,$(B')^{(3)}$\, with the unique endomorphism \,$B' \in \Sp(W)$\, with \,$B'|W^{\C} = B$\,) the complex-$9$-dimensional space \,$L$\, spanned by
$$ \Mengegr{\; \tau\,b_\nu \wedge b_2 \wedge b_3 \;,\; b_1 \wedge \tau\,b_\nu \wedge b_3 \;,\; b_1 \wedge b_2 \wedge \tau\,b_\nu \;}{\;\nu\in\{1,2,3\}\;} $$
is invariant (where again \,$(b_1,b_2,b_3)$\, is a symplectic basis of \,$W$\,); 
therefore also the complex-$6$-dimensional space \,$V^{\C} := V \cap L$\, is invariant under that action.

Consider the complex Grassmannian \,$G_2(V^{\C}) \cong G_2(\C^6)$\,. We have \,$U^{\C} := U \cap V^{\C} \in G_2(V^{\C})$\,. It is easily seen that the
isotropy group of the action of \,$\SU(W^{\C})$\, on \,$G_2(V^{\C})$\, at the point \,$U^{\C}$\, equals \,$\mathrm{S}(\Ug(W^{2}\cap W^{\C}) \times \Ug(W^{1} \cap W^{\C}))$\,,
and therefore the orbit \,$M^{\C}$\, of that action through \,$U^{\C}$\, is a totally geodesic submanifold of \,$G_2(V^{\C})$\, of type
\,$(\PP,\vi=\arctan(\tfrac12),(\C,2))$\,.

Similarly, we can construct a totally geodesic submanifold of type \,$(\PP,\vi=\arctan(\tfrac12),(\R,2))$\,; they are not maximal in \,$G_2(\C^6)$\,,
but they are maximal in a totally geodesic submanifold of type \,$(\Gtwo,(\R,3))$\,, i.e.~in a \,$G_2(\R^5)$\,. 

For this we fix a real form \,$W^{\R}$\, of \,$W^{\C}$\,. Then
\,$W^{\R} \oplus \tau(W^{\R})$\, is a real form of the unitary space \,$W$\,, and therefore \,$\bigwedge^3_{\R} W := \spn_{\R}\Menge{w_1 \wedge w_2 \wedge w_3}{w_1,w_2,w_3
\in W^{\R} \oplus \tau(W^{\R})}$\, is a real form of \,$\bigwedge^3 W$\,. The group \,$\SO(W^{\R}) \cong \SO(3)$\, acts on \,$\bigwedge^3 W$\, in the obvious way,
and \,$\bigwedge^3_{\R} W$\, is invariant under that action; therefore also the real-6-dimensional space \,$V^{\C} \cap \bigwedge^3_{\R}W$\,. 

Moreover, we have \,$B^{(3)}(\zeta) = \zeta$\, for every \,$B \in \SO(W^{\R})$\,, where
$$ \zeta := \tau\,b_1 \wedge b_2 \wedge b_3 + b_1 \wedge \tau\,b_2 \wedge b_3 + b_1 \wedge b_2 \wedge \tau\,b_3 \in V^{\C} \cap {\textstyle \bigwedge^3_{\R} W} \; . $$
(Indeed, the map \,$\alpha: \End(W^{\R}) \to \R,\; B \mapsto \tfrac{1}{\|\zeta\|^2}\,\g{B^{(3)}\zeta}{\zeta}$\, is multi-linear and alternating 
in the columns of \,$B \in \End(W^{\R})$\,,
for \,$B \in \SO(W^{\R})$\, we hence have \,$\alpha(B) = \det(B) \cdot \alpha(\id) = 1$\,, therefrom \,$B^{(3)}(\zeta) = \zeta$\, follows.)
Therefore also the real-$5$-dimensional space \,$V^{\R} := V^{\C} \cap \bigwedge^3_{\R}W \cap (\R\,\zeta)^\perp$\, is invariant under \,$\SO(W^{\R})$\,. 

We now consider the real Grassmannian \,$G_2(V^{\R}) \cong G_2(\R^5)$\,. With analogous arguments as before one sees: The orbit of the \,$\SO(W^{\R})$-action
on \,$G_2(V^{\R})$\, through the point \,$U^{\R} := U^{\C} \cap V^{\R} \in G_2(V^{\R})$\, is a totally geodesic submanifold of type \,$(\PP,\vi=\arctan(\tfrac12),(\R,2))$\,,
isometric to \,$\RP^2$\,. 

\paragraph{Type \,$\boldsymbol{(\Sph,\vi=\arctan(\tfrac13),2)}$\,.} The totally geodesic submanifolds of this type can be constructed in the following way:
\,$G_2(\C^4)$\, is holomorphically isometric to the \emph{oriented}, real 2-Grassmannian \,$G_2^+(\R^6)$\, (the 4-dimensional complex quadric), it therefore
contains totally geodesic submanifolds isometric to \,$G_2^+(\R^5)$\, (of type \,$(\Qthree)$\,). \,$G_2^+(\R^5)$\, is a two-fold covering manifold over
\,$G_2(\R^5)$\, (the covering map being given by ``forgetting'' the orientation of \,$U \in G_2^+(\R^5)$\,). The pre-image \,$M$\, under this covering map 
of a totally geodesic submanifold of \,$G_2(\R^5)$\, of type \,$(\Sph,\vi=\arctan(\tfrac12),(\R,2))$\, (isometric to \,$\RP^2$\,) 
is of course a totally geodesic submanifold of \,$G_2^+(\R^5)$\,; it turns out that \,$M$\, is connected, and therefore isometric to a \,$2$-sphere. 

Seen as a totally geodesic submanifold of \,$G_2^+(\R^5)$\, (which is isomorphic to the complex quadric \,$Q^3$\,) \,$M$\, is of type \,$(\mathrm{A})$\, of my
classification of totally geodesic submanifolds of the complex quadrics in \cite{Klein:2007-claQ}; seen as a totally geodesic submanifold of \,$G_2(\C^6)$\,,
\,$M$\, is of type  \,$(\Sph,\vi=\arctan(\tfrac13),2)$\, of the present classification.

\end{document}